\documentclass[12pt, reqno]{amsart}
\usepackage[utf8]{inputenc}
\usepackage[T1]{fontenc}
\usepackage{lmodern}
\usepackage{microtype}
\usepackage{mathtools}
\usepackage{amssymb}
\usepackage{epsfig}
\usepackage{graphicx}
\usepackage{color}
\usepackage{fullpage}
\definecolor{shadecolor}{gray}{0.875}
\usepackage{amscd}
\usepackage{comment}
\usepackage{tikz-cd}
\usepackage{ifdraft}
\usepackage{thmtools}

\usepackage[style=ieee-alphabetic, url=false, backref=true]{biblatex}
\usepackage{hyperref}
\usepackage[capitalize]{cleveref}

\numberwithin{equation}{section}

\usepackage[all]{xy}

\calclayout
\allowdisplaybreaks[3]

\declaretheorem[title=Proposition, parent=section]{prop}
\declaretheorem[title=Theorem, sibling=prop, refname={Theorem,Theorems}, Refname={Theorem,Theorems}]{theo}
\declaretheorem[title=Corollary, sibling=prop]{coro}

\declaretheorem[title=Lemma, sibling=prop]{lemm}
\declaretheorem[title=Assumption, sibling=prop]{assu}

\declaretheorem[title=Definition, style=definition, sibling=prop]{defi}
\declaretheorem[title=Conjecture, style=definition, sibling=prop]{conj}

\declaretheorem[title=Remark, style=definition, sibling=prop]{rema}
\declaretheorem[title=Example, style=definition, sibling=prop, refname={Example,Examples}, Refname={Example,Examples}]{exam}
\declaretheorem[title=Notation, style=definition, sibling=prop]{nota}
\declaretheorem[title=Construction, style=definition, sibling=prop]{const}

\crefname{equation}{Equation}{Equations}
\Crefname{equation}{Equation}{Equations}

\def\cF{{\mathcal F}}

\def\cI{{\mathcal I}}

\def\cN{{\mathcal N}}
\def\cO{{\mathcal O}}

\def\cV{{\mathcal V}}

\def\rk{{\mathrm{rk}}}

\def\bZ{{\mathbb Z}}

\def\Pic{\mathrm{Pic}}

\def\vol{\mathrm{Vol}}

\def\Mor{\mathrm{Mor}}
\def\Top{\mathrm{Top}}
\def\Nef{\mathrm{Nef}}
\def\Hilb{\mathrm{Hilb}}
\def\Supp{\mathrm{Supp}}

\def\Spec{\mathrm{Spec}}
\def\Sym{\mathrm{Sym}}
\def\Pic{\mathrm{Pic}}

\def\Tr{\mathrm{Tr}}

\newcommand{\hMor}{\widehat{\operatorname{Mor}}_e}
\newcommand{\tMor}{\widetilde{\operatorname{Mor}}_e}

 \author{Ronno Das}
 \address{Ronno Das, Department of Mathematics \\ Stockholm University}
 \email{ronnodas@gmail.com}
 
 \author{Brian Lehmann}
\address{Brian Lehmann, Department of Mathematics \\
Boston College  \\
Chestnut Hill, MA \\
02467 USA}
\email{lehmannb@bc.edu}

\author{Sho Tanimoto}
\address{Sho Tanimoto, Graduate School of Mathematics \\
Nagoya University \\
Furocho Chikusa-ku, Nagoya \\
464-8602, Japan}
\email{sho.tanimoto@math.nagoya-u.ac.jp}

 \author{Philip Tosteson}
 \address{Philip Tosteson, Department of Mathematics \\ University of North Carolina}
 \email{ptoste@unc.edu}
 
 \address{Will Sawin, Princeton University, Fine Hall, 304 Washington Rd, Princeton NJ 08540, USA}
\email{wsawin@princeton.edu}

\address{Mark Shusterman, Faculty of Mathematics and Computer Science, Weizmann
Institute of Science, 234 Herzl Street, Rehovot 76100, Israel.}
\email{mark.shusterman@weizmann.ac.il}

\makeatletter
\let\@wraptoccontribs\wraptoccontribs
\makeatother

\contrib[with an appendix by]{Will Sawin and Mark Shusterman}

\title[Manin's conjecture for quartic del Pezzo surfaces]{Homological stability and Manin's conjecture for rational curves on quartic del Pezzo surfaces}


\begin{document}

\begin{abstract}

We prove a version of Manin's conjecture (over $\mathbb{F}_{q}$ for $q$ large)  and the Cohen--Jones--Segal conjecture (over $\mathbb{C}$) for maps from rational curves to split quartic del Pezzo surfaces.  The proofs share a common method which builds upon prior work of the first and fourth authors.  The main ingredients of this method are (i) the construction of bar complexes formalizing the inclusion-exclusion principle and its point counting estimates, (ii) dimension estimates for spaces of rational curves using conic bundle structures, (iii) estimates of error terms using arguments of Sawin--Shusterman based on Katz's results, and (iv) a certain virtual height zeta function revealing the compatibility of bar complexes and Peyre's constant.
Our argument substantiates the heuristic approach to Manin's conjecture over global function fields given by Batyrev and Ellenberg--Venkatesh in this case.  
\end{abstract}

\maketitle

\tableofcontents

\section{Introduction}

One of the central problems in Diophantine geometry is to understand the distribution of rational points on algebraic varieties over global fields. Around the 1980's, Yuri Manin and his collaborators proposed an ambitious program to predict the asymptotic formula for the counting function of rational points on a smooth Fano variety defined over a number field in terms of its global geometry.  This leads to the formulation of Manin's conjecture which has been refined by various mathematicians over the last three decades. Here is a current formulation of Manin's conjecture: 

\begin{conj}[Manin's conjecture, {\cite{FMT, BM, Peyre, BT, Peyre03, LST18, LS24}}]
\label{conj:Manin}
\hfill
Let $k$ be a number field and let $X$ be a smooth Fano variety defined over $k$.
Fix an adelic metrization on the anticanonical divisor $-K_X$ and consider the corresponding height function $\mathsf H_{-K_X} : X(k) \to \mathbb R_{\geq 0}$ and the associated counting function
\[
N(Q, \mathsf H_{-K_X},T) = \#\{ x \in Q \, | \, \mathsf H_{-K_X}(x) \leq T\},
\]
where $Q \subset X(k)$ is any subset. Suppose that $X(k)$ is not thin in the sense of Serre (\cite{Serre}). Then there exists a thin set $Z \subset X(k)$ such that we have
\[
N(X(k) \setminus Z, \mathsf H_{-K_X},T)  \sim \alpha(-K_X, \Nef_1(X))\beta(X)\tau_{-K_X}(X)T(\log T)^{\rho(X)-1}
\]
as $T \to \infty$ where $\rho(X)$ is the Picard rank of $X$, $\alpha(-K_X, \Nef_1(X))$ is the alpha constant of the nef cone $\Nef_1(X)$, $\beta(X)$ is the size of the Brauer group of $X$ modulo constants, and $\tau_{-K_X}(X)$ is the Tamagawa number as introduced in \cite{Peyre}, \cite{BT}, and \cite{LS24}.
\end{conj}

Our focus is the global function field analogue of \cref{conj:Manin} for trivial families over $\mathbb F_q(t)$. 
When $X$ is a Fano variety defined over a finite field $\mathbb F_q$, Manin's conjecture for $X_{\mathbb{F}_{q}(t)}$ predicts the number of $\mathbb{F}_{q}$-rational curves $s :\mathbb P^1 \to X$ of bounded height, or equivalently, the number of $\mathbb{F}_{q}$-points on a bounded subfamily of the moduli space of rational curves $\mathrm{Mor}(\mathbb P^1, X)$.
In this context, Victor Batyrev developed an influential heuristic argument supporting Manin's conjecture (\cite{Bat88}).  
Moreover Jordan Ellenberg and Akshay Venkatesh (\cite{EV05}) suggested that Batyrev's heuristics might be justifed by combining homological stability of $\mathrm{Mor}(\mathbb P^1, X)_\alpha$ with the Grothendieck--Lefschetz trace formula.  

The approach proposed by Ellenberg and Venkatesh resonates with another conjecture on homological stability of the space of rational curves: the Cohen--Jones--Segal conjecture for Fano varieties over $\mathbb{C}$ (\cite{Segal, CJS00}).  Loosely speaking, this conjecture predicts that the homology groups of the space of pointed algebraic morphisms $\Mor_{*}(\mathbb{P}^{1},X)_{\alpha}$ should ``stabilize'', as the homology class $\alpha$ becomes more positive, to the homology groups of the space of pointed topological morphisms $\Top_{*}(\mathbb{P}^{1},X)_{\alpha}$.  However the general  relationship between the Cohen--Jones--Segal conjecture and Manin's conjecture is, for the most part, only explained heuristically in the literature, and we consider further study into this aspect to be important.

Our main theorems prove a version of Manin's conjecture (over $\mathbb F_q$ for $q$ large) and the Cohen--Jones--Segal conjecture (over $\mathbb C$) for split degree $4$ del Pezzo surfaces.   
The proofs rely on a common method, highlighting the compatibility of the two conjectures and realizing Batyrev--Ellenberg--Venkatesh's heuristics for these surfaces.  We combine algebro-geometric methods which describe the geometry of the moduli spaces of rational curves with topological methods which establish homological stability via the inclusion-exclusion principle.   A key tool in our work is the bar complex  as recently developed by the first and fourth authors in \cite{DT24} in the context of the Cohen--Jones--Segal conjecture.

\subsection{Main results}

Let $S$ be a smooth quartic del Pezzo surface over $\mathbb{F}_{q}$ such that the Picard rank of $S$ is equal to the Picard rank of $S\otimes \overline{\mathbb F}_q$.  We prove Manin's conjecture for $S$ after slightly shrinking the nef cone of curves.


To this end, let $\ell$ be a non-negative, rational, homogeneous, continuous, and piecewise linear function on the nef cone $\mathrm{Nef}_1(S)$ of $S$.
Let $\epsilon > 0$ be a small rational number. We define the rational polyhedral cone $\mathrm{Nef}_1(S)_{\ell, \epsilon}$ to be the closure of the set of non-zero real classes $\alpha \in \mathrm{Nef}_1(S)$ such that $\ell(\alpha)/(-K_S.\alpha) \geq \epsilon$. 

For a positive integer $d$, we let $N^{\mathrm{Nef}_1(S)_{\ell, \epsilon}}(\mathbb P^1, S, -K_S, d)$ denote the number of $\mathbb F_{q}$-morphisms $s: \mathbb{P}^{1} \to S$ of anticanonical degree $\leq d$ such that the class of the image is contained in $\mathrm{Nef}_1(S)_{\ell, \epsilon}$.

\begin{theo} \label{theo:intromaninconj}
There exists a uniform constant $C$ with the following properties:

Fix a sufficiently small positive rational $\epsilon > 0$ 
and let $q$ be a prime power such that $q^{\epsilon} > C$.  Let $S$ be a smooth quartic del Pezzo surface over $\mathbb{F}_{q}$ such that the Picard rank of $S$ is equal to the Picard rank of $S\otimes \overline{\mathbb F}_q$. Then there exists a non-negative, rational, homogeneous, continuous, and piecewise linear function $\ell$ on $\mathrm{Nef}_1(S)$ which takes positive values on a dense open cone $U\subset \mathrm{Nef}_1(S)$ such that
\begin{equation*}
N^{\mathrm{Nef}_1(S)_{\ell, \epsilon}}(\mathbb P^1, S, -K_S, d) \sim_{d \to \infty} (1-q^{-1})^{-1}\alpha(-K_S, \mathrm{Nef}_1(S)_{\ell, \epsilon})\tau_{-K_S}(S)q^{d}d^{5},
\end{equation*}
where  $\alpha(-K_S, \mathrm{Nef}_1(S)_{\ell, \epsilon})$ is the alpha constant of the cone $\mathrm{Nef}_1(S)_{\ell, \epsilon}$ and $\tau_{-K_S}(S)$ is the Tamagawa constant, both introduced in \cite{Peyre}.  
More explicitly, the Tamagawa constant is given as the Euler product
\begin{equation*}
\tau_{-K_{S}}(S) = q^{2} \left(1 - q^{-1} \right)^{-6} \left( \prod_{c \in |\mathbb{P}^{1}|} \left( 1- q^{-|c|} \right)^6  \frac{\# S(\mathbb F_{q^{|c|}})}{q^{2|c|}} \right),
\end{equation*}
where $|\mathbb{P}^{1}|$ is the set of closed points $c \in \mathbb{P}^{1}$ and $|c|$ is the degree of the residue field $\kappa(c)$ over $\mathbb{F}_{q}$.
\end{theo}

\begin{rema}
We prove that $C = 2^{32}$ works in this statement. Also, the function $\ell$ is uniform and independent of the choice of data $q, \epsilon, S$ by naturally identifying the nef cones of any two smooth quartic del Pezzo surfaces defined over algebraically closed fields.
\end{rema}
\begin{rema}
For any integral class $\alpha \in \mathrm{Nef}_1(S)_{\ell, \epsilon}$, we prove the convergence: $$\lim_{m \to \infty}\frac{\#\mathrm{Mor}(\mathbb P^1, S)_{m\alpha}(\mathbb F_q)}{q^{-mK_S.\alpha}} = \tau_{-K_S}(S).$$ 
 We establish this convergence with a uniform error term so that we can take the summation over lattice points in $\mathrm{Nef}_1(S)_{\ell, \epsilon}$ and get the polynomial term.
  Thus one might consider our statement as a version of Manin's conjecture using the all heights approach proposed by Peyre (\cite{Peyre17, Peyre21}) which predicts the asymptotic formula for the counting function of rational curves whose classes are $\epsilon$-away from the boundary of the nef cone. 
\end{rema}

\begin{rema}
Our conic-bundle method for studying rational curves on split degree $4$ del Pezzo surfaces can be used to analyze curves on other split del Pezzo surfaces as well.  This leads to
proofs of a similar version of Manin's conjecture for split del Pezzo surfaces of degree $\geq 5$ as well as lower bounds of the correct magnitude for split del Pezzo surfaces of degree $\leq 3$.
These applications will be explored in a subsequent paper \cite{Tan25}.
\end{rema}

\begin{rema}
Using our method, one can achieve a power saving error term for the counting function $N^{\mathrm{Nef}_1(S)_{\ell, \epsilon}}(\mathbb P^1, S, -K_S, d)$. Indeed there exists a uniform constant $0 < \zeta \leq 1$ independent of $q, \epsilon$ such that one can obtain an error term of order $O(q^{(1-\kappa)d})$ where $0 < \kappa <\min\{ \epsilon - \log_qC, \zeta\epsilon\}$. (We should mention that the polynomial part of the leading term is given by the Ehrhart polynomial of the cone.) Note that our assumption $q^{\epsilon} > C$ implies that $\epsilon - \log_qC >0$, and this saving gets larger as $q$ gets larger, but will not exceed $\zeta\epsilon$.  
\end{rema}

Since the exponential bound on $q$ depends on $\epsilon$, \cref{theo:intromaninconj} does not directly imply Manin's conjecture for the entire nef cone.  However, we are able to establish an approximate result.

\begin{coro} \label{coro:weakmaninconj}
There exists a uniform constant $R$ with the following property.

Let $S$ be a smooth quartic del Pezzo surface over a finite field $\mathbb{F}_{q}$ such that the Picard rank of $S$ is equal to the Picard rank of $S\otimes \overline{\mathbb F}_q$. Then
\begin{equation*}
\limsup_{d \to \infty}\left| \frac{N(\mathbb P^1, S, -K_S, d)}{q^{d+2} d^{5}}  - (1-q^{-1})^{-1}q^{-2}\alpha(-K_S, \Nef_{1}(S))\tau_{-K_S}(S) \right| < \frac{R}{\log(q)}.
\end{equation*}
\end{coro}

Next suppose that $S$ is a degree $4$ del Pezzo surface over $\mathbb{C}$ and that $\alpha$ is a curve class in the interior of the nef cone $\Nef_{1}(S)$ of curves.  \cite{Testa09} shows that there is a unique irreducible component $M_{\alpha}$ parametrizing morphisms $s: \mathbb{P}^{1} \to S$ of homological class $\alpha$.  We let $M_{\alpha,*}$ denote the corresponding moduli space of pointed morphisms and $\Top_{*}(\mathbb{P}^{1},X)_\alpha$ denote the space of pointed continuous maps of class $\alpha$.

Loosely speaking, \cite{CJS00} predicts that if $\alpha$ is in the interior of the nef cone the cohomology groups of $M_{m\alpha,*}$ should stabilize to the cohomology groups of $\Top_{*}(\mathbb{P}^{1},X)_{m\alpha}$ as $m$ goes to infinity.  The following statement gives a more refined statement of homological stability behaves for the spaces $M_{\alpha,*}$.

\begin{theo} \label{theo:introcjs}
Let $S$ be a smooth del Pezzo surface of degree $4$ over $\mathbb{C}$.  There exist a dense open subset cone $U \subset \Nef_{1}(S)$, a rational homogeneous continuous piecewise linear function $\ell: \mathrm{Nef}_1(S) \to \mathbb{R}$ which takes positive values on $U$, and a constant $c \in \mathbb R_{\geq 0}$ such that
for any $\alpha \in U$ the homomorphisms
\begin{equation*}
H_{i}^{\mathrm{sing}}(M_{\alpha,*},\mathbb{Z}) \to H_{i}^{\mathrm{sing}}( \Top_{*}(\mathbb{P}^{1},X)_{\alpha},\mathbb{Z})
\end{equation*}
induced by inclusion are isomorphisms in the range $i \leq \ell(\alpha) -c$.  
\end{theo}

We emphasize that the homotopy type of $\Top_{*}(\mathbb{P}^{1},X)_{\alpha}$ does not depend on $\alpha$ (because up to homotopy $\Top_*(\mathbb P^1, X)$ is acted on by $\pi_2(X)$ which permutes $\{\Top_{*}(\mathbb{P}^{1},X)_{\alpha}\}_{\alpha \in H_2^{\mathrm{sing}}(X, \mathbb Z)}$ transitively), so this statement indeed shows homological stability of the spaces of pointed rational curves.

Our proof of this theorem follows the strategy of \cite{DT24} which confirms the Cohen--Jones--Segal conjecture for quintic del Pezzo surfaces. The new ingredient in this paper is an algebro-geometric argument estimating the dimension of certain loci in the space of rational curves.  This allows us to show that our stratification of the moduli space has the expected behavior in a certain linear codimension range, allowing us to apply the results of \cite{DT24}.

\subsection{A sketch of our proofs}
\label{subsec:sketch}

Here we provide a sketch of our proof of \cref{theo:intromaninconj}.  We will study the moduli space $\Mor(\mathbb{P}^{1},S)$ by fixing a  birational map $\rho: S \to \mathbb{P}^{1} \times \mathbb{P}^{1}$ and studying the postcomposition map $\Mor(\mathbb{P}^{1},S) \to \Mor(\mathbb{P}^1,\mathbb{P}^{1} \times \mathbb{P}^{1})$. 
(In practice we must let the choice of $\rho$ vary as we consider different families of curves.)  Given a rational curve $s': \mathbb{P}^{1} \to \mathbb{P}^{1} \times \mathbb{P}^{1}$, the numerical class of its strict transform $s: \mathbb{P}^{1} \to S$ is determined by its multiplicities along the four blown-up points.  In this way each irreducible component of $\Mor(\mathbb{P}^{1},\mathbb{P}^{1} \times \mathbb{P}^{1})$ admits a stratification into locally closed subsets representing the images of various irreducible components of $\Mor(\mathbb{P}^{1},S)$.

In fact, we do not work with the space of morphisms directly but with a closely related affine bundle $E$ coming from the universal torsor description of maps to $\mathbb{P}^{1} \times \mathbb{P}^{1}$.  We can stratify $E$ into locally closed subsets which have fixed order of vanishing against the four blown-up points and the ``degeneracy loci'' in $E$.  For each irreducible component $M_{\alpha}$ of $\Mor(\mathbb{P}^{1},S)$ there is a $\mathbb{G}_{m}^{2}$-torsor $\widetilde{M}_{\alpha} \to M_{\alpha}$ and an inclusion $\widetilde{M}_{\alpha} \subset E$ realizing $\widetilde{M}_{\alpha}$ as one of the strata.  Unfortunately the pieces of this stratification are complicated and it is not clear how to directly run an inclusion-exclusion argument to analyze $\widetilde{M}_{\alpha}$.

Instead, we consider the subloci of $E$ where the four blown-up points and degeneracy loci vanish along fixed subschemes of $\mathbb{P}^{1}$.  The advantage of this approach -- and the reason for passing to $\mathbb{P}^{1} \times \mathbb{P}^{1}$ in the first place -- is that each such locus is a linear subspace of $E$.  By varying the subschemes of $\mathbb{P}^{1}$ where the points and degeneracy loci are forced to vanish, we obtain a family of linear subspaces of $E$; together these families yield a simplicial resolution of $E$.

Note that the spaces in this simplicial resolution are fibered by linear spaces over certain labeled configuration spaces.  To count points, it is crucial to understand the dimension of the linear fibers.  Using conic bundle structures, we prove that fibers of low codimension always have the expected dimension.  However, the fibers with high codimension are harder to control and are parametrized by a locus that is difficult to understand using point counting arguments.  To resolve this issue, we work with \'etale cohomology and control the contributions of high codimension strata using a general result of Sawin and Shusterman about \'etale cohomology groups of moduli spaces of curves (see \cref{appendixSS}).

Following \cite{DT24}, we implement our inclusion-exclusion argument using a 
``homological sieve'' method based on posets and the bar complex.  Since we only understand the strata of low codimension, we truncate the bar complex and prove that the cohomology of the truncated bar complex approximates the cohomology of our original space in some range.  We then carry this comparison one step further by considering a ``virtual'' bar complex  where we compute what would happen if each linear fiber had exactly the expected dimension.  The final step is to show that an application of the Grothendieck-Lefschetz trace formula to the virtual bar complex  yields Peyre's constant; we construct and study a height zeta function that allows us to compare the local factors for the two quantities.

The following subsections highlight the various technical advances needed to execute this plan.

\subsubsection{Conic bundle structures and the space of rational curves} \label{sect:introconics}

Let $S$ be a split quartic del Pezzo surface over $k = \mathbb F_q$.  For clarity we henceforth let $B = \mathbb P^1$ denote our domain curve.
\cite{BLRT21} proved that for any nef curve class $\alpha$ the moduli space $M_\alpha = \mathrm{Mor}(B, S)_\alpha$ is geometrically irreducible and has the expected dimension $-K_S.\alpha + 2$.
For each $\alpha$ we show that one can find a birational morphism $\rho : S \to \mathbb P^1 \times \mathbb P^1$ with the following properties: let $E_1, \dotsc, E_4$ be the exceptional divisors contracted by $\phi$, and $F, F'$ be general fibers of two conic fibrations $S \to \mathbb P^1 \times \mathbb P^1 \to \mathbb P^1$ respectively.   
Then we have
\[
2F.\alpha \geq \sum_{i = 1}^4 E_i.\alpha, \quad 2F'.\alpha \geq \sum_{i = 1}^4 E_i.\alpha
\]
Let us fix such a birational morphism.
Let $(p_i, p_i') \in \mathbb P^1 \times \mathbb P^1$ be the four points which are blown-up by $\beta$.
Then $M_\alpha$ can be identified with the space of morphisms $B \to \mathbb P^1 \times \mathbb P^1$ of degree $(a, a') := (F.\alpha, F'.\alpha)$ such that the multiplicity at $(p_i, p_i')$ is equal to $k_i := E_i.\alpha$. Using this description and the above inequalities, one can show that the morphism 
\[
\Psi_\alpha: M_\alpha \to U_{\mathbf{k}}, \, [s : B \to S] \mapsto (s^*E_1, s^*E_2, s^*E_3, s^*E_4)
\]
is dominant where $\mathbf{k} = (k_1, k_2, k_3, k_4)$ and
\[
U_{\mathbf{k}} = \{ (T_i) \in \prod_{i = 1}^4 \mathrm{Hilb}^{[k_i]}(B) \, | \, \text{ for any $i \neq j$, the supports of $T_i$ and $T_j$ are disjoint.} \}.
\]

We study the structure of $M_{\alpha}$ using the map $\Psi_{\alpha}$.  First, we are able to prove an upper bound for the dimension of the complement $U_{\mathbf{k}} \setminus \Psi_{\alpha}(M_\alpha)$.  Second, we show that $\Psi_\alpha$ realizes $M_{\alpha}$ as a Zariski open subset of a relative product of projective bundles over an open subset of $U_{\mathbf{k}}$.

In particular, we have an explicit upper bound of the correct magnitude for the number of rational curves of class $\alpha$.
Altogether, these results allow us to prove an ``unobstructedness'' statement showing that certain subloci of $M_{\alpha}$ have the expected dimension.

\subsubsection{The method of the bar complexes}
We have identified $M_\alpha$ with the space of morphisms $B \to \mathbb P^1 \times \mathbb P^1$ of degree $(a, a') = (F.\alpha, F'.\alpha)$ such that the multiplicity at $(p_i, p_i')$ is equal to $k_i = E_i.\alpha$. One may rephrase this in the following way: let $\mathbb P^1 \times \mathbb P^1 = \mathbb P(V_1) \times \mathbb P(V_2)$ where $V_i$ is a $2$-dimensional vector space over $k$. Let $\ell_{i, j} \subset V_j$ be the $1$-dimensional subspaces corresponding to $p_i, p_i'$ respectively for $j = 1,2$. We consider the space $\widetilde{M}_\alpha$ of sections
\[
(s, t) \in \Gamma(B, V_1 \otimes \mathcal O(a)) \oplus \Gamma(B, V_2\otimes \mathcal O(a')),
\]
such that $s, t$ are nowhere vanishing and the length of $(s, t)^{-1}(\ell_{i, 1} \otimes \mathcal O(a) \oplus \ell_{i, 2}\otimes \mathcal O(a'))$ is $k_i$ for $i = 1, \dotsc, 4$. Then there is a natural morphism
\[
\widetilde{M}_\alpha \to M_\alpha
\]
realizing $\widetilde{M}_{\alpha}$ as a $\mathbb G_m^2$-torsor over $M_\alpha$. 
So one can reduce the counting problem for the number of $k$-rational points on $M_\alpha$ to the counting problem for the number of $k$-rational points on $\widetilde{M}_\alpha$.
As mentioned before, we have a dominant morphism $\widetilde{M}_\alpha \to U_{\mathbf{k}}$, and for each $w \in U_{\mathbf{k}}$, the fiber $\widetilde{M}_{\alpha, w}$ is realized as a Zariski open subset of 
\[
E_w = \Gamma(B, V_1 \otimes \mathcal O(a) \oplus V_2 \otimes \mathcal O(a'))_w,
\]
which is the space of sections $(s, t)$ with the incidence conditions $s(w_i) \subset \ell_{i, 1}\otimes \mathcal O (a)$ and $t(w_i) \subset \ell_{i, 2} \otimes \mathcal O(a')$.  Since the fibers of $E \to U_{\mathbf{k}}$ generically
have constant rank, in this way we almost realize $\widetilde{M}_\alpha$ as a Zariski open subset of a vector bundle over $U_{\mathbf{k}}$.

Our next goal is to understand the number of $k$-rational points on the complement
\[
E \setminus \widetilde{M}_\alpha.
\]
This complement can be described as a union of locally closed strata coming from moduli spaces of curves of lower degree.  Thus we are well-situated to apply the inclusion-exclusion principle.  The formalization of this idea is the method of bar complexes developed in \cite{DT24} over $\mathbb C$.  While in \cite{DT24} the bar complexes are certain topological simplicial spaces, in this paper we develop the theory of bar complexes over an arbitrary perfect field using the language of simplicial schemes and their homotopy theory. The advantage of this perspective is that we formalize the intersection properties of stratifications and point counting becomes cleaner in this setting.  In particular, one can appeal to the Grothendieck--Lefschetz trace formula for simplicial schemes to translate \'etale cohomology computations into point counts.

We think of the Hilbert scheme $\mathrm{Hilb}(B)$, ordered by inclusion, as a \emph{poscheme}, i.e., a scheme with a poset structure. 
Consider the following finite poset $Q$, which we also think of as a poscheme (with the reduced scheme structure):
\begin{equation*}
Q:= \{ V, V_1 \oplus \{0\}, \{0\} \oplus V_2,  0 \} \cup \{ \ell_{i, 1} \oplus \ell_{i, 2}, \ell_{i, 1} \oplus \{0\}, \{0\} \oplus \ell_{i, 2}  \}_{i=1,\dotsc,4} \,,
\end{equation*}
ordered by inclusion of subspaces.
We let $\mathrm{Hilb}(B)^Q$ denote the scheme parametrizing poset homomorphisms
\[
x : Q \to \mathrm{Hilb}(B),
\]
such that $x^{-1}(B) = V$.  
Then the complement $E \setminus \widetilde{M}_\alpha$ admits a stratification into linear subspaces parametrized by the poscheme $\mathrm{Hilb}(B)^Q$.

More precisely, we embed $U_{\mathbf{k}}$ into $\mathrm{Hilb}(B)^Q$ by mapping $w \mapsto x$ with $x_{ \ell_{i, 1} \oplus \ell_{i, 2}} = w_i$, $x_{V} = B$, and $x_q = \emptyset$ for other $q \in Q$.  We then define
\[
(U_{\mathbf{k}} < \mathrm{Hilb}(B)^Q)  = \{ (w < x) \, | \, w \in U_{\mathbf{k}}, \, x \in \mathrm{Hilb}(B)^Q\},
\]
We construct our stratification as follows: for any $(w < x) \in (U_{\mathbf{k}} < \mathrm{Hilb}(B)^Q)(\overline{k})$, the subspace
\[
Z_{w < x} =  \Gamma(B, V_1 \otimes \mathcal O(a) \oplus V_2 \otimes \mathcal O(a'))_x \subset E_w,
\]
is the space of sections satisfying the incidence condition imposed by $x$.
Then these subspaces are bundled as 
\[
Z \subset (U_{\mathbf{k}} < \mathrm{Hilb}(B)^Q) \times_{U_{\mathbf{k}}} E.
\]
We stratify $(U_{\mathbf{k}} < \mathrm{Hilb}(B)^Q) = \sqcup_T \mathcal N_T$ into locally closed subsets based on the combinatorial type $T$ (see \cref{defn:combinatorial-type}).  Only certain types $T$ -- the ``saturated'' types (see \cref{defi:saturated}) -- impose non-redundant conditions on the space of sections.  We then truncate the poset: we choose a closed subscheme $P \subset (U_{\mathbf{k}} < \mathrm{Hilb}(B)^Q)$ which is obtained by taking a finite union of loci $\mathcal{N}_{T}$ where the expected codimension is not too large and then take the downward closure under the poset structure.  The bar complex $B(P, Z)$ is a simplicial scheme whose objects are described by 
\[
B(P, Z)([n]) = \{ w < x_0 \leq \dotsb \leq x_n, z \in Z_{w < x_n} \, | \, (w < x_i) \in P \}.
\]
The following theorem has been proved over $\mathbb C$ in \cite[Theorem 5.9]{DT24}, and we establish this over an arbitrary perfect field:

\begin{theo}[a weaker version of \cref{theo:5.4inDT}] \label{theo:intro54}
\label{theo:barcomplex_intro}
Let $I \in \mathbb N$. Assume that
\begin{enumerate}
\item $P \to U_{\mathbf{k}}$ is proper and $P$ is downward closed;
\item for every pair $(w < x)\in P(\overline{k})$ with $x$ saturated and for every saturated $y \in \mathrm{Hilb}(B)^Q(\overline{k})$ such that $x \prec y$,  
the fiber $Z_{w < y}$ has the expected dimension, and;
\item $P$ contains every $\mathcal N_T$ such that the expected cohomological codimension function $\kappa$ satisfies $\kappa(T) \leq I$.
\end{enumerate}
Then for all $i > 2\mathrm{dim}(E) - I -2$, the map 
\[
H^i_{\text{\'et}, c}(\mathrm{im}((Z)_{\overline{k}} \to (E)_{\overline{k}}), \mathbb Z_\ell) \to H^i_{\text{\'et}, c}(B(P, Z)_{\overline{k}}, \mathbb Z_\ell),
\]
is an isomorphism. 
\end{theo}

Assumption (2) of the above theorem is almost verified in the stability range given by $I = \lfloor \frac{1}{8}\min\{2a - \sum_i k_i, 2a' - \sum_i k_i\}-\frac{1}{2} \rfloor$ by our analysis of $M_{\alpha}$ described in the previous subsection. (More precisely, we need to shrink $U_{\mathbf{k}}$ to a Zariski open subset $U$ such that $U_{\mathbf{k}} \setminus U$ has high codimension.  Since the complement has high codimension, this is negligible for our discussions and so we ignore this issue for now.)  Then \cref{theo:intro54} guarantees that the cohomology of $E \setminus \widetilde{M}_{\alpha}$ is approximated by the corresponding bar complex.

One advantage of working with the bar complex is that we can compute its cohomology using the following spectral sequence.  An important feature of this spectral sequence is that only certain combinatorial types $T$ (i.e.~the ``essential saturated'' types; see \cref{defi:essentialtype}) can contribute.
\begin{theo}[\cref{theo:spectralsequenceforbar}] \label{theo:introspectralsequence}
There is a spectral sequence:
\[
\bigoplus_{i + j = n} E_1^{i, j} = \bigoplus_{T: \text{ essential, saturated type of $P$}} H^{n}_{\text{\'et}, c}((Z|_{\mathcal N_T})_{\overline{k}}, \mu'(T)[1]) \implies H^{n}_{\text{\'et}, c}(B(P, Z)_{\overline{k}}, \mathbb Z_\ell),
\]
where $\mu'(T)$ is a complex in the derived category of $\mathbb Z_\ell$-constructible sheaves in the pro-\'etale topology whose stalk at $(w<x)$ captures the compactly supported cohomology of the nerve of the poset interval $(w,x)$.
\end{theo}

The Grothendieck--Lefschetz trace formula allows us to compute the number of $\mathbb{F}_{q}$-points on a variety by analyzing the Frobenius action on cohomology.  Thus our goal is to understand, up to error terms, the truncated alternating sum of the Frobenius traces:
\[
 \sum_{i \geq 4a + 4a' -2\sum_j k_j -I + 10} (-1)^i \mathrm{Tr}(\mathrm{Frob} \curvearrowright H^i_{\text{\'et}, c}(B(P, Z)_{\overline{k}}, \mathbb Q_\ell)).
\]
The spectral sequence of \cref{theo:introspectralsequence} reduces this to a computation over essential saturated types:
\begin{equation*}
 \sum_{T: \text{ essential type}}   \sum_{i \geq 4a + 4a' -2\sum_j k_j -I + 10} (-1)^i \mathrm{Tr}(\mathrm{Frob} \curvearrowright H^{i}_{\text{\'et}, c}((Z|_{\mathcal N_T})_{\overline{k}}, \mu'(T)[1]\otimes \mathbb Q_\ell)).
\end{equation*}
Assuming that $Z$ has the expected dimension over $\mathcal N_T$, one can use the Grothendieck--Lefschetz trace formula for simplicial schemes to conclude that, up to error terms, the above quantity is equal to
\[
-q^{2a+ 2a'+4} \sum_{(w< x) \in (U_{\mathbf{k}} \leq Q^{JB})(k)} \mu_k(w, x)q^{-\gamma(x)},
\]
where $Q^{JB} \subset \mathrm{Hilb}(B)^Q$ is the Zariski open subset parametrizing saturated elements, $\mu_k$ is the M\"obius function for the poset $Q^{JB}(k)$, and $\gamma(x)$ is the expected codimension of the incidence condition imposed by $x$.
Putting everything together, the quantity $\#\widetilde{M}_{\alpha}(k)$ is equal to
\[
q^{2a+ 2a'+4} \sum_{(w\leq x) \in (U_{\mathbf{k}} \leq Q^{JB})(k)} \mu_k(w, x)q^{-\gamma(x)},
\]
up to error terms.

\subsubsection{The virtual height zeta function and Peyre's constants}

Let $\alpha$ be an ample curve class on a split quartic del Pezzo surface $S$.  Using the invariants $a(\alpha), a'(\alpha), k_i(\alpha)$, we assume that 
\[
2a(\alpha) - \sum_{i = 1}^4 k_i(\alpha) > 0, 2a'(\alpha) - \sum_{i = 1}^4 k_i(\alpha) > 0.
\]
We showed that 
\[
\#\widetilde{M}_{\alpha}(k) \sim q^{2a+ 2a'+4} \sum_{(w\leq x) \in (U_{\mathbf{k}} \leq Q^{JB})(k)} \mu_k(w, x)q^{-\gamma(x)}
\]
as $a, a', k_i \to \infty$. On the other hand, Manin's conjecture predicts that
\begin{align*}
\#\widetilde{M}_{\alpha}(k) &\sim (1-q^{-1})^{-4} q^{2a  + 2a' - \sum_i k_i+4} \tau_{-K_S}(S) \\
&= (1-q^{-1})^{-4} q^{2a + 2a' - \sum_i k_i + 4} \prod_{c \in |B|}(1-q^{-|c|})^6(1+6q^{-|c|} + q^{-2|c|}),
\end{align*}
as $a, a', k_i \to \infty$.
To this end, we consider the following virtual height zeta function: 
\[
\mathsf Z(\mathbf{t}) = \sum_{k_1, \dotsc, k_4 = 0}^\infty q^{\sum_{i = 1}^4k_i}\left( \sum_{(w\leq x) \in (U_{\mathbf{k}} \leq Q^{JB})(k)} \mu_k(w, x)q^{-\gamma(x)} \right)t_1^{k_1}t_2^{k_2}t_3^{k_3}t_4^{k_4}.
\]
Since the M\"obius function is multiplicative, the above zeta function can be rewritten as the following convergent Euler product: 
\begin{multline*}
\prod_{c \in |B|}\biggl( 1 - 6q^{-2|c|} + 8q^{-3|c|} -3q^{-4|c|} \\+ \sum_{i = 1}^4\sum_{d = 1}^\infty (qt_i)^{|c|d} (q^{-2d|c|}  - 2q^{-(2d+1)|c|} + 2q^{-(2d + 3)|c|} - q^{-(2d + 4)|c|})\biggr).
\end{multline*}
Then we have the following proposition:

\begin{prop}[\cref{prop:formalheightzeta}]
 There exists a $\delta > 0$ such that the series
\begin{equation*}
\prod_{i = 1}^4(1-t_i)\mathsf Z(\mathbf{t})
\end{equation*}
absolutely converges when $|t_i|\leq q^{\delta}$ for $i = 1, \dotsc, 4$.
Moreover we have
\[
\lim_{t_i \to 1} \prod_{i = 1}^4 (1-t_i) \cdot \mathsf Z(\mathbf{t}) = (1-q^{-1})^{-4}\prod_{c \in |B|}(1-q^{-|c|})^6(1+6q^{-|c|} + q^{-2|c|}).
\]
\end{prop}

Thus Abel summation 
implies that
\[
q^{\sum_i k_i}\sum_{(w\leq x) \in (U_{\mathbf{k}} \leq Q^{JB})(k)} \mu_k(w, x)q^{-\gamma(x)} \sim (1-q^{-1})^{-4}\prod_{c \in |B|}(1-q^{-|c|})^6(1+6q^{-|c|} + q^{-2|c|}),
\]
as $k_i \to \infty$.
Therefore \cref{theo:intromaninconj} has been verified modulo the error terms.

\subsubsection{Uniform estimates for error terms by Sawin--Shusterman based on Katz's results}

To handle the error terms, we use the following theorem by Sawin--Shusterman based on Katz's results:
\begin{theo}[Sawin--Shusterman in \cref{appendixSS}] \label{theo:introerrorterms}
Let $X$ be a projective variety with an ample divisor $H$.
Then there exists a constant $C > 1$ such that for any $\alpha \in N_1(X)_{\mathbb Z}$ we have
\[ 
\sum_{i=0}^\infty \dim H^i_{\text{\'et}, c} ( \mathrm{Mor}(\mathbb P^1, X, \alpha), \mathbb Q_\ell) \leq C^{H.\alpha+1}.\]
\end{theo} 
In our setting, we let $H = -K_S$. Assuming that $\alpha$ is in the cone $\Nef_1(S)_{\ell, \epsilon}$, our stability range satisfies that
$\frac{1}{32}\min\{ 2a(\alpha) - \sum_ik_i(\alpha), 2a'(\alpha)-\sum_ik_i(\alpha)\} \geq - \epsilon K_S.\alpha$.
To beat the above error estimates
we need to assume that $q^{\epsilon} > C$.  This forces us to restrict to the modified cone $\Nef_1(S)_{\ell, \epsilon}$ instead of the entire nef cone.

One might ask why one needs to consider compactly supported cohomology instead of directly conducting point counting using the inclusion-exclusion principle. A reason for this is that to show that our stratification has expected behavior in certain codimension range, we need to shrink $U_{\mathbf{k}}$ to a certain Zariski open subset $U$. Unfortunately we do not have much access to this $U$ other than knowing that $U_{\mathbf{k}} \setminus U$ has high codimension. To do a naive point counting, one needs to estimate the number of $\mathbb F_q$-points on this complement which is a formidable task from our experience. Fortunately Sawin--Shusterman's result shows that we only need to understand high degree cohomologies and these cohomologies are not affected by shrinking $U_{\mathbf{k}}$ to $U$ because the complement has high codimension. This is one of main motivations why we appeal to compactly supported cohomology and Grothendieck--Lefschetz trace formula.

\subsection{Comparison to \cite{DT24}}

\cite{DT24} studies the Cohen--Jones--Segal conjecture for blow ups of the projective space at finitely many points.  Some key breakthroughs in \cite{DT24} are:
\begin{itemize}
\item the realization that birational geometry yields a stratification of the moduli space of curves that is well-suited for the inclusion-exclusion principle using the bar complex;
\item the study of the posets $\Hilb(C)^Q$, $Q^{JC}$, their combinatorial stratifications, and a criterion for the cohomology of bar complexes to converge;
\item the introduction of spaces of positive and semi-topological maps as ``intermediate'' steps interpolating between the spaces of holomorphic and continuous maps, and; 
\item techniques for verifying that the cohomology of the bar complex agrees with the cohomology of the space of positive maps in an optimal range, via finite dimensional approximation and spaces of semi-topological maps.
\end{itemize}
In particular these innovations lead to the proof of Cohen--Jones--Segal conjecture for quintic del Pezzo surfaces.

Our results also rely on the inclusion-exclusion principle using the bar complex.  One important technical change is that our proof of \cref{theo:intromaninconj} requires the bar complex method over perfect fields -- this involves a translation of the method developed in \cite{DT24} into the language of simplicial schemes and their homotopy theory.  In addition, our proofs contain the following new ideas and developments: 
\begin{itemize}
\item the realization that the ``homological sieve'' method can be used to obtain control of error terms even when naive point counting does not apply;
\item algebro-geometric proofs of stronger unobstructedness results for the moduli spaces of rational curves using conic bundle structures on split quartic del Pezzo surfaces;   
\item   
a description of the main term of the counting function in terms of the combinatorics of posets and their M\"obius functions by applying the Grothendieck--Lefschetz trace formula for simplicial schemes to the bar complex, and;
\item the introduction of the virtual height zeta function revealing the compatibility of the combinatorial main term and Peyre's constant.
\end{itemize}
In particular, our approach leads to tight control of the 
behavior of the number of $\mathbb F_q$-points on the moduli space of rational curves of a given class.  The techniques developed in this paper open up new geometric and topological approaches to Manin's conjecture over global function fields.

\subsection{Open questions}

Here we record some possible follow-ups to this paper.

\subsubsection{Motivic Manin's conjecture}

There have been extensive activities on motivic Manin's conjecture, which concerns the convergence of the moduli space of rational curves in Grothendieck ring of varieties (see \cite{Bou09, CLL16, Bilu23, BDH22, Faisant23, Faisant25}).  It would be interesting to explore our homological sieve method in the motivic setting. 

\subsubsection{Stabilization maps}
Another interesting question is to explore the realization of stabilization maps exhibiting homological stability established by \cref{theo:introcjs}. We should emphasize that we only establish homological stability over $\mathbb C$, and in particular we did not establish stability over finite fields as Galois representations even though we show the convergence of point counting for the spaces of rational curves. It would be interesting to prove Galois equivariant homological stability; one possible approach would be using log geometry and the Kato-Nakayama formalism introduced by Vaintrob in \cite{Vaintrob19, Vaintrob21} and applied in \cite{BDPW,EL23, LL25}. We plan to come back to this issue in our future work.

\subsection{Related works}

\subsubsection{Manin's conjecture}

Manin's conjecture for smooth del Pezzo surfaces over number fields has been extensively studied. \Cref{conj:Manin} for smooth del Pezzo surfaces of degree $\geq 6$ follows from \cite{BTtoric96, BTtoric98} because these surfaces are toric. \Cref{conj:Manin} for smooth split  del Pezzo surfaces of degree $5$ has been proved over $\mathbb Q$ in \cite{delaBre} and recently over arbitrary number fields in \cite{BD24}. However, the conjecture is largely open for smooth del Pezzo surfaces of degree $\leq 4$: there is only one confirmed case of Manin's conjecture for quartic del Pezzo surfaces (\cite{DBB11}) and there is no example of proved cases of the conjecture for del Pezzo surfaces of degree $\leq 3$. \cite{BS19} established upper bounds and lower bounds of the correct magnitude for quartic del Pezzo surfaces with conic bundle structures and \cite{FLS} established lower bounds of correct magnitude for del Pezzo surfaces with the same structures. These approaches using conic bundles are inspirational to us and motivated the material in \cref{sect:introconics}.

There are fewer results on Manin's conjecture over global function fields. As mentioned before, there is a very influential heuristic argument by Batyrev \cite{Bat88} combined with Ellenberg--Venkatesh's proposal of using homological stability in this context \cite{EV05}. Manin's conjecture over global function fields has been proved for toric varieties (\cite{Bourqui03, Bourqui11}) and equivariant compactifications of forms of $\mathbb G_a^n$ (\cite{formsofGa}) using harmonic analysis and for low degree hypersurfaces (\cite{BV16, BW23, Sawin24}) and complete intersections (\cite{Glas24}) using the circle method. There are also approaches using universal torsors which were pioneered by Bourqui in \cite{BourquiD09, Bou11, Bou13}.
\cite{GH24, Glas25} obtain certain upper bounds for the counting functions of rational points on low degree del Pezzo surfaces, but this does not provide upper bounds of optimal magnitude.

\subsubsection{Cohen--Jones--Segal conjecture}

\cite{Segal} showed the cohomology of the moduli space of degree $d$ pointed rational curves on $\mathbb{P}^{n}$ stabilizes as $d$ increases to the cohomology of the space of pointed continuous maps $S^{2} \to \mathbb{P}^{n}$.  Segal's result was extended to various homogenous spaces such as projective toric varieties (\cite{Guest}), Grassmannians (\cite{Kirwan}), generalized flag varieties (\cite{BMHM94}), and varieties admitting actions by connected solvable Lie groups with dense orbits (\cite{BHM01}).  \cite{CJS00} predicted an analogous statement for all Fano varieties $X$.  More recently, Browning--Sawin used the circle method to construct a spectral sequence converging to the compactly supported cohomology of the space of maps to low degree affine hypersurfaces. They showed that the degeneration of this spectral sequence is closely related to the Cohen--Jones--Segal conjecture (\cite{BS20}).  In \cite{DT24} the first and fourth authors studied the Cohen--Jones--Segal conjecture for spaces $X$ which are the blow-ups of projective space at a collection of points.  Their results establish the conjecture under certain conditions on the numerical class $\alpha$.  In particular, their work settles the conjecture for $100\%$ of nef classes in the case of quintic del Pezzo surfaces but only applies to a certain subset of the nef cone for quartic del Pezzo surfaces.

\subsubsection{Topological approaches to analytic number theory over global function fields}

There is a vast literature on topological approaches to analytic number theory over global function fields. This has been first pioneered in \cite{EVW} which explored the Cohen--Lenstra heuristics over global function fields.   The homological stability of Hurwitz spaces plays a crucial role, and there are many subsequent works on the same topic, e.g., applications of homological stability to Malle's conjecture \cite{ETW, LL25}, heuristics on Selmer groups \cite{EL23}, and the Cohen--Lenstra heuristics \cite{LL24a, LL24b, LL25}.  An application of stable homology of the braid group to moments of families of quadratic L-functions was studied in \cite{BDPW, MPPRW}. There are fewer results applying homological stability to Manin's conjecture; for example, \cite{Banerjee} shows an approach to Manin's conjecture for toric varieties using a topological method that is also a version of inclusion-exclusion.

\bigskip

\

\noindent
{\bf Acknowledgements:}
The authors are grateful to Will Sawin and Mark Shusterman for writing an appendix to this paper.
The authors would like to thank Kazuhiro Fujiwara and Will Sawin for helpful discussions, and Brendan Hassett, J\'anos Koll\'ar, Lars Hesselholt, Aaron Landesman, Yuri Tschinkel, Akshay Venkatesh, and Melanie Wood for comments on our work.
The authors would also like to thank Remy van Dobben de Bruyn, John Calabrese, Lars Hesselholt, Anton Mellit, Dan Petersen, and Andy Putman for answering the authors' questions.
The authors would like to thank Tim Browning, Jordan Ellenberg, Lo\"is Faisant, and Aaron Landesman for their comments on a draft of this paper.

The second author was supported by Simons Foundation grant Award Number 851129.  The third author was partially supported by JST FOREST program Grant number JPMJFR212Z, by JSPS KAKENHI Grant-in-Aid (B) 23K25764, by JSPS Bilateral Joint Research Projects Grant number JPJSBP120219935, by JSPS KAKENHI Early-Career Scientists Grant number 19K14512, by Inamori Foundation, and by MEXT Japan, Leading Initiative for Excellent Young Researchers (LEADER).

\section{Preliminaries}
\label{sec:preliminaries}

\noindent
{\bf Notation}: The natural numbers $\mathbb{N}$ include $0$.

The ground field is denoted by $k$ and a variety over $k$ is an integral separated scheme of finite type over $k$.  Given a scheme $X$ over $k$, a component of $X$ will mean an irreducible component unless we state otherwise. We will denote the base change of $X\to \Spec \, k$ to $\Spec \, \overline{k}$ by $X_{\overline{k}}$.

For a projective variety $X$, we denote the space of $\mathbb R$-Cartier divisors up to numerical equivalence by $N^1(X)$ and the space of $\mathbb R$-$1$-cycles up to numerical equivalence by $N_1(X)$. Inside these real vector spaces, we denote the lattice generated by integral Cartier divisors (resp.~integral $1$-cycles) by $N^1(X)_\bZ$ (resp.~$N_1(X)_\bZ$).
The nef cones of divisors and curves are denoted by $\Nef^1(X) \subset N^1(X)$ and $\Nef_1(X) \subset N_1(X)$ respectively, and the pseudo-effective cones of divisors and curves are denoted by $\overline{\mathrm{Eff}}^1(X) \subset N^1(X)$ and $\overline{\mathrm{Eff}}_1(X) \subset N_1(X)$ respectively.
For any closed cone $\mathcal C \subset N_1(X)$, we denote $\mathcal C\cap N_1(X)_{\mathbb Z}$ by $\mathcal C_{\mathbb Z}$.

\begin{lemm} \label{lemm:smoothnesscriterion}
Let $f: X \to Y$ be a $k$-morphism of finite-type separated $k$-schemes such that $X$ is equidimensional and $Y$ is irreducible. 
Suppose that $Y$ is smooth and that every fiber of $f$ over a closed point is smooth of dimension $\dim(X) - \dim(Y)$.  Then $X$ is a smooth scheme and $f$ is a smooth morphism.
\end{lemm}

\begin{proof}
 
Let $x \in X$ be a closed point.  Since the fiber of $f$ containing $x$ is smooth, $\rk_{x} \Omega_{X/Y} = \dim(X) - \dim(Y)$.  Since $Y$ is smooth, the exact sequence
\begin{equation*}
f^{*}\Omega_{Y/k} \otimes k(x) \to \Omega_{X/k} \otimes k(x) \to \Omega_{X/Y} \otimes k(x) \to 0
\end{equation*}
shows that $\rk_{x}\Omega_{X/k} \leq \dim(X)$.  Since $\dim(X)$ is also a lower bound on the rank, equality is achieved.  Thus $X$ is smooth at every closed point, hence smooth.

By upper semicontinuity of fiber dimension on the domain, every fiber of $f$ has dimension $\dim(X)-\dim(Y)$.  Since $X,Y$ are smooth, we conclude that $f$ is flat by \cite[Tag 00R4]{stacks}. 
Finally, a flat morphism of finite-type $k$-schemes such that the fibers over closed points are smooth is a smooth morphism.
\end{proof}

\subsection{Bounds on traces of endomorphisms}

In this section, we record certain estimates for the sum of eigenvalues of an endmorphism, which we will eventually apply to Frobenius acting on étale cohomology:

\begin{defi}
Suppose that $K$ is a field equipped with a norm $|\cdot|: \overline K \to \mathbb R$.  For any finite-dimensional $K$-vector space $V$ with an endomorphism $F$, the $\mathsf L^{1}$-trace of $F$ on $V$ is
$$|F,V| = \sum_{\lambda \in \overline K}\dim((V\otimes \overline K)_{\lambda}) |\lambda|,$$ where  $(V\otimes \overline K)_{\lambda}$ denotes the generalized eigenspace of $F$ with eigenvalue $\lambda$.
\end{defi}

Note that $|\Tr(F \curvearrowright V)| \leq |F,V|$ and if $W$ is an $F$-subquotient 
then $|F,W| \leq |F,V|$.   Furthermore the $\mathsf L^{1}$-trace is additive in short exact sequences.

Our first lemma addresses the behavior of spectral sequences when taking the trace of an endomorphism.  Let $E_1^{i,j}$ be a spectral sequence of $K$-vector spaces which are nonzero for only finitely many $i,j$,  converging to $H^{i+j}$ equivariantly with respect to the action of an endomorphism $F$.  We denote the differentials in this spectral sequence by $d_{r}$.

\begin{lemm}
\label{lemm:spectral}
	For any integer $I$, we have the inequality $$|\sum_{n\geq I}(-1)^n\Tr(F\curvearrowright \bigoplus_{i+j =n} E_1^{i,j}) - \sum_{n \geq I} (-1)^n\Tr(F \curvearrowright H^n)| \leq |F, \bigoplus_{i+j = I} E_1^{i,j}|.$$  
\end{lemm}
\begin{proof}
	We write $(E_r)^m$ for the group $\bigoplus_{i+j = m} E_r^{i,j}$ and we write $B_r^n$ for $d_r((E_r)^{n-1})$, i.e.~for the total boundary group in cohomological degree $n$ 
	on the $r$-th page.  We have that 
	$$(-1)^{I-1} \Tr(F\curvearrowright B_r^I) + \sum_{n \geq I} (-1)^n \Tr(F \curvearrowright (E_{r})^n)  = \sum_{n \geq I} (-1)^n \Tr(F \curvearrowright (E_{r+1})^n).$$
	Since there are finitely many nonzero terms on the $E_{1}$ page, the spectral sequence degenerates at a finite stage $R$ so that $E_R = E_\infty$.  Inductively  we obtain
	\begin{align*}
	\sum_{n \geq I} (-1)^n\Tr(F \curvearrowright H^n) &=\sum_{n \geq I} (-1)^n\Tr(F\curvearrowright (E_R)^n) \\
	&= (-1)^{I-1}\sum_{r = 1}^R \Tr(F\curvearrowright B_r^I) + \sum_{n \geq I} (-1)^n \Tr(F\curvearrowright (E_1)^n).
	\end{align*}
	To bound $\sum_{r = 1}^R \Tr(F\curvearrowright B_r^I)$, note that there is a filtration $0 \subseteq G_1 \subseteq G_2 \subseteq  \dotsb \subseteq G_R$ of $(E_1)^I$ such that the associated graded $G_r/G_{r-1} = B_r^I$.  Hence we have that $$ \left| \sum_{r = 1}^R \Tr(F\curvearrowright B_r^I) \right| \leq \sum_{r = 1}^R |F,B_r^I|  = |F,G_R| = |F, (E_1)^I|$$  
	as desired. 
	\end{proof}
	
Our next goal is to bound the $\mathsf L^{1}$-trace under various tensor operations.  Let $W = \oplus_{(n,i) \in \mathbb{N}^{2}} W_{n,i}$ be an $\mathbb{N}^{2}$-graded $K$-vector space.  We will be interested in the case when $W$ is equipped with an endomorphism $F$ that respects the grading.

\begin{lemm}\label{singletensor}%
	Let $V,W$ be $\mathbb{N}^{2}$-bigraded vector spaces equipped with endomorphisms $F_{V}$ and $F_{W}$ which are compatible with the grading.  Suppose there are positive constants $d,b$ such that $|F_{V},V_{n,i}| \leq C_{1} d^n b^i$ and $|F_{W},W_{n,i}| \leq C_{2} d^n b^i$.  Then $$|F_{V} \otimes F_{W}, (V \otimes W)_{n,i}| \leq C_{1}C_{2} (n+1)(i+1) d^n b^i.$$ 
\end{lemm}
\begin{proof}
	We have
	\begin{align*}
	|F_{V} \otimes F_{W}, (V \otimes W)_{n,i}|  & = 
	\sum_{(n,i) = (n_1,i_2) + (n_2, i_2)}  |F_{V}, V_{n_1, i_1}| |F_{W}, W_{n_2,i_2}|  
	\end{align*}
	and the sum is indexed over a set of size $(n+1)(i+1)$.
\end{proof}

\begin{lemm}\label{tensoralgebrabound}
	Let $W$ be an $\mathbb{N}^{2}$-bigraded vector space with endomorphism $F$ which is compatible with the grading.  Suppose that $W_{0, i} = W_{n,0} = 0$ for all $n,i$ and that there are positive constants $E,d,b$ with $E \geq 1$ such that  $|F, W_{n,i}| \leq E d^n b^i$.  Then the free associative algebra $T(W) = \bigoplus_{r} W^{\otimes r}$ satisfies $|F,T(W)_{n,i}| \leq (2Ed)^n (2b)^i$.
\end{lemm}
\begin{proof}
We have 
	$$|F, T(W)_{n,i}| = \sum_{r \geq 0} \sum_{(n,i) = (n_1,i_1) + \dotsb + (n_r,i_r)} \prod_{s = 1}^r |F, W_{n_s,i_s}|.$$
	By our vanishing hypothesis on $W$, the only non-zero contributions occur when $r \leq n$.
	In particular, since we have $E^{r} \leq E^{n}$ in this range, we conclude that
	$$|F, T(W)_{n,i}| \leq \sum_{r \geq 0} \sum_{(n,i) = (n_1,i_1) + \dotsb + (n_r,i_r)}  E^n d^n b^i. $$  
	Since a positive integer $t$ has $2^{t-1}$ compositions (i.e.~ordered partitions) there are no more than $2^{n-1} 2^{i-1}$ tuples satisfying the hypotheses. 
\end{proof}

\begin{lemm}\label{symbound}
Let $W$ be an $\mathbb{N}^{2}$-bigraded vector space with an endomorphism $F$ that is compatible with the grading.  Suppose that:
\begin{itemize}
\item $W_{0,i} = 0$ for all $n$,
\item $W_{\bullet, 0} := \bigoplus_{n} W_{n,0}$ is finite dimensional, and; 
\item there are positive constants $E,d,b$ with $E \geq 1$ such that $|F, W_{n,i}| \leq E d^n b^i$.
\end{itemize}
Let  $D > 1$ be a constant such that $|F,W_{\bullet,0}| \leq D$.
Then the free symmetric algebra $\Sym(W)$ satisfies $|F,\Sym(W)_{n,i}| = O((n + 1)(i + 1) (\max(2Ed, D))^n (2b)^i)$ (where the implied constants are independent of $n,i$).
\end{lemm}
\begin{proof}
We have that $\Sym(W) = \Sym(W_{\bullet, 0}) \otimes \Sym(\bigoplus_{i > 0} W_{\bullet, i})$.  We may apply \cref{tensoralgebrabound} and appeal to the surjection $T(V) \to \Sym(V)$ to obtain that $|F, (\Sym(\bigoplus_{i > 0} W_{\bullet, i}))_{n,i}| \leq (2Ed)^n (2b)^i$.    Furthermore $$|F, \Sym(W_{\bullet, 0})_{n,0}| \leq |F, \Sym^{ \leq n}(W_{\bullet,0})| \leq \sum_{j\leq n} D^j = O(D^n).$$   So applying \cref{singletensor} to the the tensor decomposition of $\Sym(W)$ we obtain the result.
\end{proof}

\subsection{Homology of configuration space}
Given a poset $P$ with a bottom element $\hat{0}$ and a top element $\hat{1}$, we denote by $H^*(P)$ the reduced cohomology with $\mathbb Q_\ell$-coefficients of the nerve of the open interval $(\hat{0}, \hat{1})$ in $P$.  In particular, when $X$ is a finite set then the poset $\mathrm{P}(X)$ of set partitions admits a bottom and top element and thus it makes sense to discuss $H^{*}(\mathrm{P}(X))$.

\cite[Example 3.14]{Petersen} gives a description of an upper bound for the compactly supported cohomology of the ordered configuration space.  This example is stated when the base field is $\mathbb C$, however the analogous statement over arbitrary fields follows from \cite[Theorem 3.3 (ii)]{Petersen}:

\begin{theo}\label{ConfigCohom}
	There is a Frobenius and $S_n$-equivariant spectral sequence converging to $H^*_{\text{\'et}, c}(\mathrm{Conf}_n \mathbb P^1_{\overline{k}}, \mathbb Q_\ell)$ with total $E_1$ page given by  
	$$\bigoplus_{\substack{p  \text{ a set partition} \\ \text{of degree $n$}}}    \bigotimes_{\lambda \in p} \left(H_{\text{\'et}}^*(\mathbb P^1_{\overline{k}}, \mathbb Q_\ell)  \otimes H^{|\lambda|-3}(\mathrm{P}(\lambda))[-|\lambda|+1] \right).$$ 
	Summing over all $n$, we obtain a spectral sequence converging to $\oplus_{n} H^*_{\text{\'et}, c}(\mathrm{Conf}_n \mathbb P^1_{\overline{k}}, \mathbb Q_\ell)$ with total $E_1$ page
	$$\mathrm{Sym}\left( H_{\text{\'et}}^*(\mathbb P^1_{\overline{k}}, \mathbb Q_\ell)[1] \otimes \left(\bigoplus_{n} H^{n - 3}(\mathrm{P}(\{1, \dotsc, n\}))[-n] \right) \right),$$ 
	where the symmetric power takes place in the category of symmetric sequences of bigraded vector spaces with Frobenius action.
\end{theo}

We need the following upper bound on the poset homology of the partition lattice.
\begin{prop}[{\cite[Theorem 7.3]{Stanley}}]
         Assume that $K$ is a field of characteristic $0$.
	For all $n \geq 1$, there is an $S_n$-equivariant quotient
	\[
	\overline{K}S_n \to H^{n - 3}(\mathrm{P}(\{1,\dotsc,n\}),\overline{K}),
	\]
	 where $S_n$ acts by left multiplication.
\end{prop}

In particular, for any vector space $V$, there is a surjection 
$$V^{\otimes n} = \overline{K}S_n \otimes_{S_n} V^{\otimes n} \to  H^{n - 3}(\mathrm{P}(\{1,\dotsc, n\}), \overline{K}) \otimes_{S_n} V^{\otimes n}.$$
Combining this with \cref{ConfigCohom} we obtain the following upper bound on $\mathsf L^1$-traces.

\begin{prop} \label{SchurBound}
        Assume that our ground field $k$ is $\mathbb F_q$.
	Let $V$ be a bigraded vector space over $K$.
Suppose that $F$ acts as an endomorphism on $V$ and as the Frobenius on $H^*_{\text{\'et}}(\mathbb P^1_{\overline{k}}, \mathbb Q_\ell)$.
	Then there is an inequality of $\mathsf L^1$-traces 
	$$\left|F, \mathrm{Sym} \left(H^*_{\text{\'et}}(\mathbb P^1_{\overline{k}}, \mathbb Q_\ell)[1] \otimes \bigoplus_{n \geq 1} (V[-1])^{\otimes n} \right)\right| \geq \left|F,  \bigoplus_{m \geq 0}  H_{\text{\'et}, c}^*(\mathrm{Conf}_m \mathbb P^1_{\overline{k}}, \mathbb Q_\ell) \otimes_{S_m} V^{\otimes m}\right|.$$
\end{prop}

Using the estimates of the previous subsection we get the following.

\begin{prop}\label{CombinedBound}
 Assume that our ground field $k$ is $\mathbb F_q$.  
	Let $V$ be a bigraded vector space with an action of $F$ that is concentrated in homological degrees $\geq 2$ (equivalently in cohomological degrees $\leq - 2$) and in grading degree $\geq 1$.  
	We also assume that $F$ acts on $H_{\text{\'et}, c}^*(\mathrm{Conf}_m \mathbb P^1_{\overline{k}}, \mathbb Q_\ell)$ as the Frobenius. Suppose that $|F, V_{\bullet, 2}| \leq D$.  If $|F, V_{n,i}| \leq E d^n b^i$ for some $E, d ,b > 0$  
	with $E \geq 1$
	and $b<1$, then $W:= \oplus_m H_{\text{\'et}, c}^*(\mathrm{Conf}_m \mathbb P^1_{\overline{k}}, \mathbb Q_\ell) \otimes_{S_m} V^{\otimes m}$  is concentrated in homological degrees $\geq 0$ and satisfies 
	$$|F, (W_{n,i})_{S_n} | = O\left(ni \max\left(2Ed + 8Edb^2q, D\right)^n (4b)^i\right). $$ 
\end{prop}
\begin{proof}
 	By  \cref{SchurBound} it suffices to bound the $\mathsf L^1$-trace on 
	$$\mathrm{Sym}(H_{\text{\'et}}^*(\mathbb P^1_{\overline{k}}, \mathbb Q_\ell)[1] \otimes \bigoplus_{n \geq 1} (V[-1])^{\otimes n}).$$ Applying \cref{tensoralgebrabound} to $V[-1]$ we obtain that 
	\[
	|F, (\bigoplus_{m\geq 1}  (V[-1])^{\otimes m})_{n,i} | \leq (2Eba)^n (2b)^i.
	\]  
	So $|F, (\bigoplus_{m \geq 1} H_{\text{\'et}}^*(\mathbb P^1_{\overline{k}}, \mathbb Q_\ell)[1] \otimes  (V[-1])^{\otimes m} )_{n, i} | \leq (\frac{1}{2b} +   2bq) (2Ebd)^n (2b)^i $. 
	Therefore applying \cref{symbound} to $\bigoplus_{m \geq 1} H_{\text{\'et}}^*(\mathbb P^1_{\overline{k}}, \mathbb Q_\ell)[1] \otimes  (V[-1])^{\otimes m}$  (observing that the homological degree zero component of this space is $V_{\bullet,2}$) we obtain the result.
\end{proof}

\subsection{Simplicial schemes and cohomology}

We next discuss simplicial schemes and their cohomology.  We follow closely the exposition of \cite[Section 2.2]{DH24}. (For a more classical reference, one may consult Brian Conrad's notes \cite{Conrad}.)

\begin{defi} 
Let $\Delta$ be the simplex category whose objects are non-empty finite subsets of $\mathbb N$ and whose morphisms are non-decreasing maps.  A simplicial scheme is a functor $\Delta^{\mathrm{op}} \to \mathrm{sSch}$ where $\mathrm{sSch}$ is the category of separated schemes.  A morphism between simplicial schemes is a natural transformation of functors.

For each $n \in \mathbb N$, we write $[n] = \{0, \dotsc, n\} \in \Delta$.  Then a simplicial scheme is determined by its values on the full subcategory $\{ [n] \, |\,  n \in \mathbb N \} \subset \Delta^{\mathrm{op}}$. 

Given a separated scheme $X$, the constant simplicial scheme $X_{\bullet}$ is defined by setting $X_\bullet([n]) = X$ for every $n \in \mathbb N$ and setting all maps to be $\mathrm{id}_X$.  
Throughout the paper, we frequently drop $\bullet$ and identify separated schemes and constant simplicial schemes.
\end{defi}

\begin{defi} 
Let $X$ be a separated scheme and $A$ be a simplicial scheme. Then an augmentation morphism from $A$ to $X$ is a morphism $\epsilon: A \to X$ in the category of simplicial schemes. \end{defi}

For any separated scheme 
$X$, we let $\mathrm{Sh}(X)$ denote the abelian category of $\mathbb Z_\ell$-module sheaves in the pro-\'etale topology as in \cite{BS15} and let $D^b(X)$ denote the associated bounded derived category.

Suppose that $A$ is a simplicial scheme; 
as mentioned earlier we can work with $A$ by restricting our attention to the full subcategory $\{ [n] \, |\,  n \in \mathbb{N} \} \subset \Delta^{\mathrm{op}}$. A $\mathbb Z_\ell$-sheaf on $A$ is a family of $\mathbb Z_\ell$-sheaves $\mathcal F_n$ on $A([n])$ in the pro-\'etale topology equipped with a compatible family of maps $A(i)^*\mathcal F_n \to \mathcal F_m$ for $i : [n] \to [m]$.
Let $\mathrm{Sh}(A)$ be the category of $\mathbb Z_\ell$-sheaves on $A$.  This is an abelian category with enough injective objects and thus the formalism of derived categories works in this settings. We denote the corresponding bounded derived category by $D^b(A)$. (For more details, see \cite[Section 2.2]{DH24}.)

Let $X$ be a separated scheme and let $A$ be a simplicial scheme with an augmentation $\epsilon : A \to X$. We next define the following functors:
\[
\epsilon_* : \mathrm{Sh}(A) \to \mathrm{Sh}(X), \quad \epsilon^* : \mathrm{Sh}(X) \to \mathrm{Sh}(A).
\]
First, for any sheaf $\mathcal G$ on $X$ we define $(\epsilon^*\mathcal G)_n = \epsilon_n^*\mathcal G$ with the natural identification $A(i)^*\epsilon_n^* \mathcal G \cong \epsilon^*_m\mathcal G$ coming from $\epsilon_m = \epsilon_n \circ A(i)$ where $i : [n] \to [m]$. The adjoint pushforward $\epsilon_*$ is given by $\epsilon_* \mathcal F= \mathrm{Eq}(\epsilon_{0 *}\mathcal F_0 \rightrightarrows \epsilon_{1*}\mathcal F_1)$ for any sheaf $\cF$ on $A$.  
These induce the derived functors:
\[
R\epsilon_* : D^b(A) \to D^b(X), \quad  \epsilon^* : D^b(X) \to D^b(A).
\]
\begin{rema}
There is a concrete description of $R \epsilon_*$: given a complex $\cI^\bullet$ of sheaves on $A$ such that $(\cI^i)_n$ is an injective sheaf on $A_n$ for every $i,n$, we may construct a representative of $R \epsilon_*\cI^\bullet$ as follows.   Take the chain complex of cosimplicial sheaves $[n] \mapsto \epsilon_{n*} \cI^\bullet$, form a bicomplex by applying the functor that takes a cosimplicial sheaf to its (normalized) chain complex, and then take the total complex.  (This construction works because we may factor $\epsilon_*$ as a pushforward to sheaves on the constant simplicial scheme on $X_\bullet$, followed by the functor from simplicial sheaves on $X_\bullet$ to sheaves on $X$ given by taking the limit over $\Delta$.  Because the pushforward preserves injectives, $R\epsilon_*$ factors as a composite of derived functors,  and the levelwise injectives are acyclic for the pushforward to $X_\bullet$, while the derived functor of the limit over $\Delta$ is modeled by the (normalized)  total complex.)
\end{rema}

\begin{rema}\label{Simplicial_Model}
	As a consequence of the previous remark, we note that if $X = \Spec ( \overline k)$ and $A_n$ is finite \'etale over  $X$ for all $n$ (so that $A$ corresponds to a simplicial set) then $R \epsilon_* \mathbb Z_l$ is canonically isomorphic to the normalized $\mathbb Z_l$ cochains on $A$ (considered as a simplicial set).
\end{rema}

Now suppose that $\pi : X\to \Spec \, k$ is a $k$-separated scheme where $k$ is a field.  We use the derived functor to define cohomology as
\[
H^i_{\text{\'et}}(A, \mathcal F) := H^i_{\text{\'et}}(R(\pi \circ \epsilon)_*\mathcal F).
\]
When the augmentation $\epsilon : A \to X$ is proper, i.e, for every $n \in \mathbb Z_{\geq 0}$ the morphism $\epsilon_n : A_n \to X$ is proper, we define cohomology with compact supports as
\[
 H^i_{\text{\'et}, c}(A, \mathcal F) := H^i_{\text{\'et}}(R\pi_! \circ R\epsilon_*\mathcal F).
\]
Finally when we have the constant sheaf $\mathcal F = \underline{\mathbb Z}_\ell$ or $\underline{\mathbb Q}_\ell$ on $X$, we have a natural morphism
\[
\mathcal F \to R\epsilon_*\epsilon^*\mathcal F.
\]
In particular, we have an induced homomorphism between cohomologies
\[
H^i_{\text{\'et}}(X, \mathbb Z_\ell) \to H^i_{\text{\'et}}(A, \mathbb Z_\ell). 
\]
When the augmentation $\epsilon : A \to X$ is proper, we also have an induced homomorphism
\[
H^i_{\text{\'et}, c}(X, \mathbb Z_\ell) \to H^i_{\text{\'et}, c}(A, \mathbb Z_\ell). 
\]

\subsection{Poschemes and the bar complex}

Following \cite[Section 2.2]{DT24}, we recall the notion of poschemes and the construction of the bar complex.

\begin{defi} 
Let $S$ be a scheme.
A poscheme over $S$ is a finite type morphism $P \to S$ with a closed poset relation $\leq_P \subset P\times_{S} P$, i.e., for any scheme $T \to S$ the pair $(P(T), \leq_P(T))$ is a poset. 
We denote the poscheme obtained by swapping the two coordinates of $\leq_P$ by $\geq_P$. Then we have $\Delta_{P/S} = \leq_P\cap \geq_P$. In particular, $\Delta_{P/S}$ is closed so that $P \to S$ is separated.
\end{defi}

 In this paper, we will usually work with poschemes $P \to S$ which are proper and will usually assume that $\Delta_{P/S}$ is also open in $\leq_P$ so that the relation $<_{P} := \leq_{P} \! \backslash \Delta_{P/S}$ is closed.  Whenever we discuss a finite poset $Q$, we will implicitly think of it as a finite poscheme over the ground field $k$ consisting entirely of $k$-points.

\begin{defi}
Let $P$ be a poscheme over $k$. The nerve $N(P)$ associated to $P$ is the simplicial scheme which assigns to $[n]$ the scheme
\[
N(P)([n]) = \leq_P \times_P \dotsb \times_P \leq_P = \{ p_0 \leq \dotsb \leq p_n \, | \, (p_i \leq p_{i + 1}) \in \leq_P \text{ for any $i = 0, \dotsc, n-1$} \},
\]
and which assignes to a non-decreasing map $i: [n] \to [m]$ the morphism
\[
N(P)([m]) \to N(P)([n]), \, (p_0 \leq \dotsb \leq p_m) \mapsto (p_{i(0)} \leq \dotsb \leq p_{i(n)}).
\]
\end{defi}

We will frequently work not just with poschemes $P$ but also with schemes $Z$ that admit stratifications indexed by a poscheme $P$.  The following definitions make this notion precise.

\begin{defi} 
Let $P$ be a poscheme over $S$. A $P$-space consists of an 
$S$-morphism $Z \to P$ with another $P$-morphism
\[
a : \leq_P \! \times_{\pi_2} Z \to Z,
\]
where $\pi_i : \leq_P \subset P\times_{S} P \to P$ denotes the projection map to the $i$th factor and $\leq_P \! \times_{\pi_2} Z $ is considered as a $P$-scheme via $\pi_1$. For any $S$-morphism $T \to P$, we denote the base change of $Z \to P$ over $T \to P$ by $Z_T$.
\end{defi}

\begin{defi} 
Let $X$ be a scheme of finite type over $S$ and $P$ be a poscheme over $S$. A stratification of $X$ by $P$ is a closed subscheme $Z \subset P\times_{S} X$ such that the morphism
\begin{align*}
\leq_P \! \times_{\pi_2} Z & \to P\times_{S} X, \\
(p \leq q, z \in Z_q) & \mapsto (p, \rho (z)),
\end{align*}
factors through $Z$ where $\rho : Z \subset P \times_{S} X \to X$ is the projection.  
In other words, for any $p, q \in P$ with $p \leq q$ we have $Z_q \subset Z_p$.  Note that $Z$ carries the structure of a $P$-space.
\end{defi}

The analogue of the nerve construction for a $P$-space is known as the bar complex.

\begin{defi}[{\cite[Definition 2.3]{DT24}}]
Let $P$ be a poscheme over $S$ such that $\Delta_P$ is open in $\leq_P$. 
Let $Z$ be a $P$-space.
The bar complex $B(P, Z)$   
is the simplicial scheme which assigns to $[n]$ the scheme
\[
B(P, Z)([n]) = \leq_P \! \times_P \dotsb \times_P \leq_P \! \times_P Z = \{ p_0 \leq \dotsb \leq p_n, z \in Z_{p_n}\}
\]
and which assigns to a non-decreasing map $i : [n] \to [m]$ the morphism 
\[
B(P, Z)([m]) \to B(P, Z)([n]), \, (p_0 \leq \dotsb \leq p_m, z) \mapsto (p_{i(0)}\leq \dotsb \leq p_{i(n)}, a(p_{i(n)} \leq p_{m},z)).
\]
\end{defi}

\subsubsection{Relative posets} \label{sect:relativeposets}

We will also need analogues of the previous constructions in a relative setting.

\begin{defi} 
\label{defi:relativebarcomplexes}
Let $H$ be a poscheme over $k$ such that $\Delta_H$ is open in $\leq_H$. Let us fix a locally closed subscheme $W \subset H$ and define the poscheme
\[
(W < H) = \{ (w, x) \in W \times H\, | \, w < x\},
\]
with the poset structure
\[
\leq_{W<H} = \{(w < x_1, w < x_2) \, | \, x_1 \leq x_2\}.
\]
This is a poscheme over $W$. 
\end{defi}

Let $X$ be a separated scheme of finite type over $W$ and $Z \subset (W < H) \times_W X$ be a stratification of $X$ by $(W < H)$.  
Let $P \subset (W < H)$ be a closed and downward closed subscheme.
Then we define the simplicial scheme $B(P, Z)$ to be the bar complex applied to $P$ and $Z$.   
In other words
\[
B(P, Z)([n]) = B(P, Z)([n]) = \{ w <  x_0 \leq \dotsb \leq x_n, z \in Z_{w < x_n} \}.
\]
This admits the augmentation morphism 
\[
\epsilon : B(P, Z)([n]) \to X, \, (w < x_0 \leq \dotsb \leq x_n, z) \mapsto (w, \rho(z) \in X_w).
\]

\subsection{Spectral sequences for the inclusion-exclusion principle for stratified schemes}
\label{subsec:spectralsequenceforstratifiedspaces}

The following exposition is based on \cite[Section 2.1]{DT24}.
Let $k$ be a perfect field and let $X$ be a separated scheme of finite type over $k$. 

Suppose that we have a finite filtration by closed subschemes: 
\[
Z_0 \subset  Z_1 \subset \dotsb \subset  Z_n = X.
\]
This induces a co-filtration of complexes
\[
\underline{\mathbb Z}_{\ell} \to i_{n-1!}\underline{\mathbb Z}_{\ell}|_{(Z_{n-1})_{\overline{k}}} \to \dotsb \to i_{0!}\underline{\mathbb Z}_{\ell}|_{(Z_0)_{\overline{k}}}
\]
and thus a spectral sequence with $E_1$-page
\[
E_1^{i, j} = H^{i + j}_{\text{\'et}, c}((U_i)_{\overline{k}}, \mathbb Z_\ell),
\]
converging to $H^{i+j}_{\text{\'et}, c}(X_{\overline{k}}, \mathbb Z_\ell)$ where $U_i = Z_i \setminus Z_{i-1}$.

We extend this spectral sequence to spaces stratified by posets.
Suppose that we have a finite poset $P$ and a collection of closed subschemes $Z_p \subset X$ for $p \in P$ such that $Z_p \supset Z_q$ whenever $p \leq q$.
We let $S_p = Z_p \setminus \cup_{p < q} Z_q$ which is a locally closed subscheme of $X$.
Since $P$ is finite, we may choose an injective homomorphism $h : P \to (\mathbb N, \leq)$.
We define $n$ to be $\max h(P)$ and set 
$K_i : = \cup_{h(p) \geq n-i} Z_p$ with the reduced scheme structure.  We obtain the associated filtration:
\[
\emptyset = K_{-1} \subset K_0 \subset \dotsb \subset K_n \subset X
\]

First suppose that for any $x \in X$ the set $\{p \in P \, | \, x \in Z_p\}$ has a maximal element. In particular this shows that $S_p \cap S_q = \emptyset$ whenever $p \neq q$.
This assumption implies that we have $K_{i} \setminus K_{i-1} = S_{h^{-1}(n-i)}$  
and there is a spectral sequence with $E_1$-page 
\[
E_1^{i, j} =  H^{i + j}_{\text{\'et}, c}((S_{h^{-1}(n-i)})_{\overline{k}}, \mathbb Z_\ell)
\]
converging to $H^{i + j}_{\text{\'et}, c}((K_n)_{\overline{k}}, \mathbb Z_\ell)$. 
Note that for each $n$ we have 
\[
\bigoplus_{i + j = n} E_1^{i, j} = \bigoplus_{p \in P} H^{n}_{\text{\'et}, c}((S_p)_{\overline{k}}, \mathbb Z_\ell)
\]
so that this direct sum is independent of the choice of $h$. However the differentials depend on the choice of $h$.

When the above assumption fails, we can define Zariski open subsets
\[
S'_{h^{-1}(n-i)} := K_i \setminus K_{i -1} \subset S_{h^{-1}(n-i)},
\]
and we still have a spectral sequence with $E_1$-page
\begin{equation} \label{eq:stratss}
E_1^{i, j} = H^{i + j}_{\text{\'et}, c}((S'_{h^{-1}(n-i)})_{\overline{k}}, \mathbb Z_\ell),
\end{equation}
converging to $H^{i + j}_{\text{\'et}, c}((K_n)_{\overline{k}}, \mathbb Z_\ell)$.
However, note that the definition of $S'_{h^{-1}(n-i)}$, and thus also the terms in the spectral sequence, depends on the choice of $h$.

Finally we introduce a version of this spectral sequence for bar complexes.
We work in the setting of \cref{sect:relativeposets} so that $X$ is a separated scheme of finite type over $W$ with a stratification $Z$ by $(W<H)$ with $P \subset (W < H)$.
We assume that $P$ is proper over $W$.
Suppose that our poscheme $P$ over $W$ admits a stratification by a finite poset $\mathfrak T = \{ T \in \mathfrak T\}$ via the collection of closed subschemes $\mathcal Z_T \subset P \to W$ such that $\mathcal Z_T$ is downward closed with respect to the poset structure over $W$ and $\mathcal Z_T \subset \mathcal Z_{T'}$ whenever $T\leq T'$.  We also assume that for any $(w < x) \in P$, the set $\{ T \in \mathfrak T \, | \, (w < x) \in \mathcal Z_T\}$ has a minimal element.  
Again we fix an injective homomorphism $h : \mathfrak T \to \mathbb N$.
We define 
\[
\mathcal Z_i = \cup_{h(T)\leq i} \mathcal Z_T.
\]
We denote the inclusion $\mathcal Z_i \subset P$ by $\iota_i$, and we consider the sub-bar complex
\[
\iota'_i : B_i = B(\mathcal Z_i, \iota_{i}^*Z) \hookrightarrow B(P, Z),
\]
which induces a stratification of closed simplicial schemes 
\[
\emptyset = B_{-1} \subset B_0 \subset \dotsb \subset B_n = B(P, Z)
\]
This induces a co-filtration:
\[
\underline{\mathbb Z}_{\ell} \to \iota'_{n-1!}\underline{\mathbb Z}_{\ell}|_{(B_{n-1})_{\overline{k}}} \to \dotsb \to \iota'_{0!}\underline{\mathbb Z}_{\ell}|_{(B_0)_{\overline{k}}} \to 0.
\]
We denote the induced filtration by
\[
0 = F^n\underline{\mathbb Z}_{\ell} \to F^{n-1}\underline{\mathbb Z}_{\ell} \to F^{n-2}\underline{\mathbb Z}_{\ell} \to \dotsb \to F^{-1}\underline{\mathbb Z}_{\ell} = \underline{\mathbb Z}_{\ell}.
\]
Then the graded part is given by
\[
\mathrm{Gr}_F^i\underline{\mathbb Z}_{\ell} = F^{i-1}\underline{\mathbb Z}_{\ell}/F^{i}\underline{\mathbb Z}_{\ell},
\]
and this gives rise to a spectral sequence
\[
E_1^{i, j} = H^{i + j}_{\text{\'et}, c}(B(P, Z)_{\overline{k}}, \mathrm{Gr}_F^i\underline{\mathbb Z}_{\ell} ) \implies H^{i + j}_{\text{\'et}, c}(B(P, Z)_{\overline{k}}, \mathbb Z_\ell).
\]

\section{Rational curves on del Pezzo surfaces}
\label{sec:ratondel}

In this section we collect several preliminary results about the nef cone and rational curves on degree $4$ del Pezzo surfaces.  We will focus on del Pezzo surfaces satisfying the following condition:

\begin{defi}
Let $S$ be a smooth del Pezzo surface of degree $\leq 7$ over a field $k$.  We say that $S$ is split if it satisfies any of the following equivalent conditions:
\begin{itemize}
\item the pullback map $\Pic(S) \to \Pic(S_{\overline{k}})$ is an isomorphism;
\item every $(-1)$-curve on $S_{\overline{k}}$ is defined over $k$;
\item the Picard ranks $\rho(S)$ and $\rho(S_{\overline{k}})$ are equal.
\end{itemize}
\end{defi}

\subsection{Nef cone}

Suppose that $S$ is a split del Pezzo surface of degree $4 \leq d \leq 7$.  Then:
\begin{itemize}
\item The pseudo-effective cone of $S$ is generated by classes of $(-1)$-curves.
\item The nef cone of $S$ is generated by fibers of conic fibrations $S \to \mathbb{P}^{1}$ and pullbacks of the hyperplane class under birational morphisms $S \to \mathbb{P}^{2}$.
\end{itemize}
We will impose a chamber decomposition on $\Nef_{1}(S)$ using the following structures which were used previously by \cite[Corollary 2.3]{Testa09}.

\begin{lemm}
Let $S$ be a split del Pezzo surface of degree $4 \leq d \leq 7$.  For any class $\alpha \in \Nef_{1}(X)_{\mathbb{Z}}$, there is a non-negative integer $b$, a contraction of a $(-1)$-curve $\psi: S \to Y$, and a nef class $\beta \in \Nef_{1}(Y)_{\mathbb{Z}}$ such that $\alpha = -bK_{S} + \psi^{*}\beta$.

Furthermore if $\beta$ is ample then the choice of $b,\psi,\beta$ is unique.
\end{lemm}

\begin{proof}
Every codimension $1$ extremal face of $\Nef_{1}(S)$ is the pullback of $\Nef_{1}(Y)$ for some contraction of a $(-1)$-curve $\psi: S \to Y$.  Consider the subdivision of $\Nef_{1}(S)$ where the top-dimensional subcones are generated by $-K_{S}$ and a codimension $1$ extremal face.  Every nef class is contained in some cone in this subdivision, yielding an expression of the form $\alpha = -bK_{S} + \psi^{*}\beta$ for some $b \in \mathbb{R}_{\geq 0}$ and $\beta \in \Nef_{1}(Y)$.  Since $b = \alpha \cdot E$ where $E$ is the contracted curve we see that $b \in \mathbb{Z}$, and hence $\beta$ is also an integral class.  Finally if $\beta$ is ample then there is only one full-dimensional cone in this subdivision whose closure contains $\alpha$.  
\end{proof}

Applying the lemma repeatedly, we obtain:

\begin{prop} \label{prop:nefexpression}
Let $S$ be a split del Pezzo surface of degree $4$ and let $\alpha \in \Nef_{1}(X)_{\mathbb{Z}}$.  Then there is a sequence of contractions of $(-1)$-curves
\begin{equation*}
S \to S_3 \to S_2 \to S_1
\end{equation*} 
with the following property: writing $\psi_{i}: S \to S_i$ for the induced map, $H$ for the unique pullback of a hyperplane class on $\mathbb{P}^{2}$ to $S_1$, and $F_{1},F_{2}$ for the two classes of conics on $S_1$, there are non-negative integers $b, b_{3},b_{2},x,y_{1},y_{2}$ with
\begin{equation} \label{eq:alphaexpression}
\alpha = -bK_{S} - b_{3}\psi_{3}^{*}K_{S_{3}} - b_{2} \psi_{2}^{*}K_{S_{2}} + \psi_{1}^{*}(xH + y_{1}F_{1} + y_{2}F_{2}).
\end{equation}
If all these integers (except possibly $b$) are positive, then there is a unique sequence of contractions of three disjoint $(-1)$-curves such that $\alpha$ admits such an expression.
\end{prop}

\begin{defi} \label{defi:regions}
Suppose $\vec{E} = (E_{1},E_{2},E_{3})$ is an ordered triple of disjoint $(-1)$-curves.  We define $\mathcal T_{\vec{E}} \subset \Nef_{1}(S)$ to be the subcone generated by all curve classes $\alpha \in \Nef_{1}(S)_{\mathbb{Z}}$ which admit an expression as in \cref{eq:alphaexpression} for the ordered sequence of contractions of $E_{1},E_{2},E_{3}$.
\end{defi}

The arrangement of $(-1)$-curves on a split degree $4$ del Pezzo surface are described by the $5$-regular Clebsch graph so there are 960 
choices of ordered triples of disjoint $(-1)$-curves.  (The Weyl group $D_{5}$ has order $1920$; its action on this set of triples has stabilizer of order $2$ corresponding to identifying one of the two conic fibrations $S_{1} \to \mathbb{P}^{1}$.)  
\cref{prop:nefexpression} shows that as we vary $\vec{E}$ over all such choices, the union of the subcones $\mathcal T_{\vec{E}}$ is all of $\Nef_{1}(S)$ and the intersections of their interiors are empty.  We will use this subdivision of $\Nef_{1}(S)$ in later sections.

\begin{nota} \label{nota:bignota}
Fix $\alpha \in \Nef_{1}(S)$.  By an abuse of notation, we denote by $\mathcal T_{\alpha}$ any region as in \cref{defi:regions} that contains $\alpha$.  (We will only use this notation when the choice of $\mathcal T_{\alpha}$ is irrelevant.)

Retaining the notation of \cref{prop:nefexpression}, we define $\rho_{\alpha}: S \to \mathbb{P}^{1} \times \mathbb{P}^{1}$ to be the composition of $\psi_{1}$ with the contraction of the unique $(-1)$-curve yielding a map to $\mathbb{P}^{1} \times \mathbb{P}^{1}$.  Labeling the strict transform of this curve on $S$ as $E_{4}$, we write $(E_{1},E_{2},E_{3},E_{4})$ for the corresponding tuple of disjoint $(-1)$-curves on $S$.

We define $p_{\alpha}: S \to \mathbb{P}^{1}$ to be the composition of $\rho_{\alpha}$ with the projection to the first factor.  We let $F_{\alpha}$ denote a general fiber of $p_{\alpha}$.  Finally we define the non-negative integers
\begin{gather*}
h(\alpha) = -K_{S} \cdot \alpha \qquad \qquad a(\alpha)  = F_{\alpha} \cdot \alpha \qquad \qquad k_{i} = E_{i} \cdot \alpha \\
a'(\alpha)  = \frac{1}{2} \left( h(\alpha) - 2a(\alpha) + \sum_{i=1}^{4} k_{i}(\alpha) \right)
\end{gather*}
Note that $a'(\alpha)$ is the intersection of $\alpha$ against a fiber $F'_{\alpha}$ of the composition of $\rho_{\alpha}$ with the second projection map $\mathbb{P}^{1} \times \mathbb{P}^{1} \to \mathbb{P}^{1}$.  We will sometimes surpress the $\alpha$ from $\rho,p,F,h,a,a',k_{i}$ in situations where it can be understood.  We write $\mathbf{k} = (k_{1},k_{2},k_{3},k_{4})$.  In the notation of \cref{eq:alphaexpression} we have: 
\begin{gather*}
\begin{aligned}
-K_{S} \cdot \alpha & = 4b + 5b_{3} + 6b_{2} + 3x + 2y_{1} + 2y_{2} \\
F_{\alpha} \cdot \alpha &= 2b + 2b_{3} + 2b_{2} + x + y_{1}
\end{aligned}\\
E_{1} \cdot \alpha = b \qquad E_{2} \cdot \alpha = b + b_{3} \qquad E_{3} \cdot \alpha = b  + b_{3} + b_{2} \qquad E_{4} \cdot \alpha = b + b_{3} + b_{2} + x
\end{gather*}
yielding the inequalities
\begin{equation*}
h(\alpha) \geq 2a(\alpha) \geq \sum_{i=1}^{4} k_{i}(\alpha).
\end{equation*}
If $\alpha$ lies in the interior of $\mathcal{T}_{\alpha}$ then the inequalities are all strict.  The same inequalities hold for $a'(\alpha)$ by symmetry.
\end{nota}

\subsection{Rational curves}

We next address the geometry of rational curves on $S$.  Since maps to $(-1)$-curves are ignored both by Manin's Conjecture and by the Cohen-Jones-Segal conjecture, we will concentrate on rational curves whose numerical classes are nef.  Given a class $\alpha \in N_{1}(S)_{\mathbb{Z}}$, we denote by $\Mor(\mathbb{P}^{1},S)_{\alpha}$ the space of morphisms $s: \mathbb{P}^{1} \to S$ such that $s_{*}[\mathbb{P}^{1}] = \alpha$.

\begin{theo}[{\cite{Testa09}, \cite{BLRT21}}]
\label{theo: Testa}
Let $S$ be a del Pezzo surface of degree $4$ over a field $k$.  For every $\alpha \in \mathrm{Nef}_1(S)_{\mathbb{Z}}$ the moduli space $\Mor(\mathbb{P}^{1},S)_{\alpha}$ is geometrically irreducible of the expected dimension and generically parametrizes free curves.
\end{theo}

\begin{proof}
It suffices to prove the statement when $k$ is algebraically closed.  \cite[Corollary 5.2]{Testa09} addresses characteristic $0$ and \cite[Theorem 1.5]{BLRT21} addresses positive characteristic.
\end{proof}

\begin{nota}
Given a nef curve class $\alpha \in \Nef_{1}(S)_{\mathbb{Z}}$, we will denote by $M_{\alpha}$ the scheme $\Mor(\mathbb{P}^{1},S)_{\alpha}$ as in \cref{theo: Testa}.  
\end{nota}

We will verify in \cref{lemm:psismooth} that $M_{\alpha}$ is smooth over every field.

\section{Curves and conic bundle structures}
\label{sec:curvesandconicbundle}

Let $S$ be a split degree $4$ del Pezzo surface over a field $k$.  Suppose we fix a nef class $\alpha \in \Nef_{1}(S)_{\mathbb{Z}}$ along with the associated constructions in \cref{nota:bignota}.  Note that $p: S \to \mathbb{P}^{1}$ induces a morphism $p_{*}: \Mor(\mathbb{P}^{1},S) \to \Mor(\mathbb{P}^{1},\mathbb{P}^{1})$ sending $M_{\alpha}$ to the irreducible component $N_{a}$ that parametrizes morphisms of degree $a(\alpha) := F \cdot \alpha$.  In this section we will construct a factorization
\begin{equation*}
\xymatrix{
M_{\alpha} \ar[rr]^{p_{*}} \ar[dr]_{\psi_{\alpha}} & & N_{a} \\
& R_{\alpha} \ar[ur]_{\pi_{\alpha}} &
}
\end{equation*}
where $\psi_{\alpha}$ is smooth surjective with connected fibers and $\pi_{\alpha}$ is quasi-finite.  The key result in this section is that $M_{\alpha}$ is an open subset of a Zariski-locally trivial projective bundle over $R_{\alpha}$ and in turn $R_{\alpha}$ is an open subset of a Zariski-locally trivial projective bundle over an open subset of a product of projective spaces.  In particular, $M_{\alpha}$ is smooth.

\subsection{Compactification of $R_{\alpha}$} \label{sect:compactza}

In this subsection we work over a ground field different than $\mathbb{F}_{2}$.   
We fix an integer $a \geq 1$ and a ordered tuple of non-negative integers $\mathbf{k} = (k_{1},k_{2},k_{3},k_{4})$ with $k_i \leq a$ for any $i = 1, \dotsc, 4$.  

To improve clarity of exposition, we continue to write $B := \mathbb{P}^{1}$ for the domain of our maps.  The space $N_{a} = \Mor(B,\mathbb{P}^{1})_a$ 
parametrizing degree $a$ morphisms is (by definition) the open sublocus of the linear system $|\mathcal{O}(a,1)| $ on $B \times \mathbb{P}^{1}$ that parametrizes sections.  We will identify $\mathbb{P}^{2a+1}$ with the linear system $|\mathcal{O}(a,1)|  = \mathbb {P}(H^0(B \times \mathbb P^1, \cO(a,1)))$.  We fix four $k$-points $p_{1},p_{2},p_{3},p_{4} \in \mathbb{P}^{1}$ and let $B_{1},B_{2},B_{3},B_{4}$ denote the corresponding fibers of the second projection map for $B \times \mathbb{P}^{1}$.  (The $p_{i}$ will eventually represent images of the singular fibers of $S \to \mathbb{P}^{1}$.)

\begin{defi}
We define the closed subscheme
\begin{equation*}
Z_{a, \mathbf{k}} \subset \mathbb{P}^{2a+1} \times \mathrm{Hilb}^{[k_1]}(B_1) \times \mathrm{Hilb}^{[k_2]}(B_2) \times \mathrm{Hilb}^{[k_3]}(B_3) \times \mathrm{Hilb}^{[k_4]}(B_4)
\end{equation*}
as the incidence correspondence
\[
Z_{a, \mathbf{k}} = \{([C], [T_{1}], [T_{2}], [T_{3}], [T_{4}])  | \, T_{i} \subset C \cap B_{i} \text{ for $i = 1, \dotsc, 4$}\}
\]
where we have identified $\mathbb{P}^{2a+1} \cong |\mathcal{O}(1,a)|$ on $B \times \mathbb{P}^{1}$.  We call $Z_{a, \mathbf{k}}$ a \textit{configuration cover}.

We let $Z_{a, \mathbf{k}}^{\circ}$ denote the preimage of $N_{a} \subset \mathbb{P}^{2a+1}$ under the projection map. This is the space of degree $a$ maps $s : B \to \mathbb P^1$ together with divisors $T_1, \dotsc, T_4 \subset B$ such that $T_i \subset s^*(p_i)$.

\end{defi}

Our next goal is to study the properties of $Z_{a, \mathbf{k}}$ and $Z_{a, \mathbf{k}}^{\circ}$.

\begin{prop}
\label{prop:locallycompleteintersectionZ}
The space $Z_{a, \mathbf{k}}  \subset \mathbb P^{2a+1} \times \prod_{i = 1}^{4}\mathrm{Hilb}^{[k_i]}(B_i)$ is the zero locus of a section of a globally generated vector bundle of rank $\sum_{i = 1}^4k_i$. 
\end{prop}

\begin{proof}
For simplicity of notation we will write $Y = \mathbb P^{2a+1} \times \prod_{i = 1}^{r}\mathrm{Hilb}^{[k_i]}(B_i)$.
Let $\mathcal C$ be the universal family over $\mathbb{P}^{2a+1}$ parametrizing divisors of type $(a, 1)$ in $B \times \mathbb{P}^{1}$.  If we define $\mathcal L = \mathcal O_{B \times \mathbb{P}^{1} \times \mathbb P^{2a+1}}(a, 1, 1)$, then $\mathcal{C}$ is the vanishing locus of a section $s \in H^0(B \times \mathbb{P}^{1} \times \mathbb P^{2a+1}, \mathcal L)$.
We denote the projection $B \times \mathbb{P}^{1} \times Y \to B \times \mathbb{P}^{1} \times \mathbb P^{2a+1}$ by $q$.

Next let
\[
\mathcal U_i' \subset B_i \times \mathrm{Hilb}^{k_i}(B_i)
\]
be the universal family of the Hilbert scheme and let $\mathcal{U}_{i}$ be its base change to $B_i \times Y$.
We consider the projection $\rho_i : \mathcal U_i \to Y$.
Since this is a flat finite morphism of degree $k_i$ we conclude that $\mathcal E_i = (\rho_i)_*(q^*\mathcal L|_{\mathcal U_i})$ is a locally free sheaf of rank $k_i$.
Then $q^*s$ defines a global section of the locally free sheaf $\mathcal E = \oplus_{i = 1}^r \mathcal E_i$ whose zero locus is equal to $Z_{a, \mathbf{k}}$.

It only remains to show that $\mathcal{E}$ is globally generated.  It suffices to prove that each $\mathcal{E}_{i}$ is globally generated.  Note the fiber of $\mathcal{E}_{i}$ over a point $y \in Y$ is $H^{0}(T,\mathcal{O}(a))$ where $T \subset B_{i}$ is the length $k_{i}$ subscheme defined by the point in $\Hilb^{k_{i}}(B_{i})$ corresponding to $y$.
Thus the composition $H^{0}(Y,\mathcal{E}_{i}) \otimes \mathcal{O}_{Y} \to \mathcal{E}_{i} \to \mathcal{E}_{i}|_{y}$ can be identified with  
\[
H^0(B \times \mathbb{P}^{1} \times \mathbb{P}^{2a+1}, \mathcal L) = H^{0}(B \times \mathbb{P}^{1} \times Y,q^*\mathcal L) \to H^{0}(\mathcal U_{i},q^*\mathcal L|_{\mathcal U_i}) \to H^0(T, \mathcal O(a))
\]
which is the restriction map on global sections for $T \subset q(B_{i} \times \{ y \}) \subset B \times \mathbb{P}^{1} \times \mathbb{P}^{2a+1}$.  Since $k_{i} \leq a$ we see that this map is surjective, showing global generation of $\mathcal{E}_{i}$.
\end{proof}

\begin{coro} \label{coro:zacircequidim}
The space $Z_{a, \mathbf{k}}^{\circ}$ is equidimensional of dimension $2a+1$.
\end{coro}

\begin{proof}
\Cref{prop:locallycompleteintersectionZ} shows that every irreducible component of $Z_{a, \mathbf{k}}$ has dimension $\geq 2a+1$, and so the same is true for the open set $Z_{a, \mathbf{k}}^{\circ}$.  On the other hand, since the projection morphism $Z_{a, \mathbf{k}}^{\circ} \to \mathbb{P}^{2a+1}$ is quasi-finite we see that every irreducible component of $Z_{a, \mathbf{k}}^{\circ}$ has dimension $\leq 2a+1$.  
\end{proof}

We next study the restriction of the projection morphism:
\[
\phi_{a, \mathbf{k}} : Z_{a, \mathbf{k}} \to  \prod_{i = 1}^4 \mathrm{Hilb}^{[k_i]}(B_i).
\]

Note that for any $i = 1, \dotsc, 4$, $B_i$ is naturally identified with $B$ under the projection map. Using this identification we define the following Zariski open subset of $\prod_{i = 1}^{4} \mathrm{Hilb}^{[k_i]}(B_i)$:
\[
U_{\mathbf{k}} := \left\{ \left. [T_i] \in \prod_{i = 1}^{4} \mathrm{Hilb}^{[k_i]}(B_i) \,  \right| \, \text{for any $i \neq j$, the supports of $T_i$ and $T_j$ are disjoint} \right\}.
\]

\begin{lemm}
\label{lemm:configurationcover_fibration}
Assume that $2a \geq \sum_i k_i$.  Then every fiber of $\phi_{a,\mathbf{k}}: \phi_{a,\mathbf{k}}^{-1}(U_{\mathbf{k}}) \to U_{\mathbf{k}}$ is a projective space of dimension $\geq 2a + 1 - \sum_i k_i$.  Moreover we have
\[
\phi_{a, \mathbf{k}}(Z^{\circ}_{a, \mathbf{k}}) \subset U_{\mathbf{k}}.
\]
as a Zariski open subset, and the fibration
\[
\phi_{a, \mathbf{k}} : Z^{\circ}_{a, \mathbf{k}} \to \phi_{a, \mathbf{k}}(Z^{\circ}_{a, \mathbf{k}}) 
\]
is an open subset of a Zariski-locally trivial $\mathbb P^{2a + 1 - \sum_i k_i}$-bundle over $\phi_{a, \mathbf{k}}(Z^{\circ}_{a, \mathbf{k}})$.
\end{lemm}

\begin{proof}
For any $([T_i]) \in U_{\mathbf{k}}$ the linear subspace
\[
\{ [C] \in \mathbb P^{2a+1} \, | \, \text{for every $i = 1, \dotsc, 4$, $C\cap B_i$ contains $T_i$ as a closed subscheme.} \}
\]
has dimension $\geq 2a + 1 - \sum_{i = 1}^r k_i >0$, showing the first claim.  Next, for any irreducible curve $C$ in the linear series $|\mathcal{O}(1, a)|$ the supports of the $C \cap B_{i}$'s cannot share any point so that the image of $[C]$ lies in $U_{\mathbf{k}}$.  We claim that the morphism
\[
\phi_{a, \mathbf{k}}|_{Z_{a, \mathbf{k}}^\circ} : Z_{a, \mathbf{k}}^\circ \to U_{\mathbf{k}},
\]
is smooth.  Since $Z_{a, \mathbf{k}}^{\circ}$ is equidimensional by \cref{coro:zacircequidim}, smoothness follows from \cref{lemm:smoothnesscriterion} once we show that the fibers are smooth of the expected dimension.
Consider the exact sequence
\[
0 \to N_{C/B \times \mathbb P^1}(-T_1 - \dotsb -T_4) \to N_{C/B \times \mathbb P^1} \to \oplus_{i = 1}^4N_{T_i/B} \to 0.
\]
We have $N_{C/B \times \mathbb P^1}\cong \mathcal O_{C}(2a)$.  Since $2a - \sum k_{i} \geq 0$ we conclude that $h^{0}(C,N_{C/B \times \mathbb P^1}(-T_1 - \dotsb -T_4)) = 2a + 1 - \sum_{i = 1}^r k_i$ and $h^{1}(C,N_{C/B \times \mathbb P^1}(-T_1 - \dotsb -T_4)) = 0$.  In particular the obstruction space vanishes so that the moduli space of sections is smooth at $C$ of dimension $h^{0}$.  This concludes the proof of smoothness of $\phi_{a, \mathbf{k}}|_{Z^{\circ}_{a, \mathbf{k}}}$.  
Since $\phi_{a, \mathbf{k}}|_{Z^{\circ}_{a, \mathbf{k}}}$ is flat, $\phi_{a, \mathbf{k}}(Z^{\circ}_{a, \mathbf{k}})$ is an open subset of $U_{\mathbf{k}}$.   Then the natural morphism from $Z^{\circ}_{a, \mathbf{k}}$ to $\mathbb{P}^{2a+1} \times \phi_{a, \mathbf{k}}(Z^{\circ}_{a, \mathbf{k}})$ realizes each fiber as an open subset of a linear subspace of $\mathbb{P}^{2a+1}$.  This is also true of the generic fiber, thus $\phi_{a,\mathbf{k}}^{-1}( \phi_{a,\mathbf{k}}(Z^{\circ}_{a, \mathbf{k}}))$ is a Zariski-locally trivial projective bundle over $\phi_{a, \mathbf{k}}(Z^{\circ}_{a, \mathbf{k}})$.
\end{proof}

As an immediate consequence, we obtain: 

\begin{coro} \label{coro: Z_irreducible}
Suppose that $2a \geq \sum_{i} k_{i}$.  Then $Z_{a, \mathbf{k}}^{\circ}$ is smooth and geometrically irreducible of dimension $2a+1$.
\end{coro}

\subsection{Configuration covers and the Stein factorization}

We now return to our original setting: we would like to consider the Stein factorization of (a compactification of) $p_{*}:\Mor(B, S)_\alpha = M_{\alpha} \to N_{a} = \Mor(B, \mathbb P^1)_a$.  Recall that \cref{nota:bignota} describes the notation involving the nef curve class $\alpha$ and the constants $h,a,\mathbf{k}$.

Since $S$ is a split del Pezzo surface, the conic bundle $p: S \to \mathbb{P}^{1}$ will have four reducible fibers.  For $i=1,\dotsc,4$ we let $p_{i} \in \mathbb{P}^{1}$ denote the $k$-point which is the $p$-image of the exceptional divisor $E_{i}$.  We construct $Z_{a, \mathbf{k}}$ with respect to these four points as in the previous subsection.   We then define the morphism $\psi_{\alpha}: M_{\alpha} \to Z_{a, \mathbf{k}}^{\circ}$ by sending
\begin{equation*}
[s: B \to S] \mapsto ( [p \circ s], [s^{*}(E_1)], \dotsc, [s^{*}(E_4)] )
\end{equation*}
Note that we have a factoring $p_{*} = \pi \circ \psi_{\alpha}$ where $\pi: Z_{a, \mathbf{k}}^{\circ} \to \mathbb{P}^{2a+1}$ is the projection map.  In this subsection our goal is to study the morphism $\psi_{\alpha}$.

\begin{lemm} \label{lemm:psismooth}
The morphism $\psi_{\alpha}: M_{\alpha} \to Z_{a, \mathbf{k}}^{\circ}$ is smooth.
\end{lemm}

\begin{proof}

By \cref{lemm:smoothnesscriterion} it suffices to show that the fibers of $\psi_{\alpha}$ are smooth of the expected dimension.
We first analyze the fibers of $M_{\alpha} \to U_{\mathbf{k}}$.  Each fiber can be interpreted as the sublocus of $\Mor(B,S)_{\alpha}$ that sends a designated subscheme $T_{i}$ to $E_{i}$ for every $1 \leq i \leq 4$.  Thus the tangent space and obstruction space of the fiber at a point $s: B \to S$ are determined by the cohomology groups of $s^{*}T_{S}(- \log(E_{1} + E_{2} +E_{3} + E_{4}))$.  Using the birational map $\rho: S \to \mathbb{P}^{1} \times \mathbb{P}^{1}$ of \cref{nota:bignota} we can identify
\begin{align*}
s^{*}T_{S}(- \log(E_{1} + E_{2} +E_{3} + E_{4})) & \cong s^{*}\rho^{*}T_{\mathbb{P}^{1} \times \mathbb{P}^{1}}(-E_{1}-E_{2}-E_{3}-E_{4}) \\
& \cong \mathcal{O}_{B}(2a - \sum k_{i}) \oplus \mathcal{O}_{B}(2a' - \sum k_{i}).
\end{align*}
By the inequalities of \cref{nota:bignota} we see both direct summands have vanishing $H^1$. 
We conclude that $M_{\alpha} \to U_{\mathbf{k}}$ is smooth by \cref{lemm:smoothnesscriterion}.

The tangent spaces of the fibers of $Z^{\circ}_{a, \mathbf{k}} \to U_{\mathbf{k}}$ were computed in the proof of \cref{lemm:configurationcover_fibration}.  By comparing this lemma against the direct summand decomposition above, we conclude that the deformation theory of points in the fibers of $\psi_{\alpha}$ is controlled by the kernel of the map
\begin{equation*}
s^{*}T_{S}(- \log(E_{1} + E_{2} +E_{3} + E_{4})) \to s^{*}p^{*}T_{\mathbb{P}^{1}}(- \sum T_{i} ).
\end{equation*}
This map is projection onto the first factor in the direct summand decomposition above.  The kernel $\mathcal{O}_{B}(2a' - \sum k_{i})$ is a nef line bundle.  Again we are in a situation where obstructions vanish and thus can conclude smoothness by \cref{lemm:smoothnesscriterion}.
\end{proof}

\begin{defi}
We let $R_{\alpha} \subset Z_{a, \mathbf{k}}^{\circ}$ denote the open subset which is the image of the flat morphism $M_{\alpha} \to Z_{a, \mathbf{k}}^{\circ}$.  
We let $\pi_{\alpha}: R_{\alpha} \to N_{a}$ denote the restriction of the projection map $\pi: Z_{a, \mathbf{k}}^{\circ} \to N_{a} \subset \mathbb{P}^{2a+1}$.
\end{defi}


Before continuing we need an auxiliary construction.

\begin{const} \label{cons:qalpha}
We have a morphism $p'_{*}: \Mor(B,S) \to \Mor(B,\mathbb{P}^{1})$ given by composing with the birational map $\rho: S \to \mathbb{P}^{1} \times \mathbb{P}^{1}$ and the second projection map $\mathbb{P}^{1} \times \mathbb{P}^{1} \to \mathbb{P}^{1}$.  We construct $Z_{a', \mathbf{k}}^{\circ}$ as in \cref{sect:compactza}; note that all results in this section hold equally well for $Z^{\circ}_{a', \mathbf{k}}$.  In particular, defining $Q_{\alpha} = R_{\alpha} \times_{U_{\mathbf{k}}} Z^{\circ}_{a', \mathbf{k}}$ we obtain an induced morphism
\begin{equation*}
i: M_{\alpha} \to Q_{\alpha}.
\end{equation*}
\end{const}

\begin{lemm} \label{lemm:malphaopenimmer}
The morphism $i$ is an open immersion realizing $M_{\alpha}$ as an open subset of a Zariski-locally trivial $\mathbb{P}^{h - 2a +1}$ bundle over $R_{\alpha}$.
\end{lemm}

\begin{proof}
Since both $M_{\alpha}$ and $Q_{\alpha}$ are smooth and have the same dimension, it suffices to show that $i$ is \'etale and set-theoretically injective.  We show that it is \'etale by computing the induced map on tangent spaces.  Recall from the proof of \cref{lemm:psismooth} that the relative tangent space of the map $M_{\alpha} \to U_{\mathbf{k}}$ is
\begin{equation*}
H^{0}(B,s^{*}T_{S}(-\log(E_{1}+E_{2}+E_{3}+E_{4}))) \cong s^{*}p^{*}T_{\mathbb{P}^{1}}(- \sum T_{i} ) \oplus s^{*}p'^{*}T_{\mathbb{P}^{1}}(- \sum T_{i} ).
\end{equation*}
In turn these two summands are the relative tangent spaces of $R_{\alpha}$ and $Z_{a', \mathbf{k}}^{\circ}$ over $U_{\mathbf{k}}$.  This shows that the induced map on tangent spaces is an isomorphism, so $i$ is \'etale.  Injectivity follows from the fact that $\rho_{*}: \Mor(B,S) \to \Mor(B,\mathbb{P}^{1} \times \mathbb{P}^{1})$ is a bijection on the set of morphisms which meet the open locus where $\rho$ is an isomorphism.
\cref{lemm:surjectivityphi} shows that $Z'_{a', \mathbf{k}}$ is an open subset of a Zariski-locally trivial projective bundle over $U_{\mathbf{k}}$, and thus by base change $M_{\alpha}$ has the analogous property over $R_{\alpha}$.
\end{proof}

Putting everything together, we obtain:

\begin{theo}
The morphism $p_{*}: M_{\alpha} \to N_{a}$ factors as $p_{*} = \pi_{\alpha} \circ \psi_{\alpha}$ where $\psi_{\alpha}: M_{\alpha} \to R_{\alpha}$ realizes $M_{\alpha}$ as an open subset of a Zariski-locally trivial $\mathbb{P}^{h - 2a +1}$ bundle and $\pi_{\alpha}$ is quasifinite.  In particular, $M_{\alpha}$ is smooth.
\end{theo}

\subsection{Dimension bounds}
Consider the chain of morphisms $M_{\alpha} \to Z_{a, \mathbf{k}}^{\circ} \to U_{\mathbf{k}}$.  Our final goal is to estimate the dimension of the complement of the image of $M_{\alpha}$.  We will do this in several steps; the first is to analyze the image of $Z_{a, \mathbf{k}}^{\circ} \to U_{\mathbf{k}}$.

\begin{lemm}
\label{lemm:surjectivityphi}  
Consider the sublocus of $U_{\mathbf{k}}$ parametrizing closed subschemes in $\sqcup_{i} B_{i} \subset B \times \mathbb{P}^{1}$ which are not contained in any irreducible section parametrized by $|\mathcal{O}(a, 1)|$.  This locus is empty if $\sum_{i = 1}^4 k_i \leq a$ 
and has dimension at most $2\sum k_i - (2a+1)$ otherwise.

In particular, the complement of $\phi_{a, \mathbf{k}}(Z_{a, \mathbf{k}}^{\circ})$ in $U_{\mathbf{k}}$ is empty if $\sum_{i = 1}^4 k_i \leq a$ and has dimension at most $2\sum k_i - (2a+1)$ otherwise.
\end{lemm}

\begin{proof}
By \cref{lemm:configurationcover_fibration} every fiber of $\phi_{a, \mathbf{k}}$ over a point of $U_{\mathbf{k}}$ is a projective space of dimension $\geq 2a + 1 -\sum_i k_i >0$ and in particular is non-empty.  Suppose $([T_i]) \in U_{\mathbf{k}}$ is a point which is not in the image of $Z_{a, \mathbf{k}}^{\circ}$.  Thus the general divisor contained in the fiber $\phi_{a, \mathbf{k}}^{-1}([T_i])$ will admit some $\pi_{1}$-vertical components.  We let $\ell_{1}$ (respectively $\ell_{2}$) denote the intersection of $\pi_{2}^{*}\mathcal{O}(1)$ with the union of all irreducible $\pi_{1}$-vertical components of a general divisor in the fiber which do not (respectively do) intersect some $T_{i}$.  Then the dimension of $\phi_{a, \mathbf{k}}^{-1}([T_i])$ is
\[
\max \left\{0, 2a +1 -2\ell_1 -\ell_2 - \sum_i k_i \right\} + \ell_1.
\]
This needs to be greater than or equal to $2a + 1 - \sum_i k_i$.  This is only possible when $2a +1 -2\ell_1 -\ell_2 - \sum_i k_i \leq 0$ so that the dimension is exactly $\ell_{1}$.  In particular
\[
\ell_1 \geq 2a + 1 - \sum_i k_i.
\]
On the other hand, for any $i = 1, \dotsc, 4$, by taking intersection of our divisors against $B_{i}$ we see that $a-\ell_1 \geq k_i$.
Altogether we have
\[
a-\max_i \{ k_i\} \geq \ell_1 \geq 2a + 1 - \sum_i k_i.
\]
This implies that $ \sum_i k_i - \max_i \{ k_i\} \geq a + 1$, proving the claim about emptiness.

The dimension of the locus in $\mathbb{P}^{2a+1}$ parametrizing divisors such that the 
$\pi_{1}$-vertical components have degree $\ell_1 + \ell_{2}$ is
\[
2a - \ell_1 - \ell_2+1.
\]
Thus, the sublocus of $\phi_{a, \mathbf{k}}^{-1}(U_{\mathbf{k}})$ corresponding to such divisors is also bounded above by this number.  Since the argument above showed that the fiber over $([T_i])$ has dimension $\ell_1$, we have
\[
\dim_{[T_i]}(U_{\mathbf{k}} \setminus \phi_{a, \mathbf{k}}(Z_{a, \mathbf{k}}^{\circ})) \leq 2a - 2\ell_1 - \ell_2+1.
\]
Combining with the inequality
\[
\ell_1 \geq 2a+1 - \sum_i k_i,
\]
we obtain
\[
\dim(U_{\mathbf{k}} \setminus \phi_{a, \mathbf{k}}(Z_{a, \mathbf{k}}^{\circ})) \leq 2\sum_i k_i -2a -1. \qedhere
\]
\end{proof}

We next turn to the morphism $\psi_{\alpha}: M_{\alpha} \to Z_{a, \mathbf{k}}^{\circ}$.

\begin{prop}
\label{prop:surjectivity}
Consider the map $\psi_{\alpha}: M_{\alpha} \to Z_{a, \mathbf{k}}^{\circ}$.  The complement of $\psi_{\alpha}(M_{\alpha})$ is empty if $h-2a - \sum k_i \geq 0$ and has dimension at most $\sum k_{i} - (h - 2a + 1)$ otherwise.
\end{prop}

\begin{proof}
Choose a point $([C'], [T_{1}], [T_{2}], [T_{3}], [T_{4}]) \in Z_{a, \mathbf{k}}^{\circ}$ and define $Y_{C'} = C' \times_{\mathbb{P}^{1},p} S$.  Then $\rho_{\alpha}$ induces a birational morphism $\phi_{C'} : Y_{C'} \to C' \times \mathbb{P}^{1}$.  For every morphism $B \to S$ in the fiber of $M_{\alpha} \to Z_{a, \mathbf{k}}^{\circ}$ the graph can be identified with a section $C_{1}$ of $C' \times \mathbb{P}^{1}$ of relative anticanonical height $2a' = h-2a-\sum k_{i}$ such that $T_{i} \subset C_{1} \cap '(C' \times \{ p_{i} \})$.

Conversely, suppose we have a section $C_1$ of relative anticanonical height $2a'$ on $C' \times \mathbb{P}^{1}$ whose intersection with $C' \times \{ p_i\}$ contains $T_i$ for every $i$.  Then its strict transform on $Y_{C'}$ gives a section which in turn can be viewed as the graph of an element of $M_{\alpha}$.  Thus the dimension of the complement of $\psi_{\alpha}(M_{\alpha})$ is bounded by the dimension of the locus of $([T_{i}])$ which are not contained in an irreducible section of $C' \times \mathbb{P}^{1}$ of relative anticanonical height $2a'$.  \Cref{lemm:surjectivityphi} shows that this locus is empty if $\sum k_{i} \leq a'$ and has dimension at most $2 \sum k_{i} - (2a'+1)$ otherwise.  Using $a' = \frac{1}{2}(h - 2a +\sum k_{i})$ we obtain the statement.
\end{proof}

By combining \cref{lemm:surjectivityphi,prop:surjectivity} we obtain:

\begin{prop} \label{prop:finaldimestimate}
Consider the composition
\[
 \phi_{a, \mathbf{k}} \circ \psi_\alpha: M_{\alpha} \to Z_{a, \mathbf{k}}^{\circ} \to U_{\mathbf{k}}.
\]
Then the dimension of $U_{\mathbf{k}(\alpha)} \setminus  \phi_{a, \mathbf{k}} \circ \psi_\alpha(M_{\alpha})$ is at most 
\[
\max\left\{2 \sum_i k_i - (2a +1),2 \sum_i k_i  -(2a' +1) \right\}.
\]
\end{prop}

\subsection{Pointed morphisms} \label{sect:pointedmorzak}

In this section we reformulate the previous results in the setting of pointed morphisms.  Since the arguments are essentially the same, we will only outline the proofs.

In this section we let $B = \mathbb P^1$ denote our domain curve.  We fix base points $*_B$ and $*_{\mathbb P^1}$ on $B$ and $\mathbb P^1$ respectively.  For any $a \geq 0$ we let $\mathbb P^{2a}_{*}$ denote the linear system of divisors of type $(a, 1)$ in $B\times \mathbb P^1$ passing through $(*_B, *_{\mathbb P^1})$.  This is isomorphic to a hyperplane $\mathbb{P}^{2a}$ in $|(a,1)|$.
Suppose $p_1, \dotsc, p_{4} \in \mathbb P^1$ are distinct $k$-points on $\mathbb P^1$ which are different from $*_{\mathbb P^1}$.
Let $B_i = B\times \{p_i\}  \subset B \times\mathbb P^1 $.

For a fixed index $1 \leq i \leq 4$ and a fixed integer $0\leq k_i \leq a$ we define 
\[
Z_{a, i, k_i, *} = \{([C], [T]) \in \mathbb P^{2a}_* \times \mathrm{Hilb}^{[k_i]}(B_i)\, | \, T \subset C\} \subset \mathbb P^{2a}_* \times \mathrm{Hilb}^{[k_i]}(B_i).
\]
We define the configuration cover as
\[
Z_{a, \mathbf{k}, *} =  Z_{a, 1, k_1, *} \times_{\mathbb P^{2a}_*} \dotsb \times_{\mathbb P^{2a}_*} Z_{a, 4, k_4, *}\subset \mathbb P^{2a}_* \times \prod_{i = 1}^{4}\mathrm{Hilb}^{[k_i]}(B_i) 
\]
equipped with the projection $\pi_{a, \mathbf{k}, *}: Z_{a, \mathbf{k}, *} \to \mathbb P^{2a}_*$.

We denote the open locus of $\mathbb P^{2a}_{*}$ parametrizing honest sections by $N_{a, *}$ and we denote its preimage in $Z_{a, \mathbf{k}, *}$ by $Z_{a, \mathbf{k}, *}^\circ$.
We first analyze the structure of $Z_{a, \mathbf{k}, *}^\circ$
using the fibration
\[
\phi_{a, \mathbf{k}, *} : Z^\circ_{a, \mathbf{k}, *} \to  U_{\mathbf{k}, *}
\]
where
\[
U_{\mathbf{k}, *} := \left\{ \left. [T_i] \in \prod_{i = 1}^{4} \mathrm{Hilb}^{[k_i]}(B_i) \,  \right| \begin{array}{c} \text{for every $i$, the support of $T_i$ does not contain } *_{B} \textrm{ and}\\ \text{for every $i \neq j$, the supports of $T_i$ and $T_j$ are disjoint} \end{array} \right\}.
\]

\begin{lemm}
\label{lemm:configurationcover_fibrationpointed}
Assume that $2a > \sum_i k_i$. 
Then every fiber of $\phi_{a,\mathbf{k}, *}: \phi_{a,\mathbf{k}, *}^{-1}(U_{\mathbf{k}, *}) \to U_{\mathbf{k}, *}$ is a projective space of dimension $\geq 2a  - \sum_i k_i$.  Moreover we have
\[
\phi_{a, \mathbf{k}, *}(Z^{\circ}_{a, \mathbf{k}, *}) \subset U_{\mathbf{k}, *}.
\]
as a Zariski open subset, and the fibration
\[
\phi_{a, \mathbf{k}, *} : Z^{\circ}_{a, \mathbf{k}, *} \to \phi_{a, \mathbf{k}}(Z^{\circ}_{a, \mathbf{k}, *}) 
\]
is an open subset of a Zariski-locally trivial $\mathbb P^{2a - \sum_i k_i}$-bundle over $\phi_{a, \mathbf{k}, *}(Z^{\circ}_{a, \mathbf{k}, *})$.

In particular $Z^{\circ}_{a, \mathbf{k}, *}$ is smooth and geometrically irreducible.
\end{lemm}

We next return to the surface $S$ and the morphism $p_{*}: M_{\alpha, *} \to N_{a, *}$.  As before, we construct $Z_{a,\mathbf{k}, *}$ with respect to the four $k$-points which are the images in $\mathbb P^1$ of the singular fibers.   We then define the morphism $\psi_{\alpha, *}: M_{\alpha, *} \to Z_{a,\mathbf{k}, *}^{\circ}$ by sending
\begin{equation*}
[s: B \to S] \mapsto ( [p \circ s], [s^{*}(E_1)], \dotsc, [s^{*}(E_4)] ).
\end{equation*}
If we assume that $2a > \sum k_{i}$ and $2a' > \sum k_{i}$, then the same argument as in \cref{lemm:psismooth} shows that $\psi_{\alpha,*}$ is smooth.  Under these assumptions we let $R_{\alpha,*} \subset Z_{a,\mathbf{k}, *}^{\circ}$ denote its image and let $\pi_{\alpha,*}: R_{\alpha, *} \to N_{a, *}$ denote the projection map.

\begin{theo}
Assume that $2a > \sum k_{i}$ and $2a' > \sum k_{i}$.  The morphism $p_{*}: M_{\alpha, *} \to N_{a, *}$ factors as $p_{*} = \pi_{\alpha, *} \circ \psi_{\alpha, *}$ where $\psi_{\alpha, *}: M_{\alpha, *} \to R_{\alpha, *}$ realizes $M_{\alpha, *}$ as an open subset of a Zariski-locally trivial $\mathbb{P}^{h - 2a}$-bundle over the Zariski open subset $R_{\alpha, *}$ and $\pi_{\alpha, *}$ is quasifinite.  In particular, $M_{\alpha, *}$ is smooth.
\end{theo}

Finally, we turn to dimension estimates.

\begin{lemm}
\label{lemm:surjectivityphipointed}  

Assume that $2a > \sum k_{i}$.
Consider the sublocus of $U_{\mathbf{k}, *}$ parametrizing closed subschemes in $\sqcup_{i} B_{i} \subset B \times \mathbb P^1$ which are not contained in any irreducible section parametrized by $|\mathcal{O}(a, 1)|$ passing through the base point.  This locus is empty if $\sum_{i = 1}^4 k_i < a$  and has dimension at most $2\sum k_i - 2a$ otherwise.

\end{lemm}

By applying \cref{lemm:surjectivityphipointed} twice -- once to the map $\phi_{a, \mathbf{k}, *}: \phi_{a, \mathbf{k}, *}^{-1}(U_{\mathbf{k}, *}) \to U_{\mathbf{k}, *}$ and once to the map  $\psi_{\alpha, *}: M_{\alpha, *} \to Z_{a, \mathbf{k}, *}^{\circ}$ -- we obtain:

\begin{prop}
\label{prop:dimensionestimatepointed}
Assume that $2a > \sum k_{i}$ and $2a' > \sum k_{i}$.  Consider the composition
\[
\phi_{a, \mathbf{k}, *}\circ  \psi_{\alpha, *} : M_{\alpha, *} \to Z_{a, \mathbf{k}, *}^{\circ} \to U_{\mathbf{k}, *}.
\]
Then the dimension of $U_{\mathbf{k}(\alpha), *} \setminus \phi_{a, \mathbf{k}, *} \circ \psi_{\alpha, *}(M_{\alpha, *})$ is at most 
\[
\max\left\{ 2\sum_i k_i - 2a, 2\sum_i k_i  - 2a' \right\}.
\]
\end{prop}

\section{Stratifying the Hilbert scheme}
\label{sec:stratifying}

Suppose $B$ is a smooth projective curve.  Then $\mathrm{Hilb}(B)$ carries the structure of a poscheme where the relation $\leq$ corresponds to containment of closed subschemes.  \cite{DT24} studied the space $\mathrm{Hilb}(B)^{Q}$ of poset homomorphisms $Q \to \mathrm{Hilb}(B)$ and introduced a combinatorial stratification of this space.  In this section we recall their constructions and rephrase them in the setting of poschemes. The main new addition is the key example \S\ref{sect:keyexample}.   We will closely follow the exposition of  \cite[Section 3]{DT24}.

\subsection{Poschemes of divisors}

Let $k$ be a perfect field and let $B$ be a smooth curve defined over $k$.
We denote the poscheme of all closed subschemes of $B$ by $\overline{\mathrm{Hilb}}(B)$. Note that this poscheme includes a maximal element which is a $k$-point corresponding to the entire scheme $B$ as well as a minimal element which is a $k$-point corresponding to the empty set.
We denote the poscheme of finite length subschemes of $B$ by $\mathrm{Hilb}(B)$ (i.e.~$\overline{\mathrm{Hilb}}(B)$ is obtained from $\mathrm{Hilb}(B)$ by adding in the maximal element $B$). 

\begin{defi}[{\cite[Definition 3.1]{DT24}}]
Let $Q$ be a finite \'etale poscheme over $k$ with a largest element $\hat{1}$ which is a $k$-rational point.  We set
\[
\widetilde{Q} = Q \setminus\{ \hat{1}\}.
\]
We denote the poscheme of homomorphisms $Q \to \overline{\mathrm{Hilb}}(B)$ which take $\hat{1}$ to $B$ by 
\[
\overline{\mathrm{Hilb}}(B)^Q.
\]
This means that for any $k$-scheme $T$, an element $x \in \overline{\mathrm{Hilb}}(B)^Q(T)$ consists of a collection of subschemes 
\[
\{ x_q \subset B\times T \}_{q \in Q(T)},
\]
such that each $x_{q}$ is a flat subscheme over $T$ for $q \in \widetilde{Q}(T)$, $x_{\hat{1}} = B \times T$, and $p \leq q$ implies $x_p \subset x_q$.
Note that one can realize $\overline{\mathrm{Hilb}}(B)^Q$ as a closed subscheme
\[
\overline{\mathrm{Hilb}}(B)^Q \subset \mathrm{Mor}(Q, \overline{\mathrm{Hilb}}(B)).
\]
We denote the subposcheme of $\overline{\mathrm{Hilb}}(B)^Q$ consisting of homomorphisms such that the preimage of $B$ is $\hat{1}$ by $\mathrm{Hilb}(B)^Q$. 
\end{defi}

Suppose we fix an element $x \in \mathrm{Hilb}(B)^Q(\overline{k})$. Then for each $c \in B(\overline{k})$, one can define a homomorphism of posets $\widetilde{g}_c : \widetilde{Q}(\overline{k}) \to \mathbb N$ by assigning to each $q \in \widetilde{Q}(\overline{k})$ the length of $x_q$ at $c$. In this way one can define a one-to-one correspondence between elements of $\mathrm{Hilb}(B)^Q(\overline{k})$ and the set of finitely supported functions:
\[
B(\overline{k}) \to \mathrm{Hom}( \widetilde{Q}(\overline{k}), \mathbb N).
\]
We say that a homomorphism $\widetilde{g}_c: \widetilde{Q}(\overline{k}) \to \mathbb N$ is trivial if $\widetilde{g}_c(q) = 0$ for every $q \in \widetilde{Q}(\overline{k})$. The support of a function $\widetilde{g} : B(\overline{k}) \to \mathrm{Hom}( \widetilde{Q}(\overline{k}), \mathbb N)$ is the set of $c \in B(\overline{k})$ such that $\widetilde{g}_c$ is non-trivial.
When $x \in \overline{\mathrm{Hilb}}(B)^Q(\overline{k})$, we can equivalently consider this element as a function
\[
g : B(\overline{k}) \to \mathrm{Hom}( Q(\overline{k}), \mathbb N \cup \{\infty\}),
\]
by assigning $g_c(q) = \infty$ whenever $x_q = B$.

\subsubsection{A key example} \label{sect:keyexample}

We describe a key example for our analysis.  The geometric motivation for this example will be explained in \cref{subsec:blowups}.

Let $V_i$ $(i = 1, 2)$ be a $2$-dimensional vector space over $k$ and set
\[
V = V_1 \oplus V_2.
\]
We fix a non-negative integer $r$ and  for each $j = 1, 2$ we pick $r$ distinct $1$-dimensional subspaces $\ell_{i, j} \subset V_j (i = 1, \dotsc, r)$ which are defined over $k$.
Let $Q$ be the following set of subspaces of $V$: 
\begin{equation*}
\{ V, V_1 \oplus \{0\}, \{0\} \oplus V_2,  0 \} \cup \{ \ell_{i, 1} \oplus \ell_{i, 2}, \ell_{i, 1} \oplus \{0\}, \{0\} \oplus \ell_{i, 2}  \}_{i=1,\dotsc,r}.
\end{equation*}

Under inclusion this becomes a finite \'etale poscheme that admits meets, i.e., the intersection of two subspaces in $Q$ is still contained in $Q$.

\subsection{Essential saturated elements}

\subsubsection{Saturated elements}

We closely follow the discussion of \cite[Section 3.3]{DT24}.
First we introduce the following definition from \cite[Section 3.3]{DT24}:

\begin{defi}[{\cite[Definition 3.3]{DT24}}]  \label{defi:saturated}
Let $Q$ be a finite \'etale poscheme over $k$ and assume that every subset $S \subset Q(\overline{k})$ admits a meet $\wedge_{s \in S}s$. Let $B$ be a smooth projective curve over $k$. We say an element $\{x_q\}_{q \in Q(\overline{k})} \in \mathrm{Hilb}(B)^Q(\overline{k})$ is saturated if for any subset $S \subset Q(\overline{k})$, the natural containment yields an equality
\[
x_{\wedge_{s \in S} s} = \bigcap_{s \in S} x_s.
\]
In other words, the associated homomorphism $x : Q(\overline{k}) \to \mathrm{Hilb}(B)(\overline{k})$ preserves meets. Let $Q^{JB}(\overline{k}) \subset \mathrm{Hilb}(B)^Q(\overline{k})$ be the subposet of saturated elements, which is a Zariski open subset of $\mathrm{Hilb}(B)^Q(\overline{k})$.  
Then since $Q^{JB}(\overline{k})$ is Galois invariant, this open subscheme descends to the open subscheme $Q^{JB}\subset \mathrm{Hilb}(B)^Q$ defined over $k$.
\end{defi}

In terms of finitely supported functions $\widetilde{g}: B(\overline{k}) \to \mathrm{Hom}(\widetilde{Q}(\overline{k}), \mathbb N)$, being saturated means that for any $c \in B(\overline{k})$, $\widetilde{g}_c : \widetilde{Q}(\overline{k}) \to \mathbb N$ preserves meets. We write $g_c : Q(\overline{k}) \to \mathbb N\cup \{\infty\}$ for the extension of $\widetilde{g}_c : \widetilde{Q}(\overline{k}) \to \mathbb N$ by assigning $g_c(\widehat{1}) = \infty$. When $g_c$ preserves meets, its left adjoint $f_c : \mathbb N \cup \{\infty\} \to Q(\overline{k})$ is given by
\[
f_c(n) = \min\{ q\in Q(\overline{k}) \, |\, g_c(q) \geq n \}.
\]
The following definition describes the properties which characterize $f_{c}$.

\begin{defi}
Let $Q$ be a finite \'etale poscheme with meets and with a maximal element $\widehat{1}$.   A chain is a poset homomorphism $f : \mathbb N \cup \{\infty\} \to Q(\overline{k})$ with the following properties:
\begin{enumerate}
\item $f(0) = \widehat{0}$ is the minimal element of $Q$, and
\item $f$ preserves joins, i.e., for any $S \subset \mathbb{N} \cup \{\infty\}$ we have
\[
f(\vee S) = \vee f(S).
\]
This is equivalent to requiring that $f(n) = \widehat{1}$ for $n \gg 0$.
\end{enumerate}
We denote the set of chains by $\mathrm{ch}(Q(\overline{k}))$.  By the trivial chain, we will mean the chain with $f(0) = \widehat{0}$ and $f(n) = \widehat{1}$ for $n \geq 1$.  The trivial chain is denoted by $V$ or $\widehat{1}$.
\end{defi}

Given a $\widetilde{g}_c : \widetilde{Q}(\overline{k}) \to N$ that preserves meets, the corresponding $f_{c}$ constructed above is a chain.  Conversely, for any chain $f_{c}:   \mathbb N\cup \{\infty\} \to Q(\overline{k})$ its right adjoint $g_{c} : Q(\overline{k}) \to \mathbb N\cup \{\infty\}$ defined by
\[
g_{c}(q) = \max\{n \in \mathbb N \cup \{\infty\} \,|\, q \geq f_{c}(n)\}
\]
will have the property that $g_{c}^{-1}( \{ \infty \}) = \{ \widehat{1} \}$.
In this way, there is a one-to-one correspondence between the set of chains $f_{c}: \mathbb N \cup \{\infty\} \to Q(\overline{k})$ 
and the set of homomorphisms $\widetilde{g}_{c} : \widetilde{Q}(\overline{k}) \to \mathbb N$ preserving meets.

Next we recall the following definition from \cite[Section 3.3]{DT24}:

\begin{defi}[{\cite[Definition 3.8]{DT24}}]
The saturation function is the function
\[
\mathrm{sat} : \mathrm{Hilb}(B)^Q(\overline{k}) \to Q^{JB}(\overline{k}),
\]
which is the left adjoint to the inclusion functor of the subposet of saturated elements.
In other words, for any $D \in \mathrm{Hilb}(B)^Q(\overline{k})$, $\mathrm{sat}(D)$ is defined as the smallest saturated element greater than or equal to $D$. Note that this saturation function is Galois-equivariant.
In particular if $D$ is a $k$-rational point on $\mathrm{Hilb}(B)^Q$, then $\mathrm{sat}(D)$ is also a $k$-rational point on $Q^{JB}$.
\end{defi}

\begin{exam}
\label{exam:saturatedexample}
Let $Q$ be the finite \'etale poscheme defined in \cref{sect:keyexample} so that we have 
\begin{equation*}
\{ V, V_1 \oplus \{0\}, \{0\} \oplus V_2,  0 \} \cup \{ \ell_{i, 1} \oplus \ell_{i, 2}, \ell_{i, 1} \oplus \{0\}, \{0\} \oplus \ell_{i, 2}  \}_{i=1,\dotsc,r}.
\end{equation*}
In this case, an element $D$ of $\mathrm{Hilb}(B)^Q(\overline{k})$ consists of divisors 
\[
D_{V_1}, D_{V_2}, D_{\ell_{1, 1}\oplus \ell_{1, 2}}, \dotsc, D_{\ell_{r, 1}\oplus \ell_{r, 2}}, D_{\ell_{1, 1}}, \dotsc, D_{\ell_{r, 1}}, D_{\ell_{1, 2}}, \dotsc, D_{\ell_{r, 2}}, D_0.
\]
This element is saturated if and only if we have 
\[
D_0 = D_{V_1}\cap D_{V_2} = D_{\ell_{i, 1}\oplus \ell_{i, 2}} \cap D_{\ell_{j, 1}\oplus \ell_{j, 2}}  
\]
for any $i \neq j$ and
\[
D_{\ell_{i, j}} = D_{\ell_{i, 1}\oplus \ell_{i, 2}}\cap D_{V_j}
\]
for any $i, j$.
For any $D \in \mathrm{Hilb}(B)^Q(\overline{k})$, $D' = \mathrm{sat}(D)$ is given by
\begin{align*}
D_0' &= (D_{V_1}\cap D_{V_2}) \cup \bigcup_{i \neq j}(D_{\ell_{i, 1}\oplus \ell_{i, 2}} \cap D_{\ell_{j, 1}\oplus \ell_{j, 2}}),\\
D_{\ell_{i, j}}' &= D_0' \cup (D_{\ell_{i, 1}\oplus \ell_{i, 2}}\cap D_{V_j}),\\
D_{\ell_{i, 1}\oplus \ell_{i, 2}}' &= D_{\ell_{i, 1}\oplus \ell_{i, 2}} \cup D_0' \cup D_{\ell_{i, 1}}' \cup D_{\ell_{i, 2}}',\\
D_{V_j}' &= D_{V_j} \cup D_0' \cup \bigcup_{i = 1}^r D_{\ell_{i, j}}'.
\end{align*}

For $c \in B(\overline{k})$, the associated function $\widetilde{g}_c : \widetilde{Q}(\overline{k}) \to \mathbb N$ consists of values
\[
\widetilde{g}_c(V_1), \, \widetilde{g}_c(V_2), \, \widetilde{g}_c(\ell_{i, 1}\oplus \ell_{i, 2}), \, \widetilde{g}_c(\ell_{i, j}), \, \widetilde{g}_c(0),
\]
for $i = 1, \dotsc, r$ and $j = 1, 2$.
The function $\widetilde{g}_c : \widetilde{Q}(\overline{k}) \to \mathbb N$ is saturated if and only if:
\begin{itemize}
\item $\widetilde{g}_c(\ell_{i, j})$ is equal to the minimum of $\widetilde{g}_c(\ell_{i, 1}\oplus \ell_{i, 2})$ and $\widetilde{g}_c(V_j)$, and
\item $\widetilde{g}_c(0)$ is equal to the minimum of $\widetilde{g}_c(V_1)$ and $\widetilde{g}_c(V_2)$ and also equal to $\widetilde{g}_c(\ell_{i, 1}\oplus \ell_{i, 2})$ for every index $i$ but one (for which the value may be larger).
\end{itemize}

Let $\widetilde{g}_c : \widetilde{Q}(\overline{k}) \to \mathbb N$ be a saturated function.
For simplicity, let us assume that $\widetilde{g}_c(V_1) \geq \widetilde{g}_c(V_2)$ and $\widetilde{g}_c(\ell_{1, 1}\oplus \ell_{1, 2})$ is the maximum among $\widetilde{g}_c(\ell_{i, 1}\oplus \ell_{i, 2})$'s.
First let us assume that $\widetilde{g}_c(\ell_{1, 1}) = \widetilde{g}_c(V_1)$, i.e., $\widetilde{g}_c(V_1) \leq \widetilde{g}_c(\ell_{1, 1}\oplus \ell_{1, 2})$. Let $m_0 = \widetilde{g}_c(0)$ and $m_{\ell_{1, 1}}= \widetilde{g}_c(\ell_{1, 1}) - \widetilde{g}_c(0)$ and $m_{\ell_{1, 1}\oplus \ell_{1, 2}} = \widetilde{g}_c(\ell_{1, 1}\oplus \ell_{1, 2}) - \widetilde{g}_c(\ell_{1, 1})$. Then $f_c : \mathbb N\cup\{ \infty\} \to Q(\overline{k})$ is given by
\[
f_c(n) = 
\begin{cases}
0 & \text{ if } n \leq m_0\\
\ell_{1, 1} & \text{ if } m_0 < n \leq m_0 + m_{\ell_{1, 1}} \\
\ell_{1, 1}\oplus \ell_{1, 2} & \text{ if } m_0 + m_{\ell_{1, 1}} < n \leq m_0 + m_{\ell_{1, 1}} + m_{\ell_{1, 1}\oplus \ell_{1, 2}}\\
V & \text{ if } m_0 + m_{\ell_{1, 1}} + m_{\ell_{1, 1}\oplus \ell_{1, 2}} < n.
 \end{cases}
\]
We write this chain as $m_0[0] + m_{\ell_{1, 1}}[\ell_{1, 1}] + m_{\ell_{1, 1}\oplus \ell_{1, 2}}[\ell_{1, 1}\oplus \ell_{1, 2}]$. When we have $\widetilde{g}_c(\ell_{1, 1}) = \widetilde{g}_c(\ell_{1, 1}\oplus \ell_{1, 2})$, $f_c$ takes the form 
\[
f_c(n) = 
\begin{cases}
0 & \text{ if } n \leq m_0\\
\ell_{1, 1} & \text{ if } m_0 < n \leq m_0 + m_{\ell_{1, 1}} \\
V_1 & \text{ if } m_0 + m_{\ell_{1, 1}} < n \leq m_0 + m_{\ell_{1, 1}} + m_{V_1}\\
V & \text{ if } m_0 + m_{\ell_{1, 1}} + m_{V_1} < n,
 \end{cases}
\]
where $m_{V_1} = \widetilde{g}_c(V_1) -  \widetilde{g}_c(\ell_{1, 1})$. We write this chain as $m_0[0] + m_{\ell_{1, 1}}[\ell_{1, 1}] + m_{V_1}[V_1]$.
\end{exam}

\subsubsection{Essential saturated elements}

We next define a poset structure on the set of chains $\mathrm{ch}(Q(\overline{k}))$.  We emphasize that our ordering is opposite to the poset structure naturally induced by $Q(\overline{k})$; this choice guarantees that the ordering of chains is compatible with the ordering of saturated elements of $\mathrm{Hilb}(B)^Q(\overline{k})$ so that there is no confusion when passing between the two settings.

\begin{defi}[{\cite[Section 3.4]{DT24}}] 
Let $Q$ be a finite \'etale poscheme with $\widehat{1}$ and assume that every subset of $Q(\overline{k})$ admits a meet.
For any $f, f' \in \mathrm{ch}(Q(\overline{k}))$, we say $f \leq f'$ holds if for any $n \in \mathbb N \cup\{\infty\}$, we have $f'(n) \leq f(n)$ as elements in $Q(\overline{k})$.

Since $Q(\overline{k})$ has meets, $\mathrm{ch}(Q(\overline{k}))$ has joins, i.e., for any $f, f' \in \mathrm{ch}(Q(\overline{k}))$, the join $f\vee f'$ of $f$ and $f'$ is given by
\[
(f\vee f')(n) = f(n)\wedge f'(n),
\]
for any $n \in \mathbb N \cup \{\infty\}$.
We have a canonical injective homomorphism
\begin{align*}
Q(\overline{k})^{\mathrm{op}} & \hookrightarrow \mathrm{ch}(Q(\overline{k})) \\
q & \mapsto (0 \leq q \leq \widehat{1}\leq \widehat{1} \leq \dotsb)
\end{align*}
\end{defi}

Now we introduce the notion of essential saturated pairs which plays a crucial role in our analysis:

\begin{defi}[{\cite[Definition 3.11]{DT24}}] \label{defi:essentialtype}
Let $Q$ be a finite \'etale poscheme with $\widehat{1}$ and assume that every subset of $Q(\overline{k})$ admits a meet.
We also assume that the natural inclusion $Q(k) \hookrightarrow Q(\overline{k})$ is a bijection.

Let $f_0 \in \mathrm{ch}(Q(\overline{k}))$ be a chain. Let $\mathfrak S$ be the set of $f_s \in \mathrm{ch}(Q(\overline{k}))$ such that $f_0 \prec f_s$, i.e., $f_0 < f_s$ and there is no $f'\in \mathrm{ch}(Q(\overline{k}))$ satisfying $f_0 < f' < f_s$.

We say a pair of chains $f_0 \leq f$ is essential if $f$ can be obtained from $f_0$ by taking the join of some subset of $\mathfrak S$. 

Let $B$ be a smooth projective curve. 
Suppose that $W \subset Q^{JB}$ is a locally closed subset with the reduced scheme structure.
We say that $(w\leq x) \in (W\leq Q^{JB})(\overline{k})$ is essential if for every $c \in \mathrm{Supp}(w\leq x)(\overline{k})$ we have that $f_c^w \leq f_c^x$ is essential.
\end{defi}

In our setting, this notion can be made explicit:

\begin{exam}
\label{exam:essentialsaturatedexample}
Let $Q$ be the poscheme defined in \cref{sect:keyexample} and considered in \cref{exam:saturatedexample}.  We will describe the essential saturated pairs $f_{0} \leq f$ where $f_{0} = m_{\ell_{i, 1}\oplus \ell_{i, 2}}[\ell_{i, 1}\oplus \ell_{i, 2}]$.

First let us list pairs such that $f_0 \prec f$: 
\begin{align*}
m_{\ell_{i, 1}\oplus \ell_{i, 2}}[\ell_{i, 1}\oplus \ell_{i, 2}] &\prec (m_{\ell_{i, 1}\oplus \ell_{i, 2}} +1)[\ell_{i, 1}\oplus \ell_{i, 2}] \\
m_{\ell_{i, 1}\oplus \ell_{i, 2}}[\ell_{i, 1}\oplus \ell_{i, 2}] &\prec [\ell_{i, j}]  + (m_{\ell_{i, 1}\oplus \ell_{i, 2}}-1)[\ell_{i, 1}\oplus \ell_{i, 2}] \qquad (j=1,2) 
\end{align*}
Then the list of pairs $f_0 \leq f$ in $\mathrm{ch}(Q)(\overline{k})$ such that $f_0 \leq f$ is essential is the above list plus the following pairs:
\begin{align*}
m_{\ell_{i, 1}\oplus \ell_{i, 2}}[\ell_{i, 1}\oplus \ell_{i, 2}] &=m_{\ell_{i, 1}\oplus \ell_{i, 2}}[\ell_{i, 1}\oplus \ell_{i, 2}]  \\
m_{\ell_{i, 1}\oplus \ell_{i, 2}}[\ell_{i, 1}\oplus \ell_{i, 2}] &< [\ell_{i, j}]  +m_{\ell_{i, 1}\oplus \ell_{i, 2}}[\ell_{i, 1}\oplus \ell_{i, 2}] \qquad (j=1,2) \\
m_{\ell_{i, 1}\oplus \ell_{i, 2}}[\ell_{i, 1}\oplus \ell_{i, 2}] &< [0] + (m_{\ell_{i, 1}\oplus \ell_{i, 2}}-1)[\ell_{i, 1}\oplus \ell_{i, 2}]\\
m_{\ell_{i, 1}\oplus \ell_{i, 2}}[\ell_{i, 1}\oplus \ell_{i, 2}] &< [0] + m_{\ell_{i, 1}\oplus \ell_{i, 2}}[\ell_{i, 1}\oplus \ell_{i, 2}]
\end{align*}
One can perform a similar analysis for other chains $f_{0}$.  The other relevant example for us will be the trivial chain $f_{0} = V = (0\leq V \leq V \leq \dotsb)$; in this case the list of saturated $f$ such that $f_{0} \leq f$ is essential is: 
\begin{equation*}
V, \qquad [V_{j}],  \qquad  [\ell_{i,1} \oplus \ell_{i,2}], \qquad   [\ell_{i,j}], \qquad [0],
\end{equation*}
as we vary over $j = 1,2$ and $i=1,\dotsc,4$.
\end{exam}

The following proposition is the main motivation of introducing this notion in \cite{DT24}:

\begin{prop}[{\cite[Proposition 3.13]{DT24}}]
\label{prop:homotopyequivalence}
Let $Q$ be a finite \'etale poscheme with $\widehat{1}$ and assume that every subset of $Q(\overline{k})$ admits a meet.

For any $f_0 \in \mathrm{ch}(Q(\overline{k}))$, essential pairs in $(f_0 < \mathrm{ch}(Q(\overline{k})))$ form an initial subposet of the poset $(f_0 < \mathrm{ch}(Q(\overline{k})))$. In particular the inclusion from the nerve of the subposcheme of essential pairs to the nerve $N(f_0 < \mathrm{ch}(Q(\overline{k})))$ of the poscheme $(f_0 < \mathrm{ch}(Q(\overline{k})))$ induces isomorphisms of $\ell$-adic cohomologies.
\end{prop}

\begin{proof}
For any $f_0 < f \in \mathrm{ch}(Q(\overline{k}))$, the set of essential elements $f_{0} < g$ satisfying $g \leq f$ admits a maximal element $\mathrm{ess}(f)$.  Indeed, one can define $\mathrm{ess}(f)$ to be the join of all $f_{s}$ such that $f_{0} \prec f_{s} \leq f$.  
Because $\mathrm{ess}(-)$ is right adjoint to the inclusion of essential elements, the natural transformation from the identity functor to  $\mathrm{ess}(-)$ induces a deformation retraction from the nerve of the finite \'etale poscheme $(f_0 < \mathrm{ch}(Q(\overline{k})))$ to the nerve of the poscheme of essential pairs.  (Alternately, see \cite[Lemma 3.2.1]{DH24}, and note this retract is proper in the sense of \cite[Section 2.1.3]{DH24} and so  \cite[Lemma 2.2.2]{DH24} applies).
\end{proof}

\subsection{Combinatorial functions}
 Here we follow the exposition of \cite[Section 3.5]{DT24}.

\begin{defi}
Let $B$ be a smooth projective curve.
Let $Q$ be a finite \'etale poscheme with meets and $\widehat{1}$. Let $h : Q(\overline{k}) \to \mathbb N$ be a function satisfying $h(\widehat{1}) = 0$.  Then we extend this function to $h : \mathrm{ch}(Q(\overline{k})) \to \mathbb N$ by
\[
h(f) = \sum_{n \geq 1} h(f(n)).
\]
Note that this is a finite sum due to the fact that $h(\widehat{1}) = 0$.

Next let $x \in Q^{JB}(\overline{k})$ be a saturated element which corresponds to a finitely supported function
\[
B(\overline{k}) \to  \mathrm{ch}(Q(\overline{k})), \, c  \mapsto f_c.
\]
We extend the function $h : \mathrm{ch}(Q(\overline{k})) \to \mathbb N$ to $h : Q^{JB}(\overline{k}) \to \mathbb N$ by defining it as
\[
h(x) = \sum_{c \in B(\overline{k})}h(f_c).
\]
Again this is well-defined because for the trivial chain $f = \widehat{1}$, we have $h(f) = 0$.
We further extend $h$ to $h : \mathrm{Hilb}(B)^Q(\overline{k}) \to \mathbb N$ by composing with $\mathrm{sat} : \mathrm{Hilb}(B)^Q(\overline{k}) \to Q^{JB}(\overline{k})$.

Finally we extend $h$ to pairs by $h(w \leq x) = h(w) - h(x)$ or $h(x) - h(w)$ depending on the context.

\end{defi}

We construct several functions on $Q^{JB}(\overline{k})$ using this construction following \cite[Section 3.5]{DT24}.

\begin{exam}[{\cite[Example 3.14]{DT24}}]
We fix $q_0 \in \widetilde{Q}(\overline{k})$, and we denote the indicator function of $\{q_0\}$ by $m_{q_0}$.
We call the extended $m_{q_0} : Q^{JB}(\overline{k}) \to \mathbb N$ as the number of occurrences of $q_0$.
This notion coincides with the notion used in \cref{exam:saturatedexample}.

\end{exam}

\begin{exam}[{\cite[Example 3.15]{DT24}}]
\label{exam:rankfunction}
We assume that $Q(k) = Q(\overline{k})$.
Assume that $Q(\overline{k})$ is a graded poset, i.e., every maximal chain in $Q(\overline{k})$ has the same length.
We first define the function $\mathrm{rank} : Q(\overline{k}) \to \mathbb N$ which assigns to $q \in Q(\overline{k})$ the length of the maximal chain starting at $\widehat{1} \in \mathrm{ch}(Q(\overline{k}))$ and ending at $q \in \mathrm{ch}(Q(\overline{k}))$. 
Then we extend this to the function $\mathrm{rank} : Q^{JB}(\overline{k}) \to \mathbb N$. The value $\mathrm{rank}(f)$ for $f \in \mathrm{ch}(Q(\overline{k}))$ is the length of the maximal chain starting at $\widehat{1} \in \mathrm{ch}(Q(\overline{k}))$ and ending at $f \in \mathrm{ch}(Q(\overline{k}))$.  For any $x \in Q^{JB}(\overline{k})$, the value $\mathrm{rank}(x)$ is equal to the cohomological dimension of the simplicial scheme $N([\widehat{1}, x]\cap Q^{JB}(\overline{k}))$.
\end{exam}

\begin{exam}[{\cite[Example 3.16]{DT24}}]
\label{exam:supportfunction}
For any $x \in Q^{JB}(\overline{k})$, we define $|\mathrm{Supp}(x)|$ as the number of $c \in B(\overline{k})$ such that $f_c$ is non-trivial.
\end{exam}

\begin{exam}
\label{exam:codimensionfunction}
Let $Q$ be the poscheme defined in \cref{sect:keyexample} and considered in \cref{exam:saturatedexample,exam:essentialsaturatedexample}.
We have $\mathrm{rank}(\ell_{i, 1} \oplus \ell_{i, 2}) = 1$, $\mathrm{rank}(V_j) =  1$, $\mathrm{rank}(\ell_{i, j}) = 2$, and $\mathrm{rank}(0) = 3$. This poset is graded.
Next we define $\gamma : Q(\overline{k}) \to \mathbb N$ by
\[
\gamma(\widehat{1}) = 0, \, \gamma(\ell_{i, 1} \oplus \ell_{i, 2}) = 2, \, \gamma(V_j) = 2, \gamma(\ell_{i, j}) = 3, \, \gamma(0) = 4.
\]
This extends to a combinatorial function $\gamma : Q^{JB}(\overline{k}) \to \mathbb N$. In \cref{subsec:blowups} we will interpret $\gamma(x)$ as the expected codimension of the incidence condition imposed by $x \in Q^{JB}(\overline{k})$.
\end{exam}

\subsection{Combinatorial stratifications of $Q^{JB}$}
Let $Q$ be a finite \'etale poscheme with $\widehat{1}$ and meets and let $B$ be a smooth projective curve.
We describe combinatorial stratifications of $\mathrm{Hilb}(B)^Q$ and $Q^{JB}$ following the exposition of \cite[Section 4]{DT24}.

\subsubsection{Stratifications by combinatorial types}

We define combinatorial types of $x \in \mathrm{Hilb}(B)^Q(\overline{k})$ following \cite[Section 4.1]{DT24}:
\begin{defi}[{\cite[Section 4.1]{DT24}}] \label{defn:combinatorial-type}
Let $x \in \mathrm{Hilb}(B)^Q(\overline{k})$ which corresponds to a finitely supported function
\[
B(\overline{k}) \to \mathrm{Hom}(\widetilde{Q}(\overline{k}), \mathbb N), \, c \mapsto \widetilde{g}_c.
\]
Then the combinatorial type of $x$ is the multiset
\[
\mathrm{type}(x) = \{\widetilde{g}_{c_1}, \dotsc, \widetilde{g}_{c_\ell}\}
\]
where $\{c_i \in B(\overline{k} )\, | \, i = 1, \dotsc, \ell\}$ is the support of $x$.
The set of combinatorial types is the set of multisets of homomorphisms $\widetilde{g} : \widetilde{Q} \to \mathbb N$.

Note that $x \in \mathrm{Hilb}(B)^Q(\overline{k})$ is in $Q^{JB}(\overline{k})$ if and only if every $g : Q (\overline{k})\to \mathbb N\cup \{\infty\} \in \mathrm{type}(x)$ is meet preserving. In this situation, we say that the multiset $\mathrm{type}(x)$ is saturated. 
\end{defi}

Note that there is a canonical bijection between saturated types and multisubsets of $\mathrm{ch}(Q(\overline{k}))$. Also the combinatorial type of the saturation $\mathrm{sat}(x)$ of $x\in \mathrm{Hilb}(B)^Q(\overline{k})$ only depends on the combinatorial type of $x$ showing that we have the saturation function from the set of combinatorial types to the set of saturated types.

Using combinatorial types, we obtain stratifications of $\mathrm{Hilb}(B)^Q$:

\begin{defi}[{\cite[Definition 4.1]{DT24}}]
Let $T$ be a combinatorial type. Then we define a Zariski locally closed subset $\mathcal N_{T, \overline{k}} \subset \mathrm{Hilb}(B)^Q_{\overline{k}}$ such that $\mathcal N_{T, \overline{k}}(\overline{k})$ is equal to the subset of $x \in \mathrm{Hilb}(B)^Q(\overline{k})$ such that $\mathrm{type}(x) = T$. Since $\mathcal N_{T, \overline{k}}$ is Galois invariant, it descends to a Zariski locally closed subset $\mathcal N_T \subset \mathrm{Hilb}(B)^Q$ defined over $k$.
\end{defi}

Note that the locally closed subset $\mathcal N_T$ is isomorphic to a form 
of the configuration space $\mathrm{Conf}(B;T)$ of $B$ labeled by $T$. Also we have
\[
\mathcal N_T\cap Q^{JB}\neq \emptyset \iff \mathcal N_T\subset Q^{JB}\iff \text{ $T$ is saturated.}
\]

\begin{defi}[{\cite[Definition 4.2]{DT24}}]
Let $\{\widetilde{g}_1, \dotsc, \widetilde{g}_\ell\}, \{ \widetilde{h}_1, \dotsc, \widetilde{h}_s\}$ be two combinatorial types. We say $\{\widetilde{g}_1, \dotsc, \widetilde{g}_\ell\} \geq \{ \widetilde{h}_1, \dotsc, \widetilde{h}_s\}$ if there is a injective map $\sigma : \{1, \dotsc, s\} \to \{1, \dotsc, \ell\}$ such that for any $j \in \{1, \dotsc, s\}$, we have
\[
\widetilde{g}_{\sigma(j)} \geq \widetilde{h}_j.
\]

We say $\{\widetilde{g}_1, \dotsc, \widetilde{g}_\ell\} \geq_+ \{ \widetilde{h}_1, \dotsc, \widetilde{h}_s\}$ if there is a surjective map $\sigma : \{1, \dotsc, \ell\} \to \{1, \dotsc, s\}$ such that for any $j \in \{1, \dotsc, s\}$, we have
\[
\sum_{i \in \sigma^{-1}(j)} \widetilde{g}_i \geq \widetilde{h}_j.
\] 
Also when $s = 0$, i.e., $\{ \widetilde{h}_1, \dotsc, \widetilde{h}_s\} = \emptyset$, we formally insist that we have $\{\widetilde{g}_1, \dotsc, \widetilde{g}_\ell\} \geq_+ \emptyset$.

\end{defi}

Note that with this definition, the set
\[
\bigcup_{T' \leq_+ T} \mathcal N_{T'},
\]
is closed in $\mathrm{Hilb}(B)^Q$ and is downward-closed with respect to the poscheme structure of $\mathrm{Hilb}(B)^Q$.

\begin{defi}[{\cite[Definition 4.5]{DT24}}]
Let $T$ be a saturated combinatorial type. We define
\[
\mathcal S_{T} = \bigcup_{T', \mathrm{sat}(T') = T} \mathcal N_{T'}.
\]
\end{defi}

Note that $\mathcal S_{T}$ agrees with the disjoint union $\sqcup_{T', \mathrm{sat}(T') = T} \mathcal N_{T'}$. This provides a coarse stratification of $\mathrm{Hilb}(B)^Q$.
We define the relation $\leq_{+, \mathrm{sat}}$ as the transitive closure of the following relation on the set of saturated combinatorial types: 
\[
T_1 \preceq_{+, \mathrm{sat}} T_2 \text{ if there exist types $T_1', T_2'$ with $T_1' \leq_+ T_2'$ and $\mathrm{sat}(T_i') = T_i$ for $i = 1, 2$.}
\]
This defines a poset structure on the set of saturated types, and the saturation function from the set of combinatorial types to the set of saturated combinatorial types is a homomorphism with respect to $\leq_+$ and $\leq_{+, \mathrm{sat}}$ (\cite[Proposition 4.6]{DT24}).

For any saturated type $T$, we define 
\[
\mathcal Z_T = \bigcup_{T' \text{ saturated}, T' \leq_{+, \mathrm{sat}} T} \mathcal S_{T'} = \bigcup_{\widetilde{T}, \, \mathrm{sat}(\widetilde{T})\leq_{+, \mathrm{sat}}T} \mathcal N_{\widetilde{T}}.
\]
It follows that $\mathcal Z_T$ is closed in $\mathrm{Hilb}(B)^Q$ and is downward-closed with respect to the poset structure of $\mathrm{Hilb}(B)^Q$. 
Also we have 
\[
\mathcal Z_{T_1} \subset \mathcal Z_{T_2} \iff T_1 \leq_{+, \mathrm{sat}} T_2,
\]
and 
\[
\mathcal Z_T \setminus \bigcup_{T' <_{+, \mathrm{sat}}T} \mathcal Z_{T'} = \mathcal S_T.
\]
See \cite[After Proposition 4.6]{DT24} for more details.

\subsection{Variants}
As usual $Q$ is a finite \'etale poscheme with $\widehat{1}$ and meets and $B$ is a smooth projective curve.

\noindent
{\bf Relative variant:}

In \cite[Section 4.3]{DT24}, the first and the fourth authors extended the notion of combinatorial types for the poscheme
\[
(\mathrm{Hilb}(B)^Q \leq \mathrm{Hilb}(B)^Q),
\]
consisting of pairs $(w \leq x)$ where $w, x \in \mathrm{Hilb}(B)^Q$.
Indeed, one may interpret this poscheme as $\mathrm{Hilb}(B)^{Q \times I}$ where $I$ is the poscheme consisting of two $k$-rational points $\widehat{0} \leq \widehat{1}$, and we give $Q \times I$ the product poscheme structure.
Any combinatorial type of this poscheme is a multiset of the form
\[
\{\widetilde{g}_{c_1}^w\leq \widetilde{g}_{c_1}^x,  \dotsc,  \widetilde{g}_{c_\ell}^w\leq \widetilde{g}_{c_\ell}^x\}.
\]
For each such combinatorial type $T$, we have a Zariski locally closed subset $\mathcal N_T \subset (\mathrm{Hilb}(B)^Q \leq \mathrm{Hilb}(B)^Q)$ defined over $k$.
We also have the relations $\leq_+$ and $\leq_{+, \mathrm{sat}}$ on the set of combinatorial types and the set of saturated types respectively. 
Let $W\subset \mathrm{Hilb}(B)^Q$ be a locally closed subset which is a finite union of locally closed subsets $\mathcal N_{T''}$ where $T''$ is a saturated type of $\mathrm{Hilb}(B)^Q$.  
Let $T$ be a saturated combinatorial type of $(W \leq \mathrm{Hilb}(B)^Q)$.
Then we define
\[
\mathcal Z_T =  \left(\bigcup_{T', \, \mathrm{sat}(T')\leq_{+, \mathrm{sat}}T} \mathcal N_{T'} \right)\cap (W < \mathrm{Hilb}(B)^Q),
\]
which is closed and downward-closed in $(W \leq \mathrm{Hilb}(B)^Q)$.

\noindent
{\bf Pointed variant:}

Assume that $B$ is pointed, i.e., we fix a $k$-rational point $* \in B$.
Since the Cohen-Jones-Segal conjecture addresses the space of pointed maps, we need to extend our definitions to pointed maps as explained in \cite[Section 4.4]{DT24}. A combinatorial type of $\mathrm{Hilb}(B, *)^Q$ is a multiset
\[
\{\widetilde{g}_{c_0}, \widetilde{g}_{c_1}, \dotsc, \widetilde{g}_{c_\ell}\},
\]
where $c_0 = *$. Note that $\widetilde{g}_{c_0}$ could be the trivial function.
We have
\[
\{\widetilde{g}_0, \widetilde{g}_1, \dotsc, \widetilde{g}_\ell\} \geq_+\{\widetilde{h}_0, \widetilde{h}_1, \dotsc, \widetilde{h}_s\}, 
\]
if there is a pointed surjective map $\sigma : \{0\} \cup \{1, \dotsc, \ell\} \to \{0\} \cup \{1, \dotsc, s\}$ such that 
for any $j \in \{0\}\cup \{1, \dotsc, s\}$ we have $h_j \leq \sum_{i \in \sigma^{-1}(j)} g_i$.
This induces $\leq_{+, \mathrm{sat}}$ and the stratification $\mathcal Z_T$ and $\mathcal S_T$.

\section{Our common set up: blow ups of a smooth quadric surface}
\label{subsec:blowups}

We next describe the common setup for the Cohen--Jones--Segal conjecture in \cref{sec:CJS} and for Manin's conjecture in \cref{sec:Manin}.  Recall from the introduction that we will think of a split degree $4$ del Pezzo surface $S$ as a blow-up of a quadric surface $\mathbb{P}^{1} \times \mathbb{P}^{1}$ at four points.  This section explains how to relate curves on $S$ and $\mathbb{P}^{1} \times \mathbb{P}^{1}$.

Throughout this section we assume that our ground field $k$ is a perfect field.

\subsection{Blow-ups of quadric surfaces}

We consider a smooth quadric surface $\mathbb{P}^{1} \times \mathbb{P}^{1}$
with the projection maps $\pi_j:\mathbb{P}^{1} \times \mathbb{P}^{1} \to \mathbb P^1$ for $j = 1, 2$.
Let $r \geq 0$ be a positive integer and fix $r$ distinct $k$-rational points $(p_1, p_1'), \dotsc, (p_r, p_r')$ on $\mathbb{P}^{1} \times \mathbb{P}^{1}$.  Let $\rho_r : S_r \to \mathbb{P}^{1} \times \mathbb{P}^{1}$ be the blow up of $\mathbb{P}^{1} \times \mathbb{P}^{1}$ at $(p_1, p_1'), \dotsc, (p_r, p_r')$.  We will always make the following:

\begin{assu}
The points $\{ p_i \}_{i=1}^{r}$ are distinct and the points $\{ p'_i \}_{i=1}^{r}$ are distinct. 
\end{assu}

Note that this assumption is satisfied for any birational morphism $S \to \mathbb{P}^{1} \times \mathbb{P}^{1}$ from a split degree $4$ del Pezzo surface $S$.

For a fixed smooth projective curve $B$, our goal is to understand the cohomology of 
\[
\mathrm{Mor}(B, S_r).
\]
Let $E_1, \dotsc, E_r$ be the exceptional divisors above $(p_1, p_1') \dotsc, (p_r, p_r')$ and for $j=1,2$ let $F_j$ be a general fiber of $\pi_j \circ \rho_r : S_r \to \mathbb P^1$. By abuse of notation, we also denote a general fiber of $\pi_j$ by $F_j$ as well.

Given a nef curve class $\alpha$ on $S_r$, we define
\[
-K_{S_r}.\alpha = h(\alpha), \qquad a(\alpha) = F_1. \alpha, \qquad  a'(\alpha) = F_2.\alpha, \qquad E_i.\alpha = k_i(\alpha).
\]
Then we have
\[
h(\alpha) = 2a(\alpha) + 2a'(\alpha) -\sum_{i = 1}^r k_i(\alpha).
\]
Note that once we fix $a(\alpha), a'(\alpha), k_1(\alpha), \dotsc, k_r(\alpha)$, then the numerical class $\alpha$ is uniquely determined.  Furthermore this notation is consistent with the special case of degree $4$ del Pezzo surfaces established in \cref{nota:bignota}.

\subsection{Parametrizing curves}

Fix a smooth projective curve $B$ defined over $k$.  A morphism $f : B \to \mathbb{P}^{1} \times \mathbb{P}^{1}$ is exactly the same as a pair of morphisms from $B$ to the two $\mathbb{P}^{1}$ factors.  Thus there is a bijective correspondence between morphisms $f: B \to \mathbb{P}^{1} \times \mathbb{P}^{1}$ and equivalence classes of pairs of line bundles $L_{1},L_{2}$ on $B$ equipped with non-vanishing sections
\[
s \in H^0(B, L_1^{\oplus 2}), \, t \in H^0(B, L_2^{\oplus 2}),
\]
where the equivalence is defined by the rescaling action of $s$ and $t$ separately.  

Let us write $\mathbb{P}^{1} \times \mathbb{P}^{1} = \mathbb P(V_1) \times \mathbb P(V_2)$ where each $V_i$ is a $2$-dimensional vector space.  We fix two non-negative integers $a, a'$ and we set $\mathcal B = \Pic^a(B) \times \Pic^{a'}(B)$.
Let $\mathcal L_a, \mathcal L_{a'}$ be the line bundles on $B \times \mathcal B$ which are pullbacks of the universal line bundles on $B \times \Pic^a(B)$ and $B\times \Pic^{a'}(B)$ respectively. Then we set
\[
\mathcal V_1 = \mathcal L_a \otimes V_1, \, \mathcal V_2 = \mathcal L_{a'} \otimes V_2, \, \mathcal V = \mathcal V_1 \oplus \mathcal V_2,
\]
which we consider both as locally free sheaves on $B \times \mathcal B$ and as the corresponding vector bundles.  We define the scheme of sections $\Gamma(B, \mathcal V) \to \mathcal B$ which represents the functor
\[
(\phi : T \to \mathcal B) \mapsto \Gamma(B \times T, (\mathrm{id}, \phi)^*\mathcal V)
\]
There is a Zariski open subscheme of $\Gamma(B,\mathcal{V})$ parametrizing pairs of sections $(s,t) \in H^{0}(B,\mathcal{V}_{1} \oplus \mathcal{V}_{2})$ such that neither $s$ nor $t$ vanishes on $B$.  This open subscheme defines a $\mathbb{G}_{m}^{2}$-torsor over the moduli space of morphisms $f: B \to \mathbb{P}^{1} \times \mathbb{P}^{1}$ of bidegree $(a,a')$.

We next turn to morphisms $f': B \to S_{r}$.  For a nef curve class $\alpha$ on $S_{r}$, we let $M_{\alpha}$ denote the moduli space of morphisms $f': B \to S_{r}$ whose numerical class lies in $\alpha$.  Note that choosing a morphism $f' : B \to S_r$ such that $s'_*[B] = \alpha$ is equivalent to choosing a morphism $f : B \to \mathbb{P}^{1} \times \mathbb{P}^{1}$ such that
\begin{itemize}
\item $\deg \, f^*F_1 = a(\alpha)$;
\item $\deg \, f^*F_2 = a'(\alpha)$; and
\item for any $i = 1, \dotsc, r$, the local multiplicity of $f(B)$ at $(p_i, p_i')$ is given by $k_i(\alpha)$.
\end{itemize}
For each index $i=1,\dotsc,r$ there are $1$-dimensional subspaces $\ell_{i, 1} \subset V_1$ and $\ell_{i, 2} \subset V_2$ corresponding to $p_i$ and $p_i'$ respectively.  Then the local multiplicity $k_i(\alpha)$ of $f$ at $(p_i, p_i')$ is given by the length of the subscheme $(s, t)^{-1}(\ell_{i, 1}\oplus \ell_{i, 2})$ on $B$.

Let $\widetilde{M}_\alpha$ denote the locally closed subscheme of $\Gamma(B,\mathcal{V})$ such that the fiber over $(L_{1},L_{2}) \in \mathcal{B}$ parametrizes all pairs of nonvanishing sections $(s,t) \in H^{0}(B,\mathcal{V}_{1} \oplus \mathcal{V}_{2})$ such that the length of the subscheme $(s, t)^{-1}(\ell_{i, 1}\oplus \ell_{i, 2})$ is $k_i(\alpha)$ for every $i = 1, \dotsc, r$. 
 Then we have a natural morphism
\[
\widetilde{M}_\alpha \to M_\alpha
\]
which is a $\mathbb G_m^2$-torsor over $M_\alpha$.

\subsection{Stratifying the space of sections}
\label{subsubsec:stratifyingthespace of sectons}

Finally, we define a stratification of $\Gamma(B,\mathcal{V})$.  Let $Q$ be the poscheme defined by the following set of subspaces of $V$ as in \cref{sect:keyexample}: 
\begin{equation*}
\{ V, V_1 \oplus \{0\}, \{0\} \oplus V_2,  0 \} \cup \{ \ell_{i, 1} \oplus \ell_{i, 2}, \ell_{i, 1} \oplus \{0\}, \{0\} \oplus \ell_{i, 2}  \}_{i=1,\dotsc,r}.
\end{equation*}
For any subspace $K_q \in Q$, we define the corresponding subbundle $\mathcal K_q \subset \mathcal V$ over $B \times \mathcal{B}$.
Since the poscheme $Q$ is closed under intersection, we have
\[
\mathcal K_{q_1} \cap \mathcal K_{q_2} = \mathcal K_{q_1 \wedge q_2},
\]
for any $q_1, q_2\in Q$.

For any $q \in Q$ we also have the universal divisor 
\[
\mathcal D_q' \subset B \times \mathrm{Hilb}(B)^Q
\]
such that the restriction of $\mathcal D_q'$ over $D \in \mathrm{Hilb}(B)^Q(T)$ is the divisor $D_q$ on $B \times T$.
We denote the pullback of $\mathcal D_q'$ via $B \times \mathcal B \times \mathrm{Hilb}(B)^Q \to B\times \mathrm{Hilb}(B)^Q$ by $\mathcal D_q$.
We also denote the projection $B \times \mathcal B \times\mathrm{Hilb}(B)^Q\to B \times \mathcal B$ by $p$.
Then we define $\Gamma_Q(B, \mathcal V)$ as the fiber product:  
\[
\xymatrix{
  \Gamma_Q(B, \mathcal V) \ar[r] \ar[d] & \prod_{q \in \widetilde{Q}}\Gamma(\mathcal D_q, p^*\mathcal K_q) \ar[d] \\
  \Gamma(B, \mathcal V) \times \mathrm{Hilb}(B)^Q \ar[r] & \prod_{q \in \widetilde{Q}}\Gamma(\mathcal D_q, p^*\mathcal V).
}
\]
Thus $\Gamma_Q(B, \mathcal V)$ is the closed subscheme of $\Gamma(B, \mathcal V) \times \mathrm{Hilb}(B)^Q$ such that the fiber over the point $\mathfrak{b} \in \mathcal{B}$ parametrizes pairs of a section $s \in \Gamma(B, \mathcal V_\mathfrak{b})$ and an element $(x_{q}) \in \mathrm{Hilb}(B)^Q$ satisfying $s(x_q) \subset (\mathcal K_q)_{\mathfrak{b}}$ for every $q \in Q$.
Note that the scheme $ \Gamma_Q(B, \mathcal V)$ is a $(\mathcal B \times \mathrm{Hilb}(B)^Q)$-scheme.
If we have $(x_{q}) \leq (x'_{q})$ in $ \mathrm{Hilb}(B)^Q$ so that we have containments $D_q \subset D'_q$ for every $q$, then for every $\mathfrak{b} \in \mathcal{B}$ we have 
\[
  \Gamma_Q(B, \mathcal V)_{\mathfrak b, x'} \subset  \Gamma_Q(B, \mathcal V)_{\mathfrak b, x}.
\]
In this way, $\Gamma(B, \mathcal V)$ admits a stratification parametrized by the poscheme $\mathrm{Hilb}(B)^Q$.  This stratification is insensitive to taking saturations, in the sense that for any $x \in \mathrm{Hilb}(B)^{Q}(k)$ and $\mathfrak{b} \in \mathcal{B}$ we have $\Gamma(B,\mathcal{V})_{\mathfrak{b},x} = \Gamma(B,\mathcal{V})_{\mathfrak{b},\mathrm{sat}(x)}$.

This stratification is compatible with the construction of moduli spaces of curves described above.  Suppose we fix a nef numerical class $\alpha$ on $S_{r}$.  Consider the sublocus $U_{\mathbf{k}(\alpha)} = \cup_{T} \mathcal{N}_{T}$ as $T$ varies over saturated combinatorial types of the form
\begin{equation*}
T = \left\{ m_{ij}[\ell_{i,1} \oplus \ell_{i,2}] \left| \sum_{j=1}^{t_{i}} m_{ij} = k_{i}(\alpha) \right. \right\}_{\substack{i=1,\dotsc,r \\ j=1,\dotsc,t_{i}}}.
\end{equation*}
Let $E \to \mathcal B \times U_{\mathbf{k}(\alpha)}$ be the restriction of $\Gamma_Q(B, \mathcal V) \to \mathcal B \times \mathrm{Hilb}(B)^Q$ to $ \mathcal B \times U_{\mathbf{k}(\alpha)}$.
Then the moduli space $\widetilde{M}_{\alpha}$ can be described as 
\begin{equation*}
\widetilde{M}_{\alpha} = \bigcup_{w \in U_{\mathbf{k}(\alpha)}}(E_w \backslash \left( \cup_{(w < x) \in (W_\alpha < \mathrm{Hilb}(B)^Q)} )\Gamma(B,\mathcal{V})_{x} \right),
\end{equation*}
where $E_w = \Gamma_Q(B, \mathcal V)_w$.

\section{Bar complexes over perfect fields}
\label{sec:barconstruction}

In this section we prove the main technical result that compares the cohomology of a space over $U_{\mathbf{k}}$ to the cohomology of an associated bar complex.  In particular we extend \cite[Theorem 5.4]{DT24} from $\mathbb{C}$ to arbitrary perfect fields.  While the treatment in \cite{DT24} relies on the geometric realizations of simplicial spaces, our discussion is based on simplicial schemes and their sheaf theoretic constructions.

\subsection{Stratifying the bar complex}
Let $B$ be a smooth projective curve with possibly a fixed base point.
We consider the non-pointed case and pointed case simultaneously.

Let $W\subset \mathrm{Hilb}(B)^Q$ be a locally closed set which is a finite union of locally closed subsets $\mathcal N_T$ where $T$ is a saturated type. We equip $W$ with the reduced scheme structure so that $W$ is a separated reduced scheme of finite type over $k$. Let $P \subset (W < \mathrm{Hilb}(B)^Q)$ be a closed and downward closed union of locally closed sets $\mathcal S_T$ for finitely many saturated types $T$ of $(W < \mathrm{Hilb}(B)^Q)$ equipped with the reduced scheme structure.
We assume that the projection $P \to W$ is proper.
Let $\mathcal B$ be a separated scheme of finite type over $k$ and let $U \subset \mathcal B \times W$ be an open subset. We write $(U \leq \mathrm{Hilb}(B)^Q)$ and $P_U \to U$ for the pullbacks of $(W \leq \mathrm{Hilb}(B)^Q) \to W$ and $P \to W$ by $U \subset \mathcal B \times W \to W$.

Suppose that we have a flat separated morphism $E \to U$ of finite type with irreducible geometric fibers which admits a stratification by the poset $(U \leq \mathrm{Hilb}(B)^Q)$ corresponding to the closed subscheme $Z \subset E\times_U (U \leq \mathrm{Hilb}(B)^Q)$.
We assume that for any $(\mathfrak b, w < x) \in (U < \mathrm{Hilb}(B)^Q)(\overline{k})$, the stratification $Z_{\mathfrak b, w<x} \subset E_{\mathfrak b, w}$ is an irreducible closed subscheme 
and that the inclusion
\[
Z_{\mathfrak b, w < x} \supset Z_{\mathfrak b, w < \mathrm{sat}(x)}
\]
is the equality. We also suppose that when we have two saturated elements $w \leq x_1$ and $w \leq x_2$, then we have
\[
Z_{\mathfrak b, w\leq x_1}\cap Z_{\mathfrak b, w \leq x_2} = Z_{\mathfrak b, w \leq x_1\vee x_2}.
\]

We are interested in the bar complex $B(P_U, Z_{P_U})$ with simplices
\[
\{ (\mathfrak b, w) \in U, w < x_0 \leq \dotsb \leq x_m \in \mathrm{Hilb}(B)^Q, z \in Z_{\mathfrak b, w < x_m}\, | \, (w < x_i) \in P \},
\]
with the augmentation morphism $\epsilon : B(P_U, Z_{P_U}) \to E_U$ given by
\[
((\mathfrak b, w) \in U, w < x_0 \leq \dotsb \leq x_m, z)\mapsto (\mathfrak b, w, z).
\]

\subsection{Spectral sequence}
Since $B(P_{U},Z_{P_U})$ is a bar complex, we can analyze its cohomology using the spectral sequence in \cref{subsec:spectralsequenceforstratifiedspaces}.  Our first task is to reinterpret the terms in this spectral sequence as the cohomology of certain complexes of sheaves on $Z|_{U \times_{W} \mathcal{N}_{T}}$.  The advantage of this perspective is that only the essential saturated types $T$ make a non-trivial contribution.

Let $\mathfrak T$ be the set of saturated combinatorial types in $P$ with ordering $\leq_+$.
We fix an injective homomorphism $h : \mathfrak T \to \mathbb N$.
Let $\mathcal Z_i = \cup_{h(T)\leq i} \mathcal Z_T$ and consider
\[
B_i = B(U \times_W \mathcal Z_i, Z_{P_U}) \subset B(P_U, Z_{P_U}).
\]  
The closed embeddings $\iota'_{i}: B_{i} \to B(P_U, Z_{P_{U}})$ define a stratification into closed subschemes. Retaining the notation of \cref{subsec:spectralsequenceforstratifiedspaces} to define a grading of $\underline{\mathbb{Z}}_{\ell}$ associated to this stratification, we obtain the following spectral sequence:
\[
\bigoplus_{i + j = n} E_1^{i, j} = \bigoplus_{i} H^{n}_{\text{\'et}, c}(B(P_U, Z_{P_U})_{\overline{k}}, \mathrm{Gr}_F^i \underline{\mathbb{Z}}_\ell) \implies H^{n}_{\text{\'et}, c}(B(P_U, Z_{P_U})_{\overline{k}}, \underline{\mathbb{Z}}_\ell).
\]
We consider $Z|_{U \times_W \mathcal N_{h^{-1}(i)}}$ and the associated nerve 
\[
\epsilon_i: (N(-\infty, Z|_{U\times_W \mathcal S_{h^{-1}(i)}} )< Z|_{U \times_W \mathcal N_{h^{-1}(i)}}) \to Z|_{U\times_W \mathcal N_{h^{-1}(i)}}
\]
which is defined to be the simplicial scheme over $Z|_{U \times_W \mathcal N_{h^{-1}(i)}}$ with $m$-th simplex given by 
\[
\{\mathfrak b, w < x_0\leq \dotsb \leq x_m < y, z \, | \, (\mathfrak b, w) \in U, (w < y)\in \mathcal N_{h^{-1}(i)}, (w < x_m) \in \mathcal Z_{h^{-1}(i)}  \setminus \mathcal S_{h^{-1}(i)},  z \in Z_{\mathfrak b, w < y} \},
\]
and the augmentation morphism $\epsilon_i : (N(-\infty, Z|_{U\times_W \mathcal S_{h^{-1}(i)}} )< Z|_{U \times_W \mathcal N_{h^{-1}(i)}})\to Z|_{U \times_W \mathcal N_{h^{-1}(i)}}$ is
\[
\epsilon_i : (\mathfrak b, w < x_0\leq \dotsb \leq x_m < y, z) \mapsto (\mathfrak b, w < y, z).
\]
We define the complex $\mu'(h^{-1}(i))[2] \in D^b(Z|_{U \times_W \mathcal N_{h^{-1}(i)}})$ to be the cone yielding the distinguished triangle
\[
 \underline{\mathbb Z}_\ell \to R\epsilon_{i*}\epsilon_i^*\underline{\mathbb Z}_\ell \to \mu'(h^{-1}(i))[2] \to \underline{\mathbb Z}_\ell[1] \to.
\]
Then we have the following proposition:
\begin{prop} \label{prop:ssforbarconstruction}
\[
H^{i + j}_{\text{\'et}, c}(B(P_U, Z_{P_U})_{\overline{k}}, \mathrm{Gr}_F^i\underline{\mathbb Z}_{\ell} ) \cong H^{i + j}_{\text{\'et}, c}((Z|_{U \times_W \mathcal N_{h^{-1}(i)}})_{\overline{k}}, \mu'(h^{-1}(i))[1]).
\]
\end{prop}

\begin{proof}
It suffices to show that we have
\[
R\epsilon_*\mathrm{Gr}_F^i\underline{\mathbb Z}_{\ell} \cong R\pi_{i!}\mu'(h^{-1}(i))[1],
\]
where $\epsilon$ is the augmentation morphism for $B(P_{U},Z_{P_U})$ and   $\pi_i : Z|_{U \times_W \mathcal N_{h^{-1}(i)}} \to E|_U$ is the projection.
Consider the relative poscheme $(-\infty,Z|_{U \times_W \mathcal N_{h^{-1}(i)}}]$ over $Z|_{U \times_W \mathcal N_{h^{-1}(i)}}$ as defined in \cite[Definition 3.1.3]{DH24}.
Let $\xi_i : N(-\infty, Z|_{U \times_W \mathcal N_{h^{-1}(i)}}] \to B_i$ 
be the morphism of simplicial schemes defined by 
\[
(w < x_0\leq \dotsb \leq x_m \leq y, z) \mapsto (w < x_0\leq \dotsb \leq x_m, z).
\]
We have an exact sequence on $B_{i}$:
\[
0 \to \mathrm{Gr}_F^i\underline{\mathbb Z}_{\ell}  \to \underline{\mathbb Z}_{\ell}|_{(B_{i})_{\overline{k}}} \to R\iota'_{i-1!}\underline{\mathbb Z}_{\ell}|_{(B_{i-1})_{\overline{k}}}  \to 0.
\]
Let $U_m = B_i([m]) \setminus B_{i-1}([m])$. We denote $U_m \hookrightarrow B_i([m])$ by $j_m$.  
The exact sequence above allows us to identify $(\mathrm{Gr}_F^i\underline{\mathbb Z}_{\ell})_m \cong Rj_{m!}j_m^*(\mathrm{Gr}_F^i\underline{\mathbb Z}_{\ell})_m$.  This implies that the natural morphism $\mathrm{Gr}_F^i\underline{\mathbb Z}_{\ell} \to R\xi_{i*}\xi_i^*\mathrm{Gr}_F^i\underline{\mathbb Z}_{\ell}$ is an isomorphism. Indeed this follows from the fact that the map
\[
\xi_i([m])^{-1}(U_m) \to U_m
\]
is an isomorphism because the inverse morphism is given by
\[
(w < x_0\leq \dotsb \leq x_m, z) \mapsto (w < x_0\leq \dotsb \leq x_m \leq \mathrm{sat}(x_m), z). 
\]
Then our assertion follows from the computation of \cite[Lemma 4.3.3]{DH24}
with $X = Z|_{U \times_W \mathcal N_{h^{-1}(i)}}$.
\end{proof}

Thus we have the following spectral sequence:
\begin{equation*}
\bigoplus_{i + j = n}E_1^{i, j} = \bigoplus_{T: \text{ saturated types of $P$}} H^{n}_{\text{\'et}, c}((Z|_{U \times_W\mathcal N_{T}})_{\overline{k}}, \mu'(T)[1]) \implies H^{n}_{\text{\'et}, c}(B(P_U, Z_{P_U})_{\overline{k}}, \mathbb Z_\ell).
\end{equation*}
Next we record a lemma regarding the complex $\mu'(\bullet)$. 
\begin{lemm}
When $T$ is not essential, we have $\mu'(T) \cong 0$.
\end{lemm}
\begin{proof}
Let $(\mathfrak b, w < y, z) : \Spec \, \kappa \to Z|_{U \times_W \mathcal N_T}$ be a geometric point.
Since $T$ is not essential, as explained in the proof of \cref{prop:homotopyequivalence} 
there is a maximum essential element $\mathrm{ess}(y)$ with $w < \mathrm{ess}(y) < y$ and $(w < \mathrm{ess}(y)) \not\in \mathcal S_T(\kappa)$. 
Indeed, $\mathrm{ess}(y)$ is defined as the join of $y'$ such that $w \prec y' \leq y$.
This implies that the nerve 
\[
(N(-\infty, Z|_{U \times_W \mathcal S_T}) < Z|_{U \times_W \mathcal N_T})_{(\mathfrak b, w < y, z)} \to \Spec \, \kappa
\]
admits a weak center in the sense of \cite[Definition 3.1.6]{DH24}. Thus it follows from \cite[Theorem 4.0.1]{DH24} that we have an isomorphism
\[
\underline{\mathbb Z}_\ell \cong R\epsilon_{T*}\epsilon_T^* \underline{\mathbb Z}_\ell,
\]
where $\epsilon_T : (N(-\infty, Z|_{U\times_W \mathcal S_{T}} )< Z|_{U \times_W \mathcal N_{T}}) \to Z|_{U\times_W \mathcal N_{T}}$ is the augmentation morphism.
Our assertion follows.
\end{proof}

Altogether we conclude the following theorem, extending an analogous statement over $\mathbb{C}$ proved by \cite[Theorem 5.3]{DT24}.

\begin{theo}
\label{theo:spectralsequenceforbar}
There is a spectral sequence:
\[
\bigoplus_{i + j = n} E_1^{i, j} = \bigoplus_{T: \text{ essential types of $P$}} H^{n}_{\text{\'et}, c}((Z|_{U \times_W \mathcal N_T})_{\overline{k}}, \mu'(T)[1]) \implies H^{n}_{\text{\'et}, c}(B(P_U, Z_{P_U})_{\overline{k}}, \mathbb Z_\ell).
\]
\end{theo}

Finally for a later application, we record the following lemma comparing saturated and non-saturated types.
\begin{lemm}
\label{lemm:deformationretract}
The following inclusion over $Z|_{U \times_W \mathcal N_{T}}$ 
\[
(N((-\infty, Z|_{U\times_W \mathcal S_{T}} ) \cap Q^{JB})< Z|_{U \times_W \mathcal N_{T}}) \hookrightarrow (N(-\infty, Z|_{U\times_W \mathcal S_{T}} )< Z|_{U \times_W \mathcal N_{T}})
\]
induces an isomorphism
\[
R\epsilon'_{T*}\epsilon_T'^*\underline{\mathbb Z}_\ell \cong R\epsilon_{T*}\epsilon_T^*\underline{\mathbb Z}_\ell ,
\]
where $\epsilon'_T : (N((-\infty, Z|_{U\times_W \mathcal S_{T}} )\cap Q^{JB})< Z|_{U \times_W \mathcal N_{T}}) \to Z|_{U\times_W \mathcal N_{T}}$ is the augmentation morphism.
\end{lemm}
\begin{proof}
A similar proof to \cref{prop:homotopyequivalence} will justify this. 
Indeed, for an algebraically closed field $\kappa$, the map
\begin{align*}
((-\infty, Z|_{U\times_W \mathcal S_{T}} )< Z|_{U \times_W \mathcal N_{T}})(\kappa) &\to (((-\infty, Z|_{U\times_W \mathcal S_{T}} ) \cap Q^{JB})< Z|_{U \times_W \mathcal N_{T}})(\kappa), \\
(\mathfrak b, w < x < y) &\mapsto (\mathfrak b, w < \mathrm{sat}(x)  < y),
\end{align*}
is a left adjoint functor to the inclusion functor 
\begin{equation*}
(((-\infty, Z|_{U\times_W \mathcal S_{T}} ) \cap Q^{JB})< Z|_{U \times_W \mathcal N_{T}})(\kappa)  \hookrightarrow ((-\infty, Z|_{U\times_W \mathcal S_{T}} )< Z|_{U \times_W \mathcal N_{T}})(\kappa).
\end{equation*}
Because nerves take adjoint functors to homotopy equivalences, it follows that for a geometric point $(\mathfrak b, w < y , z) \in Z|_{U \times_W \mathcal N_{T}}(\kappa)$, the inclusion functor between the nerves
\[
N(((w, (\mathfrak b, w) \times_W \mathcal S_T)\cap Q^{JB}) < y)(\kappa) \hookrightarrow N((w, (\mathfrak b, w) \times_W \mathcal S_T) < y)(\kappa)
\]
is a deformation retract (see \cite[Lemma 3.2.1]{DH24}). Thus \cite[Lemma 2.2.2]{DH24} shows that the inclusion induces an isomorphism on the level of stalks:
\[
(R\epsilon'_{T*}\epsilon_T'^*\underline{\mathbb Z}_\ell)_{(\mathfrak b, w < y , z)} \cong (R\epsilon_{T*}\epsilon_T^*\underline{\mathbb Z}_\ell)_{(\mathfrak b, w < y , z)}.
\]
Our assertion follows.
\end{proof}

\subsection{The main theorem on the bar complex}

Suppose that we have the codimension function for pairs
\[
\gamma: (W\leq \mathrm{Hilb}(B)^Q) \to \mathbb N,
\]
which is given by the formula $\gamma(w \leq x) = \gamma(x) - \gamma(w)$ where  $\gamma(x)$ only depends on the combinatorial type of $x$.
(In our applications, $Q$ is the poscheme defined in \cref{subsubsec:stratifyingthespace of sectons} and considered in \cref{exam:saturatedexample,exam:essentialsaturatedexample}.
Then the codimension function is given in \cref{exam:codimensionfunction}.)
Following \cite[Section 5.4]{DT24}, we say that $Z_{\mathfrak b, w \leq x}$ is unobstructed if its codimension in $E_{\mathfrak b, w}$ is given by $\gamma(w\leq x)$. Using this, we define a function $\kappa$ on the set of combinatorial types as
\begin{equation} \label{eq:kappadef}
\kappa(T) = 2\gamma (T) - \mathrm{rank}(T) - 2|\mathrm{Supp}(T)|.
\end{equation}
Here we use the definitions $\mathrm{rank}(w\leq x) = \mathrm{rank}(x) - \mathrm{rank}(w)$ and $|\mathrm{Supp}(w\leq x)| = |\mathrm{Supp}(x-w)|$ where $\mathrm{rank}$ and $|\mathrm{Supp}|$ are defined in \cref{exam:rankfunction,exam:supportfunction}, and $x-w$ is viewed as a $0$-cycle on $B$.
In the pointed case, our convention is that 
\[
|\mathrm{Supp}(g_0 \geq h_0, g_1 \geq h_1, \dotsc, g_\ell \geq h_\ell)| =  \#\{ i>0 | g_i > h_i\}.
\]
One may think of $\kappa(T)$ as the expected cohomological codimension of the contribution of $\mathcal S_T$. 

Let $R \subset(W\leq \mathrm{Hilb}(B)^Q)$ be a closed and downwards closed union of strata $\mathcal S_T$ for $T$ a saturated type.  
We assume that $R$ is initial in its closed upwards closure $\widetilde{P}$ in the following sense: there exists a morphism $\iota_R : \widetilde{P} \to \mathrm{Hilb}(B)^Q$ such that for any $(w<y) \in \widetilde{P}(\overline{k})$ 
the element $(w< \iota_R(w<y))$ is contained in $R(\overline{k})$ and is the maximal element of the set 
$\{ (w < y') \in R(\overline{k}) \, | \, y' \leq y \}$.
We denote the pullback of $Z$ via $P_U \to U$ by $Z_{P_U}$. We have the evaluation map $Z_{P_U} \to E|_U$ and we denote its image by $\mathrm{im}(Z_{P_U} \to E|_U)$. The following theorem was proved over $\mathbb C$ in \cite[Theorem 5.9]{DT24}. We extend this to arbitrary perfect fields using a similar approach.

\begin{theo}
\label{theo:5.4inDT}
Let $I \in \mathbb Q_{\geq 0 }$.  Assume that
\begin{enumerate}
\item $P \to W$ is proper and $P$ is downward closed;
\item for every pair $(\mathfrak b, w < x)\in P_U(\overline{k})$ with $x$ being saturated and every $y \in Q^{JB}(\overline{k})$ such that $x \prec y$, 
the fibers $Z_{b, w < x}$, $Z_{b, w < y}$ are unobstructed, and; 
\item we have $P \subset R$ and $P$ contains every $\mathcal S_T$ such that $\mathcal S_T \subset R$ and $\kappa(T) \leq I$.
\end{enumerate}
Then the map on compactly supported \'etale cohomology induced by 
$$
B(P_U, Z_{P_U}) \to \mathrm{im}(Z_{P_U} \to E|_U)
$$
is cohomology connected in codimension $I+2$, i.e., for all $i > 2\mathrm{dim}(E|_U) - \lfloor I \rfloor -2$, the map 
\[
H^i_{\text{\'et}, c}(\mathrm{im}((Z_{P_{U}})_{\overline{k}} \to (E|_{U})_{\overline{k}}), \mathbb Z_\ell) \to H^i_{\text{\'et}, c}(B(P_U, Z_{P_U})_{\overline{k}}, \mathbb Z_\ell),
\]
is an isomorphism, and for 
$i = 2\dim\, E|_U - \lfloor I \rfloor -2$, it is a surjection.
\end{theo}

\begin{proof}
We replace $I$ by $\lfloor I \rfloor$ so that $I$ is an integer.
We denote the augmentation morphism $B(P_U, Z_{P_U}) \to \mathrm{im}(Z_{P_U} \to E|_U)$ by $\epsilon$. Then since $P$ is proper over $W$, $\epsilon$ is also proper, showing that we have an induced morphism
\[
H^i_{\text{\'et}, c}(\mathrm{im}((Z_{P_{U}})_{\overline{k}} \to (E|_{U})_{\overline{k}}), \mathbb Z_\ell) \to H^i_{\text{\'et}, c}(B(P_U, Z_{P_U})_{\overline{k}}, \mathbb Z_\ell).
\]
We let $P'\supset P$ be the closure, downward closure, and the saturation closure of the set of $w < y$ such that there is a $(w < x) \in P$ such that $x\prec y$. (We take these closures over $\overline{k}$ and descend to $k$.)  Now $P'$ is closed, downward closed, and proper over $W$ and it is a union of finitely many $\mathcal S_T$ where $T$ is a saturated combinatorial type.
We will stratify $\mathrm{im}(Z_{P_{U}} \to E|_{U})$ using the poset of saturated combinatorial types of $P'$.

Consider the proper projection $\Psi: Z_{P'_U} \to E|_U$.
For every saturated combinatorial type $T$ of $P'$ we define $\overline{J}_T \subset E|_U$ as the closed subset $\Psi(\overline{Z|_{U \times_W \mathcal N_{T}}}) \subset E|_U$ where the closure is taking place in $Z_{P_U}$.  
Whenever we have $T \leq T'$ (note that this relation is neither $\leq_+$ nor $\leq_{+, \mathrm{sat}}$), we have $\overline{J}_T \supset \overline{J}_{T'}$. Let $J_T = \overline{J}_T \setminus \cup_{T < T'} \overline{J}_{T'}$. Note that this stratification may not satisfy the maximal element property stated in \cref{subsec:spectralsequenceforstratifiedspaces}.

Using the spectral sequence in \cref{eq:stratss}, it suffices to show that for any saturated combinatorial type $T$ of $P'$ and any open subset $J_T' \subset J_T$ the morphism $\epsilon : \epsilon^{-1}(J'_T) \to J'_T\cap \mathrm{im}(Z_{P_U} \to E|_U)$ is 
cohomology connected in codimension $I+2$.
We fix a geometric point $(\mathfrak b, w, z) : \Spec \, \kappa \to J'_T\cap \mathrm{im}(Z_{P_U} \to E|_U)$ such that $z \in J'_T \cap E_{\mathfrak b, w}$.
Then the fiber of $\epsilon$ at $(\mathfrak b, w, z)$ is the nerve 
\[
[n] \mapsto \{w < x_0 \leq \dotsb \leq x_n \, | \, x_i \in P, z \in Z_{\mathfrak b, w < x_n} \}.
\]

Suppose that $T$ is a saturated type of $P'$ such that $\iota_R(T)$ is a type of $P$. We show that under this assumption, the poscheme 
\[
\{ x \in P \, | \, z \in Z_{\mathfrak b, w < x} \} \to \Spec \, \kappa,
\] 
admits a terminal element. Once we have this, then $\epsilon : \epsilon^{-1}(J'_T) \to J'_T\cap \mathrm{im}(Z_{P_U} \to E|_U)$ is sufficiently cohomology connected, by proper base change and the contractibility of a nerve of a poset with a terminal element, or by \cite[Theorem 4.0.1]{DH24}. 
Let $y \in (Z|_{U \times_W \mathcal N_{T}})_{w, \kappa}$ be a generic point of a curve containing 
$z$ so that $z \otimes k(y) \in Z_{\mathfrak b\otimes k(y) , w\otimes k(y)  < \rho(y)}$ specializes to $z$
where $\rho : (Z|_{U \times_W \mathcal N_{T}})_{w, \kappa} \to (\mathcal N_{T})_{w, \kappa} \to \mathrm{Hilb}(B)^Q_\kappa$ is the projection.  If 
$(w < x) \in P'(\kappa)$ such that $z \in Z_{\mathfrak b, w < x}$, then we have 
\[
z\otimes k(y)  \in Z_{\mathfrak b\otimes k(y), w\otimes k(y) < x \otimes k(y) }\cap Z_{\mathfrak b\otimes k(y) , w\otimes k(y)  < \rho(y)}= Z_{\mathfrak b\otimes k(y) , w\otimes k(y)  < x\otimes k(y)  \vee \rho(y)}.
\] 
Thus we must have $x\otimes k(y)  \leq \rho(y)$ because otherwise we have  a generic point $y' \in (Z|_{U \times_W \mathcal N_{T'}})_{w,\kappa}$ with $T < T'$ such that $z\otimes k(y')  \in Z_{\mathfrak b\otimes k(y') , w\otimes k(y')  < \rho'(y')}$ which contradicts with the fact that $z \in J_T$. Since $\iota_R(T)$ is a type of $P$, we conclude that $\iota_R(\rho(y)) \in P$ so that its specialization is the terminal element of the above poscheme, proving the claim.

Now suppose that $\iota_R(T)$ is not a type of $P$.
In this situation, we will prove that we have
\[
H^i_{\text{\'et}, c}(\epsilon^{-1}(J'_T)_{\overline{k}}, \mathbb Z_\ell) = H^i_{\text{\'et}, c}((J'_T\cap \mathrm{im}(Z_{P_U} \to E|_U))_{\overline{k}}, \mathbb Z_\ell) = 0,
\]
for $i \geq 2\mathrm{dim} \, E|_U - I - 2$.
We will only explain the vanishing of the leftmost quantity as it is the more difficult of the two.  We claim that
\[
H^i_{\text{\'et}, c}(\epsilon^{-1}(J'_T)_{\overline{k}}, \mathbb Z_\ell)  = 0.
\]  
By arguing as in the proof of \cref{theo:spectralsequenceforbar}, 
we have a spectral sequence
\[
\bigoplus_{i + j = n} E_1^{i, j} = \bigoplus_{T': \text{ essential, saturated type}} H^{n}_{\text{\'et}, c}((J'_T\cap Z|_{U \times_W \mathcal N_{T'}})_{\overline{k}}, \mu''(T')[1]) \implies H^{n}_{\text{\'et}, c}(\epsilon^{-1}(J'_T)_{\overline{k}}, \mathbb Z_\ell),
\]
where $\mu''(T')$ is defined in a way analogous to the definition of $\mu'(T')$ using the nerve
\[
\epsilon''_{T'} : (N(-\infty, J'_T\cap Z|_{U\times_W \mathcal S_{T'}} )< J'_T\cap Z|_{U \times_W \mathcal N_{T'}}) \to J'_T\cap Z|_{U\times_W \mathcal N_{T'}}
\]
and $T'$ runs over essential saturated combinatorial types of $P$.
Thus it suffices to show that we have
\[
H^{i}_{\text{\'et}, c}((J'_T\cap Z|_{U \times_W \mathcal N_{T'}})_{\overline{k}}, \mu''(T')[1]) = 0,
\]
for $i \geq 2\mathrm{dim} \, E|_U - I - 2$.  This follows from 
\[
H^{i}_{\text{\'et}, c}((J'_T\cap Z|_{U \times_W \mathcal N_{T'}})_{\overline{k}}, R\epsilon''_{T'*}(\epsilon''_{T'})^*\underline{\mathbb Z}_\ell) =  H^{i}_{\text{\'et}, c}((J'_T\cap Z|_{U \times_W \mathcal N_{T'}})_{\overline{k}}, \mathbb Z_\ell) =0,
\]
for $i \geq 2\mathrm{dim} \, E|_U - I - 2$. 
Again we will only prove the vanishing of the leftmost quantity as it is the more difficult of the two.  Since
\[
H^{i}_{\text{\'et}, c}((J'_T\cap Z|_{U \times_W \mathcal N_{T'}})_{\overline{k}}, R\epsilon''_{T'*}(\epsilon''_{T'})^*\underline{\mathbb Z}_\ell)
= H^i_{\text{\'et}, c}((N(-\infty, J'_T\cap Z|_{U\times_W \mathcal S_{T'}} )< J'_T\cap Z|_{U \times_W \mathcal N_{T'}})_{\overline{k}}, \mathbb Z_\ell)
\]
it suffices to prove the right-hand side is $0$ for $i \geq 2\mathrm{dim} \, E|_U - I - 2$. 
If $T' \not\leq T$, then arguing as before using the generic point of a curve through $z$
we conclude $J'_T\cap Z|_{U \times_W \mathcal N_{T'}} = \emptyset$.
So let us assume that $T' \leq T$.
We take a maximal chain of saturated elements
\[
T' \prec T_1 \prec T_2 \prec \dotsb \prec \iota_R(T)
\]
Let $m$ be the minimum index such that $T_m$ is not a type of $P$.
Then it follows from (3) that $\kappa(T_m) > I$.
Now note that the cohomological dimension of $(N(-\infty, J'_T\cap Z|_{U\times_W \mathcal S_{T'}} )< J'_T\cap Z|_{U \times_W \mathcal N_{T'}})_{\overline{k}}$ is at most 
\begin{align*}
2\mathrm{dim}\, E|_U - 2\gamma(T_m) + \mathrm{rank}(T')-2 + 2|\mathrm{Supp}(T_m)| & \leq 2\mathrm{dim}\, E|_U - \kappa(T_m) -2 \\
& < 2\mathrm{dim}\, E|_U - I -2,
\end{align*}
because the cohomological dimension of the fiber $N(w, y)$ over $(J'_T\cap Z|_{U \times_W \mathcal N_{T'}})_{\overline{k}}$ is $\mathrm{rank}(T')-2$  and the dimension of $J'_T\cap Z|_{U \times_W \mathcal N_{T'}}$ is at most $2\mathrm{dim}\, E|_U - 2\gamma(T_m)+ 2|\mathrm{Supp}(T_m)|$.
This completes the proof.
\end{proof}

\subsection{The cohomological virtual bar complex}
 
Suppose that $E \to U$ is a finite rank vector bundle and $Z_{\mathfrak b, w < x} \subset E_w$ is a linear subspace.
As we discussed earlier, the fibers $Z_{\mathfrak{b},w<x}$ need not have the expected dimension.  For this reason, we must truncate the bar complex $B(P_{U},Z_{P_{U}})$ so that we only work with the locus of fibers whose dimensions we understand.  In this section, we define the ``cohomological virtual bar complex'' which records how the bar complex  would behave if every fiber $Z_{\mathfrak{b},w<x}$ had the expected dimension.  Note that we do not attempt to construct an actual space with this property; rather, we define a complex of sheaves which would compute the cohomology of such an object if it existed.  This construction interpolates between the ``actual'' behavior of moduli spaces of curves and the ``desired'' behavior predicted by Manin's Conjecture.

We work in the setting of \cref{subsubsec:stratifyingthespace of sectons}.
Let $T$ be a saturated combinatorial type of $(U_{\mathbf{k}} \leq \mathrm{Hilb}(B)^Q)$.
We fix a Zariski open subset $U \subset \mathcal B \times U_{\mathbf{k}}$.
Then we consider the following nerve
\[
(N(-\infty, U \times_W \mathcal S_T) < U \times_W \mathcal N_T),
\]
whose $m$-th simplex is given by
\[
\{\mathfrak b, w < x_0\leq \dotsb \leq x_m < y \, | \, (\mathfrak b, w) \in U, (w < y)\in \mathcal N_{T}, (w < x_m) \not\in \mathcal S_{T} \},
\]
with the augmentation morphism $\hat{\epsilon}_T :  (N(-\infty, U \times_W \mathcal S_T) < U \times_W \mathcal N_T) \to U \times_W \mathcal N_T$
\[
\hat{\epsilon}_T :  (\mathfrak b, w < x_0\leq \dotsb \leq x_m < y) \mapsto (\mathfrak b, w < y).
\]
We define the complex $\mu(T)[2]$ in $D^b(U \times_W \mathcal N_T)$ by 
\[
 \underline{\mathbb Z_\ell} \to R\hat{\epsilon}_{T*}\hat{\epsilon}_T^*\underline{\mathbb Z}_\ell \to \mu(T)[2] \to \underline{\mathbb Z}_\ell[1] \to.
\]
We define the cohomological virtual bar complex $\mathcal C(U)$ to be
\begin{align*}
\mathcal A_{a, a', \mathbf{k}}(U) &= \bigoplus_{T: \text{ essential type}} \mathcal A_{a, a', \mathbf{k}}(U, T) \\
&= \bigoplus_{T} \bigoplus_i H^i_{\text{\'et}, c}((U \times_W \mathcal N_T)_{\overline{k}}, \mu(T)(-n_{a, a', \mathbf{k}, T})[-2n_{a, a', \mathbf{k}, T}+1]\otimes \mathbb Q_\ell),
\end{align*}
where $n_{a, a', \mathbf{k}, T} = 2a + 2a' +4 -2\sum_i k_i -\gamma(T)$.
When $U =  \mathcal B \times W$, we write $\mathcal A_{a, a', \mathbf{k}} = \mathcal A_{a, a', \mathbf{k}}(U)$ and so on.

From now on we assume that our ground field $k$ is $\mathbb F_q$ and $U = \mathcal B \times U_{\mathbf{k}}$.  The following function plays a crucial role in our analysis:
\begin{defi}[M\"obius function]
Let $U_{\mathbf{k}}$ and $Q^{JB}$ be as before.
Let $w \in U_{\mathbf{k}}(k)$ and $w \leq x \in Q^{JB}(k)$.
Then we define the M\"obius function $\mu_k(w, x)$ using the following formula inductively:
\[
\mu_k(w, w) = 1, \, \mu_k(w, x) = -\sum_{w \leq y < x, y \in Q^{JB}(k)} \mu_k(w, y).
\]
\end{defi}

\begin{lemm}
\label{lemm:Mobius}
Let $w, w_1, w_2 \in U_{\mathbf{k}}(k)$ and $x, x_1, x_2 \in Q^{JB}(k)$ such that $w \leq x, w_1 \leq x_1, w_2 \leq x_2$.
Then we have the following properties:
\begin{enumerate}
\item Suppose that $x_1$ and $x_2$ have disjoint supports. Then we have
\[
\mu_k(w_1\vee w_2, x_1\vee x_2) = \mu_k(w_1, x_1)\mu_k(w_2, x_2).
\]
\item Suppose that $w \leq x$ is not essential. Then we have $\mu_k(w, x) = 0$.
\end{enumerate}
\end{lemm}

\begin{proof}
(1) This assertion follows from the fact that we have $[w_1\vee w_2, x_1\vee x_2] \cong [w_1, x_1]\times[w_2, x_2]$ and \cite[Proposition 1.2.1]{Wachs}.

(2) \Cref{prop:homotopyequivalence} shows that the subposet of essential elements $q$ with $w < q \leq x$ admits a maximal element $\mathrm{ess}(x)$.  
In particular, if $x$ is not essential then $(w,x)$ consists of the union of $(w,\mathrm{ess}(x)]$ with a chain from $\mathrm{ess}(x)$ to $x$, implying the desired result.
\end{proof}

Our goal here is to prove the following proposition:
\begin{prop}
\label{prop:virtualcount}
Let $T$ be an essential saturated type.
Then we have
\[
 \sum_i (-1)^i \mathrm{Tr}(\mathrm{Frob} \curvearrowright \mathcal A_{a, a', \mathbf{k}}^i(T)) = -(\#\Pic^0(B)(k))^2q^{2a+ 2a'+4} \sum_{(w< x) \in \mathcal N_T(k)} \mu_k(w, x)q^{-\gamma(x)},
\]
where $\mathrm{Frob}$ is the Frobenius action. 
\end{prop}
\begin{proof}

We can apply the Grothendieck-Lefschetz trace formula to the perfect complex
\[
R\hat{\epsilon}_{T*}\hat{\epsilon}_T^*\underline{\mathbb Q}_\ell(-n_{a, a', \mathbf{k}, T})[-2n_{a, a', \mathbf{k}, T}];
\]
indeed, this is justified by the comparison theorem \cite[Proposition 5.5.4]{BS15} 
and then applying the standard Grothendieck--Lefschetz formula for $\ell$-adic constructible complexes.
This yields
\begin{multline*}
\sum_i (-1)^i \mathrm{Tr}(\mathrm{Frob} \curvearrowright H^i_{\text{\'et}, c}((\mathcal B\times\mathcal N_T)_{\overline{k}}, R\hat{\epsilon}_{T*}\hat{\epsilon}_T^*\underline{\mathbb Q}_\ell(-n_{a, a', \mathbf{k}, T})[-2n_{a, a', \mathbf{k}, T}]))\\
= \sum_{(\mathfrak b, (w < x)) \in (\mathcal B \times \mathcal N_T)(k)} q^{2a+ 2a'+4 - \gamma(x)}\mathrm{Tr}(\mathrm{Frob} \curvearrowright (R\hat{\epsilon}_{T*}\hat{\epsilon}_T^*\underline{\mathbb Q}_\ell)_{(\mathfrak b, (w < x))})
\end{multline*}
where we have used the equality $\gamma(T) = \gamma(x) - 2\sum_i k_i$.
Now the fiber of the nerve
\[
\hat{\epsilon}_T :  (N(-\infty, \mathcal B \times \mathcal S_T) < \mathcal B \times \mathcal N_T) \to \mathcal B \times \mathcal N_T,
\]
at $(\mathfrak b, (w, x))$ is the nerve
\[
N(w < x) \to \Spec \, k.
\]
Then it follows from \cref{lemm:deformationretract} that
\begin{align*}
\mathrm{Tr}(\mathrm{Frob} \curvearrowright (R\hat{\epsilon}_{T*}\hat{\epsilon}_T^*\underline{\mathbb Q}_\ell)_{(\mathfrak b, (w < x))} &= \sum_i (-1)^i\mathrm{Tr}(\mathrm{Frob}\curvearrowright H^i_{\text{\'et}}(N(w < x)_{\overline{k}}, \mathbb Q_\ell))\\
&= \sum_i (-1)^i\mathrm{Tr}(\mathrm{Frob}\curvearrowright H^i_{\text{\'et}}(N((w < x)\cap Q^{JB})_{\overline{k}}, \mathbb Q_\ell)).
\end{align*}
Let $N_n(w, x) = \# \{ x_0 < \dotsb < x_n \in (N(w, x)\cap Q^{JB})(k)\}$. Then the Grothedieck-Lefschetz trace formula for simplicial schemes 
implies that we have 
\[
\sum_i (-1)^i\mathrm{Tr}(\mathrm{Frob}\curvearrowright H^i_{\text{\'et}}(N((w < x)\cap Q^{JB})_{\overline{k}}, \mathbb Q_\ell)
= \sum_i (-1)^i N_i(w, x).
\]
(This identity follows from \cref{Simplicial_Model} by computing the trace on the normalized cochain complex of the semisimplicial set $N((w < x)\cap Q^{JB})_{\overline{k}}$.)   Then it follows from \cite{Hall34} (see also \cite[(2.5) in p. 559]{Gree82}) that we have
\[
\sum_i (-1)^i N_i(w, x) = 1 + \mu_k(w, x).
\]
By Lang's theorem, the torsor $\Pic^d(B)$ is isomorphic to $\Pic^0(B)$. Thus when we sum over all $k$-points $(\mathfrak b, (w < x))$ we obtain $(\#\Pic^0(B)(k))^2$ times the sum over all $k$-points $(w<x)$.  We conclude that
\begin{multline*}
\sum_i (-1)^i \mathrm{Tr}(\mathrm{Frob} \curvearrowright H^i_{\text{\'et}, c}((\mathcal B\times\mathcal N_T)_{\overline{k}}, R\hat{\epsilon}_{T*}\hat{\epsilon}_T^*\underline{\mathbb Q}_\ell(-n_{a, a', \mathbf{k}, T})[-2n_{a, a', \mathbf{k}, T}]))\\
= (\#\Pic^0(B)(k))^2\sum_{(w < x) \in \mathcal N_T(k)} q^{2a+ 2a'+4 - \gamma(x)}(1 + \mu_k(w, x))
\end{multline*}
Finally we consider the triangle: 
\begin{multline*}
 \underline{\mathbb Q_\ell}(-n_{a, a', \mathbf{k}, T})[-2n_{a, a', \mathbf{k}, T}] \to R\hat{\epsilon}_{T*}\hat{\epsilon}_T^*\underline{\mathbb Q}_\ell(-n_{a, a', \mathbf{k}, T})[-2n_{a, a', \mathbf{k}, T}] \to \\
 \mu(T)(-n_{a, a', \mathbf{k}, T})[-2n_{a, a', \mathbf{k}, T}+2]\otimes \mathbb Q_\ell \to \underline{\mathbb Q}_\ell(-n_{a, a', \mathbf{k}, T})[-2n_{a, a', \mathbf{k}, T}+1].
\end{multline*}
The additivity of the alternating sums of the Frobenius traces in exact sequences proves our assertion.
\end{proof}

\subsection{Convergence}
In this section we obtain exponential bounds of the error term for cohomological virtual bar complex.
We denote by $\langle d \rangle$ the operation $(d)[2d]$ of shifting by $2d$ and twisting by $d$.

Given the poset $(U_{\mathbf{k}} \leq \mathrm{Hilb}(B)^Q)$, we have the 
cohomological virtual bar complex  associated to $\alpha = (a, a', k_1, \dotsc, k_r)$:  
$$\mathcal A_{\alpha} := \bigoplus_{T: \text{ essential type} } 
H^*_{\text{\'et}, c}(\mathcal N_{T, \overline{k}}, \mu(T)[1]\otimes\mathbb Q_\ell \langle -{E}(\alpha) + 2\sum_{i} k_i(\alpha)+\gamma(T)\rangle),$$ where $E(\alpha)= 2(a + a') + 4 $ is  the 
dimension of the space of sections of the bundle $\mathcal O(a)^{\oplus 2} \oplus \mathcal O(a')^{\oplus 2}$ over $\mathbb P^1$.

We wish to bound the $\mathsf L^1$-trace of Frobenius on the truncation into cohomological degrees $\leq  2{\rm edim}(\alpha) - I$,  where ${\rm edim}(\alpha) = -K_S.\alpha + 4$ and 
$I$ grows linearly in the degree of $\alpha$.    More precisely, we want to show that the sum of the absolute values of the generalized eigenvalues of the Frobenius is $o(q^{{\rm edim}(\alpha)})$.     
Since shifting homological degrees does not affect $\mathsf L^1$-traces,  it suffices to prove that the $\mathsf L^1$-trace of Frobenius on the truncation of the graded vector space 
\begin{multline*}
\mathcal A_{\alpha}\langle {\rm edim}(\alpha) \rangle\\
= \bigoplus_{T:  \text{ ess. type }}
H^*_{\text{\'et}, c}(\mathcal N_{T, \overline{k}}, \mu(T)[|\mathrm{Supp}(T)|]\otimes \mathbb Q_\ell \langle  \gamma(T) +2\sum_i k_i(\alpha) +  {\rm edim}(\alpha)  - E(\alpha)\rangle)
\end{multline*}
to cohomological degrees $\leq - I - |\mathrm{Supp}(T)|$ is $o(1)$.   Note that ${\rm edim}(\alpha) - E(\alpha)+2\sum_ik_i(\alpha)$ is equal to $\sum_{i} k_i$, and so this space is independent of $a, a'$.   

Accordingly, we rename it $\mathcal C_{\mathbf{k}, T}$, and summing over all  $\mathbf{k} = (k_1, \dotsc, k_r)$ and $T$ we obtain a bigraded vector space $\mathcal C := \bigoplus_{\mathbf{k}, T} \mathcal C_{\mathbf{k}, T}$ where  
the first grading is given by $\sum_{i = 1}^r k_i + |\mathrm{Supp}(T)|$ and the second grading is the homological grading.   We may decompose $\mathcal C$ into a sum over saturated types $T$ of $(\mathrm{Hilb}(B)^Q\leq\mathrm{Hilb}(B)^Q)$, where the summand associated to $T$ is either zero or $$H^*_{\text{\'et}, c}(\mathcal N_{T, \overline{k}}, \mu(T)[|\mathrm{Supp}(T)|]\otimes \mathbb Q_\ell (  \gamma(T) + k(T))[2(\gamma(T) + k(T) )],$$ placed in grade $k(T) + |\mathrm{Supp}(T)|$ where $k(T) = \sum_{i = 1}^rk_i$.

Recall that an essential saturated type $T$ corresponds to a multiset of elements in $(\mathrm{ch}(Q)\leq\mathrm{ch}(Q))$.  Suppose $T$ contains the elements $$(f^{w_1}\leq f^{x_1}), \dotsc, (f^{w_t} \leq f^{x_t}),$$ where $(f^{w_i} \leq f^{x_i})$ occurs with multiplicity $m_i$.  Then we have
$$\mathcal N_T \cong \mathrm{Conf}_{m_1 + \dotsb + m_t} \mathbb P^1/(S_{m_1} \times \dotsb \times S_{m_t}).$$
We have that the locally constant pro-\'etale sheaf on $\cN_T$
$$\mathbb H^*(\mu(T)[|\mathrm{Supp}(T)|]) := \bigoplus_{i} \mathbb H^i(\mu(T)[|\mathrm{Supp}(T)|])$$ 
pulls back to the constant local system $$\bigotimes_{i = 1}^t  H^*(\mu(w_i, x_i)[1])^{\otimes m_i},$$ on $\mathrm{Conf}_{m_1 + \dotsb + m_t} \mathbb P^1$,  which is $S_{m_1} \times \dotsb  \times S_{m_t}$-equivariant with respect to the Koszul sign rule.   Hence by Galois descent, we have 
\[\begin{multlined}[\displaywidth]
H^*_{\text{\'et}, c} (\mathcal N_{T, \overline{k}}, \mathbb H^*(\mu(T)[|\mathrm{Supp}(T)|]) \otimes \mathbb Q_\ell\langle  \gamma(T) + k(T)\rangle) = \\
 \left(H^*_{\text{\'et}, c}(\mathrm{Conf}_{m_1+ \dotsb + m_t} \mathbb P^1_{\overline{k}}, \mathbb Q_\ell) \otimes \bigotimes_{i = 1}^t  (H^*(\mu(w_i, x_i)[1]\otimes \mathbb Q_\ell)\langle \gamma(w_i \leq x_i) + k(w_i) \rangle  )^{\otimes m_i} \right)^{S_{m_1} \times \dotsb \times S_{m_t}}
 \end{multlined}\]
   where $k(w)$ denotes the length of $w$.  Note that to compute the $\mathsf L^1$-trace we may replace ${S_{m_1} \times \dotsb \times S_{m_t}}$-invariants by coinvariants since we are using characteristic zero coefficients. 
   
   Accordingly, for an essential pair $$(f^w \leq f^x) \in (\mathrm{ch}(Q)\leq \mathrm{ch}(Q)),$$ such that $f^w = k_i [\ell_{i,1} \oplus \ell_{i,2}]$, we define $$L(w \leq x) :=  H^*(\mu(w, x)[1]\otimes \mathbb Q_\ell) \langle  \gamma(w \leq x)  + k(w)\rangle.$$    We let $\mathcal J$ be the set of all such essential pairs excluding  
   $(V \leq V)$ and define $$L = \bigoplus_{w \leq x \in \mathcal J} L(w \leq x),$$ where the sum ranges over all such pairs.  We have that the first grading for $L$ is given by placing $L(w \leq x)$ in degree $k(w) + 1$ and the second grading is by homological degree.
   Then by taking the direct sum over $T$, we have that $\mathcal C$ is isomorphic to 
 \begin{multline*}
    \bigoplus_{m: \mathcal J \to \mathbb N \text{ fin. supp.}} \left(H^*_{\text{\'et}, c}(\mathrm{Conf}_{\sum_{j \in \mathcal{J}} m(j)} \mathbb P^1_{\overline{k}}, \mathbb Q_\ell) \otimes \bigotimes_{j \in \mathcal J} L(j)^{\otimes m(j)}\right)_{\prod_{j\in J} S_{m(j)}}  \\
    =  \bigoplus_{m \geq 0}  (H^*_{\text{\'et}, c}(\mathrm{Conf}_m \mathbb P^1_{\overline{k}}, \mathbb Q_\ell) \otimes L^{\otimes m})_{S_m}  .
 \end{multline*}

   We set $r = 4$ and apply \cref{CombinedBound}: we claim that the $\mathsf L^1$-trace of Frobenius on $L_{n,i}$ is bounded by $E d^n b^i$ and the $\mathsf L^1$-trace of Frobenius on $L_{\bullet, 2}$ is bounded by $D$ where 
   $$D = 4, d = 1, b= q^{-1/2}, \text{ and } E = c_1$$ 
   for some constant $c_1$ which is uniform. Indeed, for $L(w \leq x)$ the maximal cohomological degree is bounded above by $-2$ 
   and the minimum cohomological degree is bounded below by $-2(\gamma(w \leq x) + k(w))$. Note that $H^*(\mu(w,w))$ is the reduced cohomology of the empty set, so 
   only the $(-1)$-degree cohomology is non-trivial.  
As we vary over $(w \leq x) \in \mathcal{J}$, the posets $(w, x)$ only have finitely many isomorphism types.  Furthermore the number of pairs which contribute to $L_{n,i}$ is uniformly bounded.  Since everything is twisted by $q^{-(\gamma(w\leq x) + k(w))}$, one can find a uniformly bounded $c_{1}$ (using the dimensions of $H^i(\mu(w, x)[1]\otimes \mathbb Q_\ell)$ where $(w < x)$ are essential pairs, we can verify that $c_1 = 8$ works) such that $E = c_{1}$, $d = 1$, and $b= q^{-1/2}$ yield the bound on $|{\rm Frob},  L_{n,i}|$. 
   The part $L$ lying in cohomological degree $-2$ only comes from $w = x = [\ell_{i,1} \oplus \ell_{i,2}]$ and since there are four possible choices of $[\ell_{i,1} \oplus \ell_{i,2}]$ we see that $D = 4$ works as the desired upper bound for $|{\rm Frob},  L_{\bullet,2}|$.   Thus we conclude from \cref{CombinedBound} that 
   $$|{\rm Frob}, (\mathcal C)_{n,i}|  =O(i{d'}^n (4/\sqrt{q})^i),$$ for any $d' > \max\{4, 10c_1 \}$ independent of $q$. 
   (Here in the $O$-term we have absorbed the factor $n$ from \cref{CombinedBound} into the exponential term $d'^{n}$ at the cost of increasing $d'$ slightly.) Shifting, twisting back, and applying geometric series we obtain:
   
   \begin{theo}
   \label{theo:errorforvirtual}
   	Let $\alpha = (a, a', k_1, \dotsc, k_4)$ and fix any $0 < \eta < 1$.  
    For any constant $d' > 80$ we have that if $q$ is large enough so that $d'(4/\sqrt{q})^{(1-\eta)} < 1$ then for any $I > 0$ we have
	$$|{\rm Frob}, \tau_{\leq {2\rm edim}(\alpha)- I} \mathcal A_\alpha| =O\left(q^{{\rm edim}(\alpha)} (d')^{\sum_{i = 1}^4 k_i} (4/{\sqrt{q}})^{(1-\eta)I} \right),$$
    where $\tau$ is the truncation functor and the implied constants are independent of $\alpha, q$.
  
   \end{theo}

\section{Cohen-Jones-Segal Conjecture}
\label{sec:CJS}

In this section we prove \cref{theo:introcjs}.  We will adhere closely to the argument of \cite[Sections 6-9]{DT24}.  We prove the Cohen-Jones-Segal conjecture through a sequence of comparisons:
\begin{enumerate}
\item we compare $M_{\alpha,*}$ to the semi-topological model of \cite[Sections 6-7]{DT24};
\item we compare the semi-topological model to the space of pointed positive  continuous maps of \cite[Introduction]{DT24}; 
\item we compare the space of pointed positive continuous maps to the space of pointed continuous maps.
\end{enumerate}
Throughout we will focus on the pointed versions of moduli spaces of curves (although the arguments also work in the unpointed case as described in \cite{DT24}). In this section, we work over $\mathbb C$ with the Euclidean topology and all (co)homologies are singular (co)homologies.

\subsection{Comparing algebraic maps to semi-topological models}

One important contribution of \cite{DT24} is the introduction of semi-topological models for the space of maps.  Loosely speaking, these semi-topological models will parametrize continuous maps satisfying certain constraints on orders of vanishing near specified points.  In this section we compare the space of algebraic maps to a semi-topological model.  The key ingredient is our approximation result \cref{theo:algebraicvssemitopological}.

\subsubsection{Semi-topological models}

\begin{defi}[{\cite[Definition 6.1]{DT24}}]
Let $U \subset \mathbb C$ be an open subset with a coordinate $z$ and let $w : \mathfrak K \to \mathrm{Hilb}^{[n]}(U)$ be a continuous
family of divisors of $U$ parametrized by a compact set $\mathfrak K$. 
Fix a continuous family of polynomials $\rho_{w(\mathfrak k)}(z)$ such that for every $\mathfrak k \in \mathfrak K$, the divisor of zeros for $\rho_{w(\mathfrak k)}(z)$ coincides with the divisor $w(\mathfrak k)$.
We say that a family of functions $f_{\mathfrak k} : U \to \mathbb C, \mathfrak k \in \mathfrak K$ vanishes to order $\geq w$ if for every $u \in U$ and $\mathfrak k \in \mathfrak K$, there exists a constant $a_{u, \mathfrak k} \in \mathbb C$ such that $f_{\mathfrak k}(z) = a_{u, \mathfrak k}\rho_{w(\mathfrak k)}(z) + o(|\rho_{w(\mathfrak k)}(z)|)$ as $z \to u$ where $o(|\rho_{w(\mathfrak k)}(z)|)$ is uniform in $\mathfrak K$. If $a_{u, \mathfrak k} \neq 0$ for all $u \in U, \mathfrak k \in \mathfrak K$, we say that $f$
 vanishes to order exactly $w$.

These notions are independent of the choice of polynomials $\rho_{w(\mathfrak k)}(z)$ and the choice of coordinate $z \in U \subset \mathbb C$ by \cite[Proposition 6.2.(2)]{DT24} and \cite[The proof of Proposition 6.4]{DT24}.
\end{defi}

\begin{rema}
Even when $f$ vanishes to order $\geq w$ there may not be a $w'$ with $w_{\mathfrak k}' \geq w_{\mathfrak k}$ for every $\mathfrak k \in \mathfrak K$ such that $f$ vanishes exactly to order $w'$.
\end{rema}

Next suppose we fix a Riemann surface $B$ and a family of divisors on $B$ parametrized by a continuous map $w : \mathfrak K \to \mathrm{Hilb}^{[n]}(B)$ where $\mathfrak K$ is a compact set. For each $\mathfrak k_0 \in \mathfrak K$, we can choose
\begin{itemize}
\item a neighborhood of $w(\mathfrak k_0)$ of the form $\prod_{i = 1}^s \mathrm{Hilb}^{[n_i]}(U_i)$ where the $U_i$'s are disjoint open subsets of $B$ which are biholomorphic to open subsets in $\mathbb C$ and $n = \sum_{i = 1}^s n_i$; and
\item an open neighborhood  $V$ of $\mathfrak k_0$ in $\mathfrak K$ such that $w|_V$ factors through $\prod_{i = 1}^s\mathrm{Hilb}^{[n_i]}(U_i)$ with $i$th component map denoted by $w_i: V \to \mathrm{Hilb}^{[n_i]}(U_i)$.
\end{itemize}
We say that $f : \mathfrak K \times B \to \mathbb C$ is a family of functions which vanish to order $\geq w$ if for every $\mathfrak k_0 \in \mathfrak K$ and for some (equivalently every) choice of neighborhoods $U_i$'s and $V$ as above, we have the property that for each $i$ and for every compact subset $\mathfrak K' \subset V$ the restriction $f|_{\mathfrak K' \times U_i} \to \mathbb C$ vanishes to order $\geq w_{i}$.  We define a family of functions which vanish to order exactly $w$ in the analogous way. 

\begin{defi}[{\cite[Definition 6.5]{DT24}}]
Let $M$ be a complex manifold and $N \subset M$ be a complex submanifold of codimension $r$.
Let $B$ be a Riemann surface and let $w: \mathfrak K\to \mathrm{Hilb}^{[n]}(B)$ be a family of divisors parametrized by a compact subset $\mathfrak K$. Suppose that we have a continuous map
\[
f : \mathfrak K \times B \to M,
\]
such that for any $\mathfrak k \in \mathfrak K$ and any point $p$ in the support of $w(\mathfrak k)$, we have $f_{\mathfrak k}(p) \in N$.
For any $\mathfrak k_0 \in \mathfrak K$, one can obtain a collection of data:
\begin{itemize}
\item for every $p \in \mathrm{Supp}(w(\mathfrak k_0))$, there exists a neighborhood $\mathcal U_p \subset M$ of $f_{\mathfrak k_0}(p) \in N$ and holomorphic local coordinates $s_1^p, \dotsc, s_r^p$ on $\mathcal U_p$ whose vanishing locus is $N\cap \mathcal U_p$;
\item for the various $p \in \mathrm{Supp}(w(\mathfrak k_0))$ we can choose disjoint closed balls $\mathbb B_{p} \subset B$ with center $p$ and open subsets $U_{p} \subset \mathbb B_{p}$ containing $p$ such that $w(\mathfrak k_{0})$ has an open neighborhood of the form $\prod_{p} \mathrm{Hilb}^{[n_p]}(U_p)$ where $n_p$ is the multiplicity of $w(\mathfrak k_0)$ at $p$, and;

\item there is an open neighborhood $V \subset \mathfrak K$ of $\mathfrak k_0$ such that for any $\mathfrak k \in V$, we have $f_{\mathfrak k}(\mathbb B_p) \subset \mathcal U_p$ and $w|_V$ factors through $\prod_{p} \mathrm{Hilb}^{[n_p]}(U_p)$ with component maps $w_{p}: V \to \mathrm{Hilb}^{[n_p]}(U_p)$. 
\end{itemize}
We say that $f : \mathfrak K \times B \to M$ intersects $N$ to order $\geq w$ if for some (equivalently every) choice of $\mathcal U_p, \mathbb B_p, U_p, V$ and $s_1^p, \dotsc, s_r^p$, and every compact subset $\mathfrak K' \subset V$ , $s_i^{p}\circ (f|_{\mathfrak K' \times U_p})$ 
vanishes to order $\geq w_p$ for any $p$ and $i$.

Furthermore we say that $f$ intersects $N$ to order exactly $w$ if additionally for any $\mathfrak k \in \mathfrak K$, the support of $f_{\mathfrak k}^{-1}(N)$ is equal to the support of $w(\mathfrak k)$, and for every $p$, there exists $i$ such that $s_i^{p}\circ(f|_{\mathfrak K' \times U_p})$  
vanishes to order exactly $w_p$.

Finally, we say that $f$ intersects $N$ to order $>w$ if it intersects $N$ to order $\geq w$ but not to order exactly $w$.  (Note that this does not imply that $f$ intersects $N$ to order $\geq w'$ for some $w'$ such that $w'_{\mathfrak k} \geq w_{\mathfrak k}$ for every $\mathfrak k$.)

\end{defi}

Let $M$ be a complex manifold and $N_1, \dotsc, N_r \subset M$ be complex submanifolds. Consider the set
\[
\mathfrak X = \{ (w_i)_{i = 1}^r \in \prod_{i = 1}^r\mathrm{Hilb}^{[n_i]}(B), f \in \mathrm{Top}(B, M) \, | \, \text{ for any $i$, $f$ intersects $N_i$ to order $\geq w_i$.} \}.
\]
By \cite[Proposition 6.9]{DT24} there is a unique compactly generated topology on $\mathfrak X$ such that for all compact Hausdorff spaces $\mathfrak K$, continuous maps $f: \mathfrak K \to \mathfrak X$ are in bijection with pairs of continuous maps $f : \mathfrak K \times B \to M$ and $w_i : \mathfrak K \to \mathrm{Hilb}^{[n_i]}(B)$ such that $f$ intersects $N_i$ to order $\geq w_i$ for any $i = 1, \dotsc, r$.

\subsubsection{Applying to blow-ups of $\mathbb{P}^{1} \times \mathbb{P}^{1}$}
\label{subsubsec:oursituation}

We consider the setup described in \cref{subsec:blowups}.
Let $U_{\mathbf{k}} \subset \prod_{i = 1}^r\mathrm{Hilb}^{[k_i]}(B)$ be the open set parametrizing tuples $(w_i)$ of effective divisors with disjoint supports. Consider the map $E^{\mathrm{top}} \to \mathcal{B} \times U_{\mathbf{k}}$ where $E^{\mathrm{top}}$ is
\[
E^{\mathrm{top}} = \{ ((L_1, L_2), w, \text{ sections $B \to \oplus_{j = 1}^2 L_j\otimes V_j$ intersecting $\oplus_{j = 1}^2 L_j \otimes \ell_{i, j}$ to order $\geq w_i$ for any $i$} )\}.
\]
Arguing as in \cite[Proposition 6.12]{DT24}, one can identify $E^{\mathrm{top}}$ as a locally trivial bundle of infinite dimensional separable Banach spaces over $\mathcal B \times U_{\mathbf{k}}$.  We define a stratification of this topological bundle:

\begin{defi}[{\cite[Definition 6.13]{DT24}}]
For each $(L = (L_1, L_2), (w\leq x)) \in \mathcal B \times (U_{\mathbf{k}} \times \mathrm{Hilb}(B)^Q)$, let $Z_{L, w\leq x}$ be the subspace of $E_{L, w}$ consisting of sections $s \in E_{L, w}$ such that for all $c \in B$, we have
\begin{itemize}
\item if the length of $w_{\ell_{i, 1} \oplus \ell_{i, 2}}$ at $c$ is $l$ and the length of $x_{\ell_{i, 1} \oplus \ell_{i, 2}}$ at $c$ is greater than $l$, then $s$ intersects $\oplus_{j = 1}^2 L_j \otimes \ell_{i, j}$ to order $> l$;
\item if $x_{V_j}$ is supported at $c$, then our section $s_j$ vanishes at $c$ where $s_j$ is the section $s_j: B \to L_j \otimes V_j$ coming from $s$.  
\end{itemize}
\end{defi}
It is easy to see that we have (i) $Z_{L, w\leq x_1} \supset Z_{L, w \leq x_2}$ when $x_1 \leq x_2$, (ii) $Z_{L, w \leq x_1\vee x_2} = Z_{L, w \leq x_1} \cap Z_{L, w \leq x_2}$, and (iii) $Z_{L, w \leq x} = Z_{L, w \leq \mathrm{ess}(x)}$ when $x$ is saturated. Using the description of the saturation function in \cref{exam:saturatedexample}, one may also prove that $Z_{L, w \leq x} = Z_{L, w \leq \mathrm{sat}(x)}$.  Finally by arguing as in \cite[Proposition 6.14]{DT24}, one can prove that the subspace
\[
Z^{\mathrm{top}} =\{ L = (L_1, L_2), (w\leq x) \in (U_{\mathbf{k}}\leq \mathrm{Hilb}(B)^Q), s \in Z_{L, w \leq x}\}
\]
is a closed subset in $\mathcal B \times  (U_{\mathbf{k}}\leq \mathrm{Hilb}(B)^Q) \times_{\mathcal B \times U_{\mathbf{k}}} E^{\mathrm{top}}$ so that $Z^{\mathrm{top}}$ is a stratification of $E^{\mathrm{top}}$.

\subsubsection{Comparing algebraic maps and semi-topological models}
\label{subsec:algebraicvssemitopological}

Fix base points $*_B \in B$ and $*_S \in S$. We assume that $*_S$ is a general point on $S$. Let $Q$ be the poset described in \cref{subsubsec:stratifyingthespace of sectons}.
Let $Q_{r+1}\supset Q$ be the poset as in \cref{subsubsec:stratifyingthespace of sectons}, but with the additional $2$-dimensional vector space $\ell_{r + 1, 1}\oplus \ell_{r + 1, 2}$ corresponding to $*_S$.
Then we define
\[
U_{\mathbf{k}, *} \subset U_{\mathbf{k}} \times \{*_B\},
\]
as the Zariski open subset of pairs 
$(w, *_B)$ with disjoint supports.
Let $\overline{R} \subset \mathrm{Hilb}(B)^{Q_{r + 1}}$ be the closed subposet consisting of divisors $D$ such that $D_{\ell_{r + 1, 1}\oplus \ell_{r + 1, 2}} - D_{0}$ is supported on $*_B$ with multiplicity at most $1$ and $D_{\ell_{r + 1, j}} = D_0$ for $j = 1, 2$.
We define $R = (U_{\mathbf{k}, *} < \overline{R})$.
Note that $R$ is initial in   
$(U_{(\mathbf{k}, 1)} < \mathrm{Hilb}(B)^{Q_{r+1}})$, and this is one of the main reasons why the proof applies in this setting.

Let $\widetilde{\mathcal M}_{a, a', \mathbf{k}, *}$ be the space parametrizing tuples $L = (L_1, L_2) \in \mathcal B = \Pic^a(B) \times \Pic^{a'}(B)$, $w \in U_{\mathbf{k}, *}$, and a continuous section $s = (s_1, s_2) \in \Gamma^{\mathrm{top}}(B, L_1 \otimes V_1 \oplus L_2 \otimes V_2)$ intersecting to $\ell_{i, 1}\oplus \ell_{i, 2}$ to order exactly $w_i$ for any $i = 1, \dotsc, r$ and to order $\geq w_i$ for $i = r+1$.
This is an open subset of our topological bundle $E^{\mathrm{top}}_{U_{\mathbf{k}, *}}$. 
The $2$-dimensional torus $(\mathbb C^\times)^2$ acts on this space by acting as $(\lambda_1, \lambda_2)\cdot (L, w, s = (s_1, s_2)) = (L, w, (\lambda_1 s_1, \lambda_2 s_2))$, and we let
\[
\mathcal M_{a, a', \mathbf{k}, *} = \widetilde{\mathcal M}_{a, a', \mathbf{k}, *}/(\mathbb C^\times)^2.
\]
We also define $\widetilde{M}_{a, a', \mathbf{k}, *}$ as the space parametrizing tuples $L = (L_1, L_2) \in \mathcal B$, $w \in U_{\mathbf{k}}$, and a algebraic section $s = (s_1, s_2) \in \Gamma^{\mathrm{alg}}(B, L_1 \otimes V_1 \oplus L_2 \otimes V_2)$ intersecting to $\ell_{i, 1}\oplus \ell_{i, 2}$ to order exactly $w_i$ for any $i = 1, \dotsc, r$ and to order $\geq w_i$ for $i = r+1$. This variety is a Zariski open set of the algebraic fibration $E^{\mathrm{alg}}_{U_{\mathbf{k}, *}} \subset E^{\mathrm{top}}_{U_{\mathbf{k}, *}}$ over $\mathcal B \times U_{\mathbf{k}, *}$ parametrizing algebraic sections, and it admits an analogous algebraic stratification $Z^{\mathrm{alg}}_{R} = Z^{\mathrm{top}}_{R} \times_{E^{\mathrm{top}}} E^{\mathrm{alg}}$.  We also have
\[
M_{a, a', \mathbf{k}, *} \cong \widetilde{M}_{a, a', \mathbf{k}, *}/(\mathbb C^\times)^2.
\]

Before stating the main theorem of this section, let us recall the following definition:

\begin{defi}
   Let $I$ be a non-negative integer and $f : X \to Y$ be a continuous map of topological spaces. We say $f$ is homology $I$-connected if for any $i < I$, $f$ induces an isomorphism
   \[
   H_i^{\mathrm{sing}}(X, \mathbb Z) \cong H_i^{\mathrm{sing}}(Y, \mathbb Z),
   \]
   and when $i = I$, $f$ induces an surjection
   \[
   H_i^{\mathrm{sing}}(X, \mathbb Z) \hookrightarrow H_i^{\mathrm{sing}}(Y, \mathbb Z).
   \]
\end{defi}

The following theorem is essentially \cite[Theorem 7.4]{DT24}, but it is adapted to our situation:

\begin{theo}
\label{theo:algebraicvssemitopological}
Let $I \in \mathbb Q_{\geq 0}$. Let $P \subset R$ be a poset which is a closed and downward closed union of finitely many combinatorial types such that
\begin{itemize}
\item $P$ is proper over $U_{\mathbf{k}, *}$;
\item there exists a Zariski proper closed subset $F \subset U_{\mathbf{k}, *}$ such that the real codimension of $F$ in $U_{\mathbf{k}, *}$ is $> I$ and the real codimension of $\pi^{-1}(F)$ in $\mathrm{Mor}_*(B, S, (a, a', \mathbf{k}))$ is also $> I$ where $\pi : \mathrm{Mor}_*(B, S, (a, a', \mathbf{k}))  \to U_{\mathbf{k}, _*}$ is the projection. We set $V_{\mathbf{k}, *} = U_{\mathbf{k}, *}\setminus F$;
\item For every pair $w < x \in P\cap (V_{\mathbf{k}, *} \leq \mathrm{Hilb}(B)^Q)$ and every essential $w < y \in (V_{\mathbf{k}, *} \leq Q^{JB})$ such that $x \prec y$, $\Gamma^{\mathrm{alg}}(B, L_1\otimes V_1 \oplus L_2\otimes V_2)_y$ is unobstructed for every $L = (L_1, L_2) \in \mathcal B$, and;
\item $P$ contains all types $T$ with $\kappa(T)\leq I$ where $\kappa(T)$ is defined by \cref{eq:kappadef},
and all minimal types of $U_{\mathbf{k}} < \mathrm{Hilb}(B)^Q$, i.e., all $x \in \mathrm{Hilb}(B)^Q(\overline{k})$ such that $w \prec x$ for some $w \in U_{\mathbf{k}}(\overline{k})$. 
\end{itemize}
We also assume that $M_{a, a', \mathbf{k}, *}$ is smooth.
Then the map $M_{a, a', \mathbf{k}, *} \hookrightarrow \mathcal M_{a, a', \mathbf{k}, *}$ is homology $\lfloor I\rfloor$-connected.
\end{theo}

First we record the following lemma:

\begin{lemm}
If $\widetilde{M}_{a, a', \mathbf{k}, *} \hookrightarrow \widetilde{\mathcal M}_{a, a', \mathbf{k}, *}$ is homology $\lfloor I \rfloor$-connected, then the map 
$$
M_{a, a', \mathbf{k}, *} \hookrightarrow \mathcal M_{a, a', \mathbf{k}, *}
$$ 
is also homology $\lfloor I \rfloor$-connected.
\end{lemm}
\begin{proof}
For notational brevity let $G = (\mathbb{C}^\times)^2$.
Let $EG$ (as usual) be a contractible space with a free $G$ action and set $BG := EG/G$, e.g. $EG = (\mathbb{C}^\infty - 0)^2$,  making $BG = (\mathbb{C}P^\infty)^2$.
Since $\widetilde{M}_{a, a', \mathbf{k}, *} \to M_{a, a', \mathbf{k}, *}$ and $\widetilde{\mathcal{M}}_{a, a', \mathbf{k}, *} \to \mathcal{M}_{a, a', \mathbf{k}, *}$ are principal $G$-bundles, we can replace the quotients by the associated Borel constructions (which have the same homotopy type) to get a map of fiber bundles
\[(\widetilde{M}_{a, a', \mathbf{k}, *} \times EG)/G \to (\widetilde{\mathcal{M}}_{a, a', \mathbf{k}, *} \times EG)/G\]
over $BG$ with respective fibers $\widetilde{M}_{a, a', \mathbf{k}, *}$ and $\widetilde{\mathcal{M}}_{a, a', \mathbf{k}, *}$.
Now the required isomorphism follows from comparing the corresponding Serre spectral sequences.
\end{proof}

\begin{proof}[Proof of  \cref{theo:algebraicvssemitopological}:]
By the previous lemma it suffices to prove that $\widetilde{M}_{a, a', \mathbf{k}, *} \hookrightarrow \widetilde{\mathcal M}_{a, a', \mathbf{k}, *}$ is homology $\lfloor I \rfloor$-connected.
First note that $\widetilde{M}_{a, a', \mathbf{k}, *} \hookrightarrow \widetilde{\mathcal M}_{a, a', \mathbf{k}, *}$ is a $\mathcal B\times U_{\mathbf{k}, *}$-map. By comparing two Leray spectral sequences over $\mathcal B\times U_{\mathbf{k}, *}$ it suffices to show that for a basis of open sets $U \subset \mathcal B\times U_{\mathbf{k}, *}$, we have that
\[
\widetilde{M}_{a, a', \mathbf{k}, *}|_U \hookrightarrow \widetilde{\mathcal M}_{a, a', \mathbf{k}, *}|_U,
\]
is homology $\lfloor I \rfloor$-connected.

Next we consider the following diagram

\centerline{
\xymatrix{ 
\widetilde{M}_{a, a', \mathbf{k}, *}|_{U\cap V_{\mathbf{k}}}  \ar[d] \ar[r] & \widetilde{\mathcal M}_{a, a', \mathbf{k}, *}|_{U \cap V_{\mathbf{k}}} \ar[d] \\
\widetilde{M}_{a, a', \mathbf{k}, *a}|_U  \ar[r]& \widetilde{\mathcal M}_{a, a', \mathbf{k}, *}|_U
}
}
It follows from our assumptions that the vertical maps are homology $\lfloor I \rfloor$-connected. 
Indeed, the assertion for the left map follows from our smoothness and codimension assumptions.
The assertion for the right map follows from a similar argument after using \cite[Proposition 6.16]{DT24} to pass the computation to a finite rank subbundle of $E^{\mathrm{top}} \to {U_{\mathbf{k}, *}}$.
Thus it suffices to show that $\widetilde{M}_{a, a', \mathbf{k}, *}|_{U\cap V_{\mathbf{k}}}  \hookrightarrow  \widetilde{\mathcal M}_{a, a', \mathbf{k}, *}|_{U \cap V_{\mathbf{k}}}$ is homology $\lfloor I \rfloor$-connected.
Then the rest of the proof is same as the proof of \cite[Theorem 7.1]{DT24}: we can apply the same discussion starting from Step (3) of the proof of \cite[Theorem 7.1]{DT24}.  
\end{proof}

\subsection{Comparing semi-topological models to the space of positive continuous maps}

The next step is to compare semi-topological models to the space of pointed positive continuous maps.  
We first recall the notion of pointed positive continuous maps from \cite{DT24}.

\begin{defi}[{\cite[Introduction]{DT24}}]
Let $X$ be a smooth projective variety equipped with a finite set of disjoint smooth divisors $\{E_{i} \}$.  Fix a homological curve class $\alpha$.  We define $\Top_{*}^{+}(\mathbb{P}^{1},X)_{\alpha} \subset \Top_{*}(\mathbb{P}^{1},X)_{\alpha}$ to be the subspace consisting of pointed continuous maps $s: \mathbb{P}^{1} \to X$ such that:
\begin{enumerate}
\item $s^{-1}(E_{i})$ is discrete for every $i$, and
\item $s$ has positive local intersection multiplicity against $E_{i}$ at every point of $s^{-1}(E_{i})$.
\end{enumerate}
\end{defi}

Following \cite[Section 8]{DT24} we show that $\mathcal M_{a, a', \mathbf{k}, *}$ introduced in \cref{subsec:algebraicvssemitopological} is weakly homotopy equivalent to the space
$
\mathcal T_{a, a', \mathbf{k}, *},
$
parametrizing pairs $(w, s)$ of $w \in U_{\mathbf{k}}$ and pointed positive continuous maps $s : B \to S$ of the class $(a, a', \mathbf{k})$ such that
\begin{itemize}
\item for any $i = 1, \dotsc, r$, $s^{-1}(E_i)$ is discrete and the local intersection multiplicity at every point of $s^{-1}(E_i)$ is positive, and;
\item at a point of $s^{-1}(E_i)$ the local intersection multiplicity of $s$ along $E_i$ is equal to the multiplicity of $w$ at that point.
\end{itemize}

Let $\mathcal M'_{a, a', \mathbf{k}, *}$ be the space of pairs $(w, s)$ of $w \in U_{\mathbf{k}}$ and pointed continuous maps $s: B \to \mathbb P^1 \times \mathbb P^1$ of class $(a, a')$ such that $s$ intersects $(p_i, p'_i)$ to order exactly $w_i$. Then we have a canonical continuous map
\[
\mathcal M_{a, a', \mathbf{k}, *} \to \mathcal M'_{a, a', \mathbf{k}, *}.
\]
Let $\mathcal T'_{a, a', \mathbf{k}, *} \subset \mathcal T_{a, a', \mathbf{k}, *}$ be the space of pairs $(w, s)$ of $w \in U_{\mathbf{k}}$ and continuous maps $s : B \to S$ such that $s$ intersects $E_i$ to order exactly $w_i$.
Then by taking the strict transform, the argument of \cite[Proposition 8.3]{DT24} shows that $\mathcal M'_{a, a', \mathbf{k}, *}$ is homeomorphic to $\mathcal T'_{a, a', \mathbf{k}, *}$.
In this way we have a continuous map
\[
\mathcal M_{a, a', \mathbf{k}, *} \to \mathcal M'_{a, a', \mathbf{k}, *} \cong \mathcal T'_{a, a', \mathbf{k}, *} \hookrightarrow \mathcal T_{a, a', \mathbf{k}, *}.
\]
Then applying the arguments of \cite[Section 8.2]{DT24}, one can prove that
\[
\mathcal M_{a, a', \mathbf{k}, *} \to \mathcal T_{a, a', \mathbf{k}, *} 
\]
is a weak homotopy equivalence (see \cite[Theorem 8.2]{DT24}).

Finally arguing as in \cite[Proposition 8.11]{DT24} one can show that the projection
\[
\mathcal T_{a, a', \mathbf{k}, *}  \to \mathrm{Top}^+_*(B, S)_{a, a', \mathbf{k}}, \, (w, s) \mapsto s,
\]
is a homeomorphism. Thus we obtain the following theorem:
\begin{theo}[{\cite[Theorem 8.12]{DT24}}]
\label{theo:semitopologicalvstopological}
We have a weak homotopy equivalence
\[
\mathcal M_{a, a', \mathbf{k}, *}  \simeq \mathrm{Top}^+_*(B, S, (a, a', \mathbf{k})).
\]
\end{theo}


\subsubsection{Homological stability and the space of positive continuous maps}
\label{subsubsec:homologicalstability}

We assume that $r = 4$.
Let $a, a', \mathbf{k}$ be non-negative integers such that
\[
2a -\sum_{i = 1}^4k_i > 0, \, 2a' -\sum_{i = 1}^4k_i > 0
\]
We define $I \in \mathbb Q_{\geq 0}$ by
\[
I=   \frac{1}{8} \min\left\{2a -\sum_{i = 1}^4k_i,  2a' -\sum_{i = 1}^4k_i\right\} - \frac{1}{2}.
\]
For $w<x \in (U_{\mathbf{k}, *} < \overline{R})$ we define 
\[
E(w < x) := 4m_0(x)  +3 \sum_{i, j} m_{\ell_{i, j}}(x) + 2\sum_j m_{V_j}(x) +2\sum_i m_{\ell_{i, 1} \oplus \ell_{i, 2}}(x) - 2\sum k_i.
\]
Let $P \subset (U_{\mathbf{k}, *} < \overline{R})$ be the subposet consisting of $w < x$ such that $E(w < x) \leq 2I$.
The following lemma shows that $P$ is proper over $U_{\mathbf{k}}$:
\begin{lemm}
\label{lemm:Eproper}
For any positive constant $\mathfrak T$ the set of $x \in \mathrm{Hilb}(B)^Q$ such that $M(x) := 4m_0(x)  +3 \sum_{i, j} m_{\ell_{i, j}}(x) + 2\sum_j m_{V_j}(x) + 2\sum_i m_{\ell_{i, 1} \oplus \ell_{i, 2}}(x) \leq \mathfrak T$ is closed and downward closed. 
\end{lemm}
\begin{proof}
Downwards closure is immediate because the coefficient of $m_q$ in the linear function in the statement of the lemma increases as $q$ decreases.  
To prove that the set is closed, it suffices to show that if $x'$ is a specialization of $x$ then $M(x') \leq M(x)$.   We will do this via the geometric interpretation of $M(x)$ as a codimension of a kernel of a certain morphism between coherent sheaves.  

We may assume that the total multiplicity of the divisors appearing in $x$ is $N$.   (The same will hold for $x'$.)  Then choose $d \in \mathbb N$ sufficiently so large that 
$$\Gamma(B, \cO(d) \otimes V_1 \oplus \cO(d) \otimes V_2)\to \Gamma(D, \cO(d) \otimes V),$$ is surjective for all divisors $D$ of size $\leq N$. (For example we can take $d = N$.)   In this case we claim that $M(x)$ equals the codimension of $$\Gamma_Q(B,\cV)_x \subseteq \Gamma(B,\cO(d) \otimes V_1 \oplus \cO(d) \otimes V_2).$$   By construction of $\Gamma_Q(B,\cV)_x$ as the fiber product of morphisms between vector bundles, this codimension decreases on specialization, so the same will hold for $M(x)$.    

To establish the claim suppose that $x$ is supported at $p_1, \dotsc, p_r$ and the maximum of the divisors supported at $p_i$ is $\leq n_i$.    By our choice of $d$, we have that $\Gamma(B,\cO(d) \otimes V_1 \oplus \cO(d) \otimes V_2) $ surjects onto $\prod_{i =1}^r \Gamma(n_i p_i, V)$, so it suffices to determine the codimension of the space of incidence conditions imposed by $x$ on $\prod_{i =1}^r \Gamma(n_i p_i, V )$, which splits as a product over $p_i$. The incidence conditions imposed at a single point $p_i$ by $g^x_{p_i}$ only depend on ${\rm sat}(g^x_{p_i})$ which corresponds to a chain $f^x_{p_i}$.  Interpreting $\Gamma(n_i p_i, V)$ as the space of $n_i-1$ jets of $V$, we find that a section $s \in \Gamma(n_i p_i, V)$ satisfies the incidence conditions of $f^x_{p_i}$ if and only if the $l$-th Taylor coefficient of $s$ lies in $K_{f^x_{p_i}(l)}$ for all $l$.    Computing the codimension of this linear space and adding over all $p_i$ yields $M(x)$.  
\end{proof}

Let $M_{a, a', \mathbf{k}, *}$ be the space of pointed rational curves of class $(a, a', \mathbf{k})$ on $S$.
Then we have a natural projection
\[
M_{a, a', \mathbf{k}, *} \to U_{\mathbf{k}, *}.
\]
Let $F_1 \subset U_{\mathbf{k}, *}$ be the complement of the image of the above projection.
It follows from \cref{prop:dimensionestimatepointed} that the real codimension of $F_1$ is greater than $16I$.

Let $(w < x)$ be an element of $P$ and $(w < y)$ be an essential saturated pair such that $x \prec y$.
Let us consider the space
\[
\Gamma(B, \mathcal O(a)\otimes V_1 \oplus \mathcal O(a')\otimes V_2)_y.
\]
We define
\begin{align*}
a_y & = a- m_0(y) - m_{V_2}(y) - \sum_{i = 1}^4m_{\ell_{i, 2}}(y) \\
 a'_y &= a- m_0(y) - m_{V_1}(y) - \sum_{i = 1}^4m_{\ell_{i, 1}}(y) \\
 k_{y, i, j} & = m_{\ell_{i, j}}(y) + m_{\ell_{i, 1} \oplus \ell_{i, 2}}(y) \\
y_j(\ell_{i, j}) & =  \left(\sum_{c \in B} (m_{\ell_{i, j}}(f_{y, c}) + m_{\ell_{i, 1} \oplus \ell_{i, 2}}(f_{y, c}))[c] \right), \quad y_j(0) = \emptyset.
\end{align*}
Then the above space can be identified with
\begin{equation}
\label{equation:P^1incident}
\Gamma(B, \mathcal O(a_y)\otimes V_1)_{y_1}\oplus \Gamma(B, \mathcal O(a'_y)\otimes V_2)_{y_2},
\end{equation}
where $\Gamma(B, \mathcal O(a_y)\otimes V_1)_{y_1}$ is the space of sections $s \in \Gamma(B, \mathcal O(a_y)\otimes V_1)$ such that $s(y_1(\ell_{i, 1})) \subset \ell_{i, 1}\otimes \mathcal O(a_y)$ for any $i = 1, \dotsc, 4$ and we define $\Gamma(B, \mathcal O(a'_y)\otimes V_2)_{y_2}$ in a similar way.
We consider 
\[
E_j(w < y) = 2m_0(y) + 2m_{V_{j'}}(y) +2 \sum_{i = 1}^4m_{\ell_{i, j'}}(y) + \sum_{i = 1}^4 (m_{\ell_{i, j}}(y) + m_{\ell_{i, 1} \oplus \ell_{i, 2}}(y)) - \sum k_i
\]
where $\{j, j'\} = \{1, 2\}$. Since this function is monotone increasing with respect to $y$ we have $E_1(w < y) \geq 0$ and $E_2(w < y)\geq 0$. Moreover we have $E = E_1 + E_2$. Thus we have $E_j \leq 2I + 2$,
\begin{align*}
&2a_y - \sum_{i = 1}^4 k_{y, i, 1}\\ 
&= 2a - \sum k_i -\left(2m_0(y) + 2m_{V_2}(y) +2 \sum_{i = 1}^4 m_{\ell_{i, 2}}(y) + \sum_{i = 1}^4 (m_{\ell_{i, 1}}(y) + m_{\ell_{i, 1} \oplus \ell_{i, 2}}(y)) - \sum k_i\right) \\ 
&\geq 6I+2,
\end{align*}
and 
\begin{align*}
&2a'_y - \sum_{i = 1}^4 k_{y, i, 2}\\ 
&= 2a - \sum k_i -\left(2m_0(y) + 2m_{V_1}(y) +2 \sum_{i = 1}^4 m_{\ell_{i, 1}}(y) + \sum_{i = 1}^4 (m_{\ell_{i, 2}}(y) + m_{\ell_{i, 1} \oplus \ell_{i, 2}}(y)) - \sum k_i\right)\\
&\geq 6I+2.
\end{align*}
The space $\Gamma(B, \mathcal O(a_y)\otimes V_1)_{y_1}$ is unobstructed if $y_1$ is contained in the image of 
\[
Z_{a_y, \mathbf{k}_{y, 1}}^\circ \to U_{\mathbf{k}_{y,1}}.
\]
By combining the equations above with \cref{prop:dimensionestimatepointed} we see that the complement of the image has dimension at most $\sum_{i = 1}^{4} k_{y_1, i, 1} -6I-2$. Thus the dimension of the set $w$ such that $w < y$ where $y$ corresponds to a point $y_{1}$ in the complement of the image of $Z_{a_y, \mathbf{k}_{y, 1}}^\circ \to U_{\mathbf{k}_{y,1}}$ is at most
\[
m_0(y) + m_{V_1}(y) + m_{V_2}(y) +  \sum_{i = 1}^4m_{\ell_{i, 2}}(y) + \sum_{i = 1}^4k_{y, i, 1} -6I-2 \leq \sum k_i - 4I.
\]
Let $F_2$ be the union of closures of such loci while $T$ runs over all combinatorial types we consider. Note that the real codimension of $F_2$ in $U_{\mathbf{k}}$ is greater than or equal to $8I$.
Similarly we construct $F_3$ for the space $\Gamma(B, \mathcal O(a'_y)\otimes V_2)_{y_2}$.
Let $F = F_1 \cup F_2 \cup F_3$. 

Next we analyze $\kappa(T)$ in our situation defined by \cref{eq:kappadef}.  
It follows from our definition that we have
\[
2\gamma(T) = 8m_0(T) + 6 \sum_{i, j} m_{\ell_{i, j}(T)} + 4\sum_j m_{V_j}(T) + 4\sum_{i} m_{\ell_{i, 1} \oplus \ell_{i, 2}}(T), 
\]
and
\[
\mathrm{rank}(T) = 3m_0(T) + 2 \sum_{i, j} m_{\ell_{i, j}}(T) + \sum_j m_{V_j}(T) + \sum_{i} m_{\ell_{i, 1} \oplus \ell_{i, 2}}(T).
\]
Thus we have
\[
\kappa(T) = 5m_0(T) + 4 \sum_{i, j} m_{\ell_{i, j}}(T) + 3\sum_j m_{V_j}(T) + 3\sum_{i} m_{\ell_{i, 1} \oplus \ell_{i, 2}}(T) - 2|\mathrm{Supp}(T)|.
\]
Using this one can prove that when we have $x_1 \prec x_2$, we have $\kappa(x_2) \geq \kappa(x_1) + 1$. This shows that $\kappa(T) \geq \mathrm{rank}(T)$.  
Thus we have $E(T) \leq 2 \, \mathrm{rank}(T) \leq 2\kappa (T)$.  
By applying \cref{theo:algebraicvssemitopological,theo:semitopologicalvstopological} we conclude the following theorem:
\begin{theo}
\label{theo:mainCJS}
Let $B = \mathbb P^1$ and $S$ be a del Pezzo surface of degree $4$.
We fix an integral nef class $\alpha$ on $S$.
Suppose that we have a birational morphism $\beta : S \to \mathbb P^1 \times \mathbb P^1$ satisfying
\[
2a -\sum_{i = 1}^4k_i > 0, \, 2a' -\sum_{i = 1}^4k_i > 0.
\]
We define $I \in \mathbb N$ by
\[
I = \frac{1}{8} \min\left\{2a -\sum_{i = 1}^4k_i,  2a' -\sum_{i = 1}^4k_i\right\} -\frac{1}{2}.
\]
Then $M_{\alpha, *} \hookrightarrow \mathrm{Top}^+_*(B, S)_{\alpha}$ is homology $\lfloor I \rfloor$-connected.
\end{theo}

Finally we prove the Cohen--Jones--Segal conjecture for rational curves on quartic del Pezzo surfaces. 
To compare the topology of the space of pointed positive continuous maps and the space of pointed continuous maps, we will deform $S$ to a smooth toric surface where the desired comparison can be done more easily.  The comparison result for toric surfaces will be proved in \cref{appe:toric}.

\begin{proof}[Proof of \cref{theo:introcjs}] 
Let $B = \mathbb{P}^{1}$.  
Let $\ell(\alpha) = \min\{\lfloor I \rfloor -1, k_{i}(\alpha) \, | \, i = 1, \dotsc, 4\}$.
By \cref{theo:mainCJS} it suffices to show that
\[
\mathrm{Top}^+_*(B, S)_{\alpha} \hookrightarrow \mathrm{Top}_*(B, S)_{\alpha},
\]
induces isomorphisms of homology groups in degree $\leq \min\{k_i(\alpha)\}$.
\Cref{prop:alphainvariance} shows that we may prove the desired property after replacing $\alpha$ by $\alpha + \gamma$ where $\gamma$ is the class of a nef curve on $S$ that is disjoint from each of the four exceptional divisors.  In particular, by repeatedly adding on the nef classes $F_{\alpha}, F'_{\alpha}$ defined in \cref{nota:bignota}, we may assume that
\[
a - k_1 - k_2, a - k_3 - k_4, a'-k_1 - k_3, a' - k_2 - k_4,
\]
are all positive.
Let $Y$ be the toric surface obtained by blowing up the four torus invariant points on $\mathbb P^1 \times \mathbb P^1$.  By considering a continuous deformation of the four blown-up points we can construct a homeomorphism from $Y$ to $S$ that identifies the four exceptional divisors $E_{1,1},E_{1,2},E_{2,1},E_{2,2}$ on $Y$ with the four exceptional divisors $E_{1},E_{2},E_{3},E_{4}$ on $S$.  Moreover this homeomorphism preserves the orientations of these real fourfolds as well as their exceptional divisors.  
By \cref{theo:torichomequiv}, the inclusion
\[
\mathrm{Top}^+_*(B, Y)_{\alpha} \hookrightarrow \mathrm{Top}_*(B, Y)_{\alpha},
\]
induces isomorphisms of homology groups of degree $\leq \min\{k_i(\alpha)\}$.  Since $Y$ and $S$ admit a homeomorphism that is compatible with the exceptional divisors, our assertion follows.
\end{proof}

\section{Manin's Conjecture}
\label{sec:Manin}

Finally, we turn to a proof of \cref{theo:intromaninconj}.  Throughout this section we work over a finite field $k = \mathbb{F}_{q}$.  We refer to \cref{nota:bignota} for our standard notations involving curve classes on our surface.

Suppose that $S$ is a split degree $4$ del Pezzo surface and that $\alpha$ is a nef curve class on $S$.  For convenience of notation we set $B = \mathbb P^1$.  We define the counting function as
\[
N(B, S, -K_S, d) = \sum_{\alpha \in \Nef_{1}(S)_{\mathbb Z},\, -K_S.\alpha \leq d} \# M_{\alpha}(k).
\]
Note that by working with nef curve classes we have removed the exceptional set consisting of $(-1)$-curves.  Given any subcone $\mathcal{C} \subset \Nef_{1}(X)$, we also define
\[
N^{\mathcal{C}}(B, S, -K_S, d) = \sum_{\alpha \in \mathcal{C}_{\mathbb Z},\, -K_S.\alpha \leq d} \# M_{\alpha}(k).
\]

\subsection{Upper bounds}
We first prove a uniform upper bound on $\# M_{\alpha}(k)$ of the expected order (but not necessarily with the expected leading constant).

\begin{lemm}
\label{lemm:upperboundsofpointcountingI}
Let $S$ be a split degree $4$ del Pezzo surface over a finite field $k = \mathbb F_q$ and let $\alpha$ be a nef curve class on $S$.  
Then
\[
\# M_{\alpha}(k) \leq \frac{1}{(1-q^{-1})^{6}}q^{h(\alpha) + 2}.
\]
\end{lemm}

\begin{proof}
First we consider the following fibration:
\[
\psi_\alpha : M_{\alpha} \to Z_{a(\alpha), \mathbf{k}(\alpha)}^{\circ}.
\]
For each $k$-rational point $x$ on $\psi_\alpha(M_{\alpha})$, the fiber $\psi_{\alpha}^{-1}(x)$ is a Zariski open subset of the projective space of dimension $h(\alpha) -2a(\alpha) + 1$ by \cref{lemm:malphaopenimmer}. This implies that we have
\[
\# M_{\alpha}(k) \leq \#Z_{a(\alpha), \mathbf{k}(\alpha)}^{\circ}(k) \frac{q^{h(\alpha) -2a(\alpha) + 2}-1}{q-1}\leq \#Z_{a(\alpha), \mathbf{k}(\alpha)}^{\circ}(k)q^{h(\alpha) -2a(\alpha) +1}\frac{1}{1-q^{-1}}
\]
Next we consider
\[
\phi_{a(\alpha), \mathbf{k}(\alpha)} : Z_{a(\alpha), \mathbf{k}(\alpha)}^{\circ} \to U_{\mathbf{k}(\alpha)}.
\]
For any $k$-rational point $x$ on $\phi_{a(\alpha), \mathbf{k}(\alpha)}(Z_{a(\alpha), \mathbf{k}(\alpha)}^{\circ} )$, the fiber $\phi_{a(\alpha), \mathbf{k}(\alpha)}^{-1}(x)$ is a Zariski open subset of the projective space of dimension $2a(\alpha) + 1 - \sum_{i = 1}^r k_i(\alpha)$ by \cref{lemm:configurationcover_fibration}. Thus we conclude 
\[
\#Z_{a(\alpha), \mathbf{k}(\alpha)}^{\circ}(k) \leq \# U_{\mathbf{k}(\alpha)}(k)q^{2a(\alpha) + 1 - \sum_{i = 1}^r k_i(\alpha)}\frac{1}{1-q^{-1}}.
\]
Finally $U_{\mathbf{k}(\alpha)}$ is a Zariski open subset of $\prod_{i = 1}^{4} \mathrm{Hilb}^{[k_i(\alpha)]}(B)$, hence we conclude
\[
\# U_{\mathbf{k}(\alpha)}(k) \leq q^{\sum_{i = 1}^{4} k_i(\alpha)}\frac{1}{(1-q^{-1})^4}.
\]
Altogether our assertion follows.
\end{proof}

By taking a sum as we vary $\alpha$, we can obtain an upper bound on the counting function.  We first need to recall the definition of the $\alpha$-constant introduced in this context by \cite{Peyre}.  We will use the description given in \cite[5.2 Definition]{Jahnel14} but will use $N_{1}(X)_{\mathbb{Z}}$ instead of the dual lattice $N^{1}(X)_{\mathbb{Z}}^{\vee}$.

\begin{defi}[{\cite[5.2 Definition]{Jahnel14}}] \label{defi:alphaconstant}
Let $k$ be a field and $X$ be a smooth Fano variety defined over $k$. Consider $N_{1}(X)$ with the lattice $\Lambda := N_1(X)_{\mathbb Z}$.  We equip $N_{1}(X)$ with the Lebesgue measure normalized so that the fundamental domain for $\Lambda$ has volume $1$.  For any closed rational polyhedral subcone $\mathcal{C} \subset \Nef_{1}(X)$, we define the $\alpha$-constant as 
\[
\alpha(-K_X,\mathcal{C}) = \rho(X) \vol \{ \alpha \in \mathcal{C} \, | \, -K_X.\alpha \leq 1\}.
\]
This is equal to the top coefficient of the Ehrhart quasi-polynomial of the polytope $\mathcal{C}_{-K_X = 1}$.  When we apply this definition to the entire nef cone, we will simply write $\alpha(-K_{X})$.  
\end{defi}

\begin{theo}
\label{theo:Upperbounds}
Let $S$ be a smooth split degree $4$ del Pezzo surface and let $\mathcal{C} \subset \Nef_{1}(X)$ be a closed rational polyhedral subcone. Then we have
\[
\limsup_{d\to \infty}\frac{N^{\mathcal{C}}(B, S, -K_S, d) }{ q^{d}d^{5}} \leq \frac{\alpha(-K_S,\mathcal{C})q^2}{(1-q^{-1})^{7}}
\]

\end{theo}
\begin{proof}
It follows from \cref{lemm:upperboundsofpointcountingI} 
that  
\[
N^{\mathcal{C}}(B, S, -K_S, d)  \leq \frac{q^2}{(1-q^{-1})^{6}}\sum_{\alpha \in \mathcal{C}_{\mathbb Z},\, -K_S.\alpha \leq d} q^{-p^*K_S.\alpha}.
\]
Let $V \subset N_{1}(X)$ denote the affine hyperplane $V =  \{ \alpha \subset N_{1}(X) | -K_{X} \cdot \alpha = 1 \}$ equipped with the $\mathbb{Z}$-structure coming from $\Lambda \cap V$.  Let $P(d)$ denote the Ehrhart quasi-polynomial for the compact rational polyhedron $\mathcal{C} \cap V$; thus $P(d)$ is a quasi-polynomial of degree $\dim(V) = 5$ and with constant leading coefficient $\alpha(-K_{X},\mathcal{C})$.  We have
\begin{align*}
\frac{q^2}{(1-q^{-1})^{6}}\sum_{\alpha \in \mathcal{C}_{\mathbb Z},\, -K_S.\alpha \leq d} q^{-K_S.\alpha} & = \frac{q^2}{(1-q^{-1})^{6}} \sum_{j = 1}^{d} q^{j} P(j) \\
& \leq \frac{q^{d+2}}{(1-q^{-1})^{7}} \sup_{j = 1,\dotsc,d} \{ P(j) \}.
\end{align*}
Dividing by $q^{d}d^{5}$ and letting $d$ go to infinity we obtain the desired statement.  
\end{proof}

\subsection{Stratifications in the algebraic settings}
Let $(w < x) \in (U_{\mathbf{k}} < \mathrm{Hilb}(B)^Q)(\overline{k})$. We define a linear space $Z_{w < x}^{\mathrm{alg}}$ as
\[
Z_{w < x}^{\mathrm{alg}} := \Gamma(B_{\overline{k}}, \mathcal O(a)\otimes V_1 \oplus \mathcal O(a')\otimes V_2)_x \subset E_w = \Gamma(B_{\overline{k}}, \mathcal O(a)\otimes V_1 \oplus \mathcal O(a')\otimes V_2)_w,
\]
where $\Gamma(B_{\overline{k}}, \mathcal O(a)\otimes V_1 \oplus \mathcal O(a')\otimes V_2)_x$ is the space of $(s, t) \in \Gamma(B_{\overline{k}}, \mathcal O(a)\otimes V_1 \oplus \mathcal O(a')\otimes V_2)$ such that for any $q\in Q$, $(s, t)(x_q) \subset \mathcal K_q$.
As $(w<x)$ vary, these assemble into a stratification
\[
Z^{\mathrm{alg}} = \{ ((w < x), z) \,| \, (w < x) \in (U_{\mathbf{k}} < \mathrm{Hilb}(B)^Q), z \in Z^{\mathrm{alg}}_{w < x} \}  \subset (U_{\mathbf{k}} < \mathrm{Hilb}(B)^Q) \times_{U_{\mathbf{k}}} E,
\]
where $E \to U_{\mathbf{k}}$ is the fibration parametrizing $E_w$ for $w \in U_{\mathbf{k}}$.
Note that this stratification is different from the one defined in \cref{sec:CJS}. However, when $(w < x)$ is essential, we have
\[
Z^{\mathrm{alg}}_{w < x} = Z^{\mathrm{top}}_{w < x} \cap E^{\mathrm{alg}}_{w}.
\]
For this reason, the two stratifications are essentially the same from the point of view of the inclusion-exclusion principle.

\subsection{Counting curves via the bar complex }
\label{subsec:countingcurves}

Here we follow the discussion of \cref{subsubsec:homologicalstability}.  Let $h, a, a', \mathbf{k}$ be non-negative integers such that
\[
h = 2a + 2a' -\sum_i k_i, \, 2a -\sum_{i = 1}^4k_i > 0, \, 2a' -\sum_{i = 1}^4k_i > 0
\]
We define $I \in \mathbb Q_{\geq 0}$ by
\[
I = \frac{1}{8} \min\left\{2a -\sum_{i = 1}^4k_i,  2a'  -\sum_{i = 1}^4k_i\right\} - \frac{1}{2}.
\]
Define
\[
E(w < x) := 4m_0(x)  +3 \sum_{i, j} m_{\ell_{i, j}}(x) + 2\sum_j m_{V_j}(x) +2\sum_i m_{\ell_{i, 1} \oplus \ell_{i, 2}}(x) - 2\sum k_i
\]
and let $P \subset (U_{\mathbf{k}} < \mathrm{Hilb}(B)^Q)$ be the subposet consisting of $(w < x)$ such that $E(w<x) \leq 2I$.
\Cref{lemm:Eproper} shows that $P$ is proper over $U_{\mathbf{k}}$.
Let $M_{a, a', \mathbf{k}}$ be the space of rational curves of class $(a, a', \mathbf{k})$ on $S$.
Then we have a natural projection
\[
M_{a, a', \mathbf{k}} \to U_{\mathbf{k}}
\]
Let $F_1 \subset U_{\mathbf{k}}$ be the complement of the image of the above projection.
It follows from \cref{prop:finaldimestimate} that the codimension of $F_1$ is greater than $8I$.

Let $(w < x)$ be a saturated element of $P(\overline{k})$ and $(w < y)$ be a saturated pair such that $x \prec y$.
Let us consider the space
\[
\Gamma(B_{\overline{k}}, \mathcal O(a)\otimes V_1 \oplus \mathcal O(a')\otimes V_2)_y.
\]
We define
\begin{align*}
a_y &= a- m_0(y) - m_{V_2}(y) - \sum_{i = 1}^4m_{\ell_{i, 2}}(y), \\
a'_y &= a- m_0(y) - m_{V_1}(y) - \sum_{i = 1}^4m_{\ell_{i, 1}}(y), \\
k_{y, i, j} &= m_{\ell_{i, j}}(y) + m_{\ell_{i, 1} \oplus \ell_{i, 2}}(y),
\end{align*}
as before.
We also define
\[
y_j (\ell_{i, j})= \left(\sum_{c \in B(\overline{k})} (m_{\ell_{i, j}}(f_{y, c}) + m_{\ell_{i, 1} \oplus \ell_{i, 2}}(f_{y, c}))[c] \right), \, y_j(0) = \emptyset.
\]
Then the space $\Gamma(B_{\overline{k}}, \mathcal O(a)\otimes V_1 \oplus \mathcal O(a')\otimes V_2)_y$ is identified with
\[
\Gamma(B_{\overline{k}}, \mathcal O(a_y)\otimes V_1)_{y_1}\oplus \Gamma(B_{\overline{k}}, \mathcal O(a'_y)\otimes V_2)_{y_2}, 
\]
where the definition of this space can be found in \eqref{equation:P^1incident}.
We consider 
\[
E_j(w < y) = 2m_0(y) + 2m_{V_{j'}}(y) +2 \sum_{i = 1}^4 m_{\ell_{i, j'}}(y) + \sum_{i = 1}^4 (m_{\ell_{i, j}}(y) + m_{\ell_{i, 1} \oplus \ell_{i, 2}}(y)) - \sum k_i
\] 
where $\{j, j'\} = \{1, 2\}$. Since this function is monotone increasing with respect to $y$ we have $E_1(w < y) \geq 0$ and $E_2(w < y)\geq 0$. Moreover we have $E = E_1 + E_2$. Thus we have $E_j \leq 2I+2$. Following the discussion of \cref{subsubsec:homologicalstability}, we have
\[
2a_y - \sum_{i = 1}^r k_{y, i, 1} \geq 6I + 2, \quad
2a'_y - \sum_{i = 1}^r k_{y, i, 2}  \geq 6I + 2
\]
The space $\Gamma(B_{\overline{k}}, \mathcal O(a_y)\otimes V_1)_{y_1}$ is unobstructed if $y_1$ is contained in the image of 
\[
Z_{a_y, \mathbf{k}_{y, 1}}^\circ \to U_{\mathbf{k}_{y,1}}.
\]
The complement of the image of this map has dimension at most $\sum_{i = 1}^4k_{y, i, 1} -6I-2$. Arguing as in \cref{subsubsec:homologicalstability}, the loci of the corresponding $w$ such that $w < y$ has dimension at most $\sum k_i - 4I.$
Let $F_2$ be the union of the Galois orbits of the closures of all such loci over $\overline{k}$ while $T$ runs over all combinatorial types we consider.  
Note that the codimension of $F_2$ in $U_{\mathbf{k}}$ is greater than or equal to $4I$.
Similarly we construct $F_3$ for the space $\Gamma(B_{\overline{k}}, \mathcal O(a'_y)\otimes V_2)_{y_2}$.
Let $F = F_1 \cup F_2 \cup F_3$. Also note that as in  \cref{subsubsec:homologicalstability}  we have $E(T) \leq 2\kappa (T)$.

Let $\widetilde{M}_{a, a', \mathbf{k}} \to M_{a, a', \mathbf{k}}$ be the $\mathbb G_m^2$-torsor as before.
We are interested in
\[
\#M_{a, a', \mathbf{k}}(k) = (q-1)^{-2}\#\widetilde{M}_{a, a', \mathbf{k}}(k).
\]
By the Grothendieck--Lefschetz trace formula, we have
\[
\#\widetilde{M}_{a, a', \mathbf{k}}(k) = \sum_{i} (-1)^i \mathrm{Tr}(\mathrm{Frob} \curvearrowright H^i_{\text{\'et}, c}((\widetilde{M}_{a, a', \mathbf{k}})_{\overline{k}}, \mathbb Q_\ell)).
\]
The goal of this subsection is to estimate this quantity.  When $i$ is small, we bound the Frobenius trace using the general results of Sawin and Shusterman in \cref{appendix-main}.  When $i$ is large, we will first compare these cohomology groups to the bar complex, then to the cohomological virtual bar complex, at the cost of adding error terms at each step.  Our desired statement (\cref{theo:Maincountingbybar}) arises naturally from the cohomological virtual bar complex via \cref{prop:virtualcount}.

Since a degree $4$ del Pezzo surface $S$ is the intersection of two quadrics in $\mathbb{P}^{4}$, the constant in the proof of \cref{appendix-main} is $C_1 = 2^{32}$.  Then by combining the result of Sawin and Shusterman in \cref{appendix-main}, Deligne's estimates, and Leray's spectral sequence, we have
\begin{equation} \label{eq:mtildesmalli}
\left|\sum_{i< 4a + 4a' -2\sum_j k_j -I + 10} (-1)^i \mathrm{Tr}(\mathrm{Frob} \curvearrowright H^i_{\text{\'et}, c}((\widetilde{M}_{a, a', \mathbf{k}})_{\overline{k}}, \mathbb Q_\ell)) \right| = O(q^{2a + 2a' - \sum k_j -I/2 + 5} C^h_1).
\end{equation}
Let $U = U_{\mathbf{k}} \setminus F$. Then since the codimension of $F$ is greater than or equal to $4I$ and since $\widetilde{M}_{a,a',\mathbf{k}}$ is flat over $U_{\mathbf{k}}$, the inclusion $\widetilde{M}_{a,a',\mathbf{k}}|_{U} \to \widetilde{M}_{a,a',\mathbf{k}}$ induces an isomorphism on \'etale cohomology with compact support in codimension $\leq 8I-2$.  In particular, 
we have
\begin{multline*}
\sum_{i \geq 4a + 4a' -2\sum_j k_j -I + 10} (-1)^i \mathrm{Tr}(\mathrm{Frob} \curvearrowright H^i_{\text{\'et}, c}((\widetilde{M}_{a, a', \mathbf{k}})_{\overline{k}}, \mathbb Q_\ell))\\ =  \sum_{i \geq 4a + 4a' -2\sum_j k_j -I + 10} (-1)^i \mathrm{Tr}(\mathrm{Frob} \curvearrowright H^i_{\text{\'et}, c}((\widetilde{M}_{a, a', \mathbf{k}}|_U)_{\overline{k}}, \mathbb Q_\ell)).  
 \end{multline*} 

 Then assuming $I \geq 1$ we have
\[
E|_U = (\widetilde{M}_{a, a', \mathbf{k}}|_U) \sqcup \mathrm{im}(Z^{\mathrm{alg}}_{P_U} \to E|_U).
\]
This implies that we have
\begin{equation}\label{eq:splittingtildem}
\begin{multlined}[0.9\displaywidth] 
\sum_{i \geq 4a + 4a' -2\sum_j k_j -I + 10} (-1)^i \mathrm{Tr}(\mathrm{Frob} \curvearrowright H^i_{\text{\'et}, c}((\widetilde{M}_{a, a', \mathbf{k}}|_U)_{\overline{k}}, \mathbb Q_\ell))\\
\begin{multlined}=\sum_{i \geq 4a + 4a' -2\sum_j k_j -I + 10} (-1)^i \mathrm{Tr}(\mathrm{Frob} \curvearrowright H^i_{\text{\'et}, c}((E|_U)_{\overline{k}}, \mathbb Q_\ell)) \\
 \phantom{=} - \sum_{i \geq 4a + 4a' -2\sum_j k_j -I + 10} (-1)^i \mathrm{Tr}(\mathrm{Frob} \curvearrowright H^i_{\text{\'et}, c}(\mathrm{im}(Z^{\mathrm{alg}}_{P_U} \to E|_U)_{\overline{k}}, \mathbb Q_\ell)) \\
  + O(q^{2a + 2a' - \sum_j k_j -I/2 +6} C^h_1),
  \end{multlined}
\end{multlined}
\end{equation}
where the last error term is coming from a bound on the Frobenius trace of a subspace of $$H^{4a + 4a' -2\sum_j k_j -\lceil I \rceil + 10}_{\text{\'et}, c}((\widetilde{M}_{a, a', \mathbf{k}}|_U)_{\overline{k}}, \mathbb Q_\ell)$$ obtained by applying \cref{appendix-main}.  
Let us analyze
\[
\sum_{i \geq 4a + 4a' -2\sum_j k_j -I + 10} (-1)^i \mathrm{Tr}(\mathrm{Frob} \curvearrowright H^i_{\text{\'et}, c}(\mathrm{im}(Z^{\mathrm{alg}}_{P_U} \to E|_U)_{\overline{k}}, \mathbb Q_\ell)).
\]
By \cref{theo:5.4inDT} with $R = P$, this quantity is equal to
\[
 \sum_{i \geq 4a + 4a' -2\sum_j k_j -I + 10} (-1)^i \mathrm{Tr}(\mathrm{Frob} \curvearrowright H^i_{\text{\'et}, c}(B(P_U, Z^{\mathrm{alg}}_{P_U})_{\overline{k}}, \mathbb Q_\ell)).
\]
By \cref{prop:ssforbarconstruction} we have a spectral sequence:
\begin{align*}
\bigoplus_{i + j = n} E_1^{i, j} &= \bigoplus_{T: \text{ essential type of $P$}} H^{n}_{\text{\'et}, c}((Z^{\mathrm{alg}}|_{U \times_{U_{\mathbf{k}}} \mathcal N_T})_{\overline{k}}, \mu'(T)[1])\\ 
& \implies H^{n}_{\text{\'et}, c}(B(P_U, Z^{\mathrm{alg}}_{P_U})_{\overline{k}}, \mathbb Z_\ell).
\end{align*}
It follows from \cref{lemm:spectral} that
\begin{equation*}
\begin{multlined}[\displaywidth]
\Big| \sum_{i \geq 4a + 4a' -2\sum_j k_j -I + 10} (-1)^i \mathrm{Tr}(\mathrm{Frob} \curvearrowright H^i_{\text{\'et}, c}(B(P_U, Z^{\mathrm{alg}}_{P_U})_{\overline{k}}, \mathbb Q_\ell))- \\ 
 \sum_{T: \text{ ess. type of $P$}}   \sum_{i \geq 4a + 4a' -2\sum_j k_j -I + 10} (-1)^i \mathrm{Tr}(\mathrm{Frob} \curvearrowright H^{i}_{\text{\'et}, c}((Z^{\mathrm{alg}}|_{U \times_{U_{\mathbf{k}}} \mathcal N_T})_{\overline{k}}, \mu'(T)[1]\otimes\mathbb Q_\ell))\Big|  \leq N
\end{multlined}
\end{equation*}
where $N$ is the sum of the absolute values of eigenvalues of the Frobenius for
\[
\bigoplus_{T: \text{ essential type of $P$}} H^{ \lceil 4a + 4a' -2\sum_j k_j -I +10\rceil}_{\text{\'et}, c}((Z^{\mathrm{alg}}|_{U \times_{U_{\mathbf{k}}} \mathcal N_T})_{\overline{k}}, \mu'(T)[1]\otimes\mathbb Q_\ell).
\]
Our next goal is to bound $N$.  Since $Z^{\mathrm{alg}}$ has the expected codimension over $U \times_{U_{\mathbf{k}}} \mathcal N_T$, for any $i \geq 4a + 4a' -2\sum_j k_j -I+10$ we have
\[
H^{i}_{\text{\'et}, c}((Z^{\mathrm{alg}}|_{U \times_{U_{\mathbf{k}}} \mathcal N_T})_{\overline{k}}, \mu'(T)[1]\otimes \mathbb Q_\ell) = H^i_{\text{\'et}, c}((U \times_{U_{\mathbf{k}}} \mathcal N_T)_{\overline{k}}, \mu(T)(-n_{a, a', \mathbf{k}, T})[-2n_{a, a', \mathbf{k}, T}+1]\otimes \mathbb Q_\ell)
\]
where as before $n_{a, a', \mathbf{k}, T} = 2a + 2a' +4 -2\sum_i k_i -\gamma(T)$.
The cohomological dimension of
\[
H^i_{\text{\'et}, c}((F \times_{U_{\mathbf{k}}} \mathcal N_T)_{\overline{k}}, \mu(T)(-n_{a, a', \mathbf{k}, T})[-2n_{a, a', \mathbf{k}, T}+1]\otimes \mathbb Q_\ell)
\]
is at most
\begin{multline*}
2\sum_j k_j - 8I + 2|\Supp(T)| + \mathrm{rank}(T) + 4a + 4a' + 9 - 4\sum_j k_j - 2\gamma(T)\\ 
\leq 4a + 4a' - 2\sum_j k_j - \kappa(T) - 8I +9
\end{multline*}
This implies that for $i \geq 4a + 4a' -2\sum_j k_j -I+10$, we have
\begin{multline*}
H^i_{\text{\'et}, c}((U \times_{U_{\mathbf{k}}} \mathcal N_T)_{\overline{k}}, \mu(T)(-n_{a, a', \mathbf{k}, T})[-2n_{a, a', \mathbf{k}, T}+1]\otimes \mathbb Q_\ell)\\ 
\cong H^i_{\text{\'et}, c}((\mathcal N_T)_{\overline{k}}, \mu(T)(-n_{a, a', \mathbf{k}, T})[-2n_{a, a', \mathbf{k}, T}+1] \otimes \mathbb Q_\ell).
\end{multline*}
Thus the following proposition follows from \cref{theo:errorforvirtual} applied with $\eta = 1/2$ combined with the inequality $h \geq \max \{ \sum_{i=1}^{4} k_{i}, I \}$: 
\begin{prop}
For any constants $C_2 > 160$ and $C_3 > (160)^{4}$ if $q > C_3$ we have
\[
N = O( q^{2a + 2a' - \sum_j k_j -I/4 +4} C^h_2).
\]
\end{prop}
Thus we conclude that 
\begin{multline*}
\sum_{i \geq 4a + 4a' -2\sum_j k_j -I + 10} (-1)^i \mathrm{Tr}(\mathrm{Frob} \curvearrowright H^i_{\text{\'et}, c}(B(P_U, Z^{\mathrm{alg}}_{P_U})_{\overline{k}}, \mathbb Q_\ell)) \\ 
= \sum_{T: \text{ ess. type of $P$}}   \sum_{i \geq 4a + 4a' -2\sum_j k_j -I + 10} (-1)^i \mathrm{Tr}(\mathrm{Frob} \curvearrowright H^{i}_{\text{\'et}, c}(((\mathcal N_T)_{\overline{k}}, \mu(T)(-n_{a, a', \mathbf{k}, T})[-2n_{a, a', \mathbf{k}, T}+1])\otimes \mathbb Q_\ell)\\
+ O(q^{2a + 2a' - \sum_j k_j -I/4+6} \max\{C_1, C_2\}^h).
\end{multline*}
When $T$ is an essential type which is not a type of $P$, we have $\kappa (T) > I$. This implies that the cohomological dimension of 
\[
H^{i}_{\text{\'et}, c}(((\mathcal N_T)_{\overline{k}}, \mu(T)(-n_{a, a', \mathbf{k}, T})[-2n_{a, a', \mathbf{k}, T}+1]\otimes \mathbb Q_\ell)
\]
is at most
\[
4a + 4a' + 9 - 2\sum_j k_j -2\gamma(T) + \mathrm{rank}(T) + 2|\Supp(T)| < 4a + 4a' + 9- 2\sum_j k_j -I.
\]
We conclude that 
\begin{align*}
&\sum_{T: \text{ ess. type of $P$}}   \sum_{i \geq 4a + 4a' -2\sum_j k_j -I + 10} (-1)^i \mathrm{Tr}(\mathrm{Frob} \curvearrowright H^{i}_{\text{\'et}, c}(((\mathcal N_T)_{\overline{k}}, \mu(T)(-n_{a, a', \mathbf{k}, T})[-2n_{a, a', \mathbf{k}, T}+1]\otimes \mathbb Q_\ell))\\
&= \sum_{T: \text{ ess. type}}   \sum_{i \geq 4a + 4a' -2\sum_j k_j -I + 10} (-1)^i \mathrm{Tr}(\mathrm{Frob} \curvearrowright H^{i}_{\text{\'et}, c}(((\mathcal N_T)_{\overline{k}}, \mu(T)(-n_{a, a', \mathbf{k}, T})[-2n_{a, a', \mathbf{k}, T}+1]\otimes \mathbb Q_\ell)).
\end{align*}
Then using \cref{theo:errorforvirtual} again, we obtain
\begin{multline*}
\sum_{i \geq 4a + 4a' -2\sum_j k_j -I + 10} (-1)^i \mathrm{Tr}(\mathrm{Frob} \curvearrowright H^i_{\text{\'et}, c}(B(P_U, Z^{\mathrm{alg}}_{P_U})_{\overline{k}}, \mathbb Q_\ell))\\
= \sum_{T: \text{ essential type}} \sum_i \mathrm{Tr}(\mathrm{Frob} \curvearrowright H^i_{\text{\'et}, c}((\mathcal N_T)_{\overline{k}}, \mu(T)(-n_{a, a', \mathbf{k}, T})[-2n_{a, a', \mathbf{k}, T}+1]\otimes \mathbb Q_\ell))\\
+ O(q^{2a + 2a' - \sum_j k_j -I/4 +6} \max\{C_1, C_2\}^h).
\end{multline*}
Now it follows from \cref{prop:virtualcount} that 
\begin{multline*}
\sum_{T: \text{ essential type}} \sum_i \mathrm{Tr}(\mathrm{Frob} \curvearrowright H^i_{\text{\'et}, c}((\mathcal N_T)_{\overline{k}}, \mu(T)(-n_{a, a', \mathbf{k}, T})[-2n_{a, a', \mathbf{k}, T}+1]\otimes \mathbb Q_\ell))\\
\begin{split}
&= \sum_{T: \text{ essential type}}-q^{2a+ 2a'+4} \sum_{(w< x) \in \mathcal N_T(k)} \mu_k(w, x)q^{-\gamma(x)} \\
&= -q^{2a+ 2a'+4} \sum_{(w< x) \in (U_{\mathbf{k}} < Q^{JB})(k)} \mu_k(w, x)q^{-\gamma(x)},
\end{split}
\end{multline*}
We conclude that
\begin{multline*}
\sum_{i \geq 4a + 4a' -2\sum_j k_j -I + 10} (-1)^i \mathrm{Tr}(\mathrm{Frob} \curvearrowright H^i_{\text{\'et}, c}(B(P_U, Z^{\mathrm{alg}}_{P_U})_{\overline{k}}, \mathbb Q_\ell))\\
\qquad \qquad = -q^{2a+ 2a'+4} \sum_{(w< x) \in (U_{\mathbf{k}} < Q^{JB})(k)} \mu_k(w, x)q^{-\gamma(x)} + O(q^{2a + 2a' - \sum_j k_j -I/4 +6} \max\{C_1, C_2\}^h).
\end{multline*}
Similarly, one can prove that
\begin{multline*}
\sum_{i \geq 4a + 4a' -2\sum_j k_j -I + 10} (-1)^i \mathrm{Tr}(\mathrm{Frob} \curvearrowright H^i_{\text{\'et}, c}((E|_U)_{\overline{k}}, \mathbb Q_\ell)) \\ 
= q^{2a+ 2a'+4}\sum_{w \in U_{\mathbf{k}}(k)} q^{-\gamma(w)} + O(q^{2a + 2a' - \sum_j k_j -I/4 +6} \max\{C_1, C_2\}^h).
\end{multline*}
Combining the previous two equations with \cref{eq:mtildesmalli,eq:splittingtildem}, we obtain:
\begin{theo}
\label{theo:Maincountingbybar}
Assume that $q > C_3$.
Then we have
\[
\#\widetilde{M}_\alpha(k) = q^{2a + 2a'+4} \sum_{(w\leq x) \in (U_{\mathbf{k}} \leq Q^{JB})(k)} \mu_k(w, x)q^{-\gamma(x)} + O(q^{2a + 2a' - \sum_j k_j -I/4 +6} \max\{C_1, C_2\}^h),
\]
where the implied constant does not depend on $q, a, a', \mathbf{k}$.

\end{theo}

\subsection{The residue of the virtual bar complex}
\label{subsec:virtualpointcounts}

We work with the set up described in \cref{sec:barconstruction}.
We assume that our ground field $k$ is a finite field $\mathbb F_q$ and that $B = \mathbb{P}^{1}_{k}$.

Let $\alpha$ be an ample curve class on a split quartic del Pezzo surface $S$ such that $\alpha$ lies in the interior of the region $\mathcal T_{\alpha}$ defined in \cref{nota:bignota}.  Using the invariants $a(\alpha), a'(\alpha), k_i(\alpha)$ of \cref{nota:bignota}, we have:
\[
2a(\alpha) - \sum_{i = 1}^4 k_i(\alpha) > 0, 2a'(\alpha) - \sum_{i = 1}^4 k_i(\alpha) > 0.
\]

By \cref{theo:Maincountingbybar}, we would like to understand the limiting behavior of
\[
q^{2a + 2a' +4} \sum_{(w\leq x) \in (U_{\mathbf{k}} \leq Q^{JB})(k)} \mu_k(w, x)q^{-\gamma(x)}
\]
as $a, a', k_i \to \infty$. To this end, we consider the following virtual height zeta function: 
\begin{equation*}
\mathsf Z(\mathbf{t}) = \sum_{k_1, \dotsc, k_4 = 0}^\infty q^{\sum_{i = 1}^4k_i}\left( \sum_{(w\leq x) \in (U_{\mathbf{k}} \leq Q^{JB})(k)} \mu_k(w, x)q^{-\gamma(x)} \right)t_1^{k_1}t_2^{k_2}t_3^{k_3}t_4^{k_4}.
\end{equation*}
It follows from \cref{lemm:Mobius}.(1) that the above zeta function can be rewritten as the following convergent Euler product: 
\begin{align*}
&\mathsf Z(\mathbf{t}) \\ &= \prod_{c \in |B|}\left(\sum_{x \in \mathrm{ch}(Q(k))} \mu_k(\hat{1}, x)q^{-\gamma (x)|c|} + \sum_{i = 1}^4\sum_{d_i = 1}^\infty (qt_i)^{|c|d_i}\sum_{x \in \mathrm{ch}(Q(k)), \,  d_i[\ell_{i, 1} \oplus \ell_{i, 2}]\leq x} \mu_k(d_i[\ell_{i, 1} \oplus \ell_{i, 2}], x)q^{-\gamma (x)|c|} \right)
\end{align*}
where $|B|$ is the set of closed points on $B$ and $|c|$ is the degree $[k(c):k]$. Then it follows from \cref{lemm:Mobius}.(2) and \cref{exam:essentialsaturatedexample} that the above Euler product is equal to
\begin{equation*}
\prod_{c \in |B|}\left( 1 - 6q^{-2|c|} + 8q^{-3|c|} -3q^{-4|c|}+ \sum_{i = 1}^4\sum_{d_i = 1}^\infty (qt_i)^{|c|d_i} (q^{-2d|c|}  - 2q^{-(2d_i+1)|c|} + 2q^{-(2d_i + 3)|c|} - q^{-(2d_i + 4)|c|})\right).
\end{equation*}
It is clear from the formula that this Euler product absolutely converges when there is a $\delta > 0$ such that $|t_i| \leq q^{-\delta}$ for $i = 1, \dotsc, 4$. The following proposition verifies the basic properties of $\mathsf Z(\mathbf{t})$.

\begin{prop}
\label{prop:formalheightzeta}

We have the following Euler product:
\begin{multline*}
\prod_{i = 1}^4Z_{\mathbb P^1}(q^{-1}t_i)^{-1}\mathsf Z(\mathbf{t}) \\
 = \prod_{c \in |B|}\left(1 +O(q^{-2|c|})(1 + \sum_i t_i^{|c|})+  \sum_{\substack{\ell_1, \ell_2, \ell_3, \ell_4 \in \{0,1\} \\ \ell_1+\ell_2 + \ell_3 + \ell_4 > 1}}O(q^{-( \ell_1 + \ell_2 + \ell_3 + \ell_4)|c|})t_1^{|c|\ell_1}t_2^{|c|\ell_2}t_3^{|c|\ell_3}t_4^{|c|\ell_4} \right),
\end{multline*}
where $Z_{\mathbb P^1}(t)$ is the Hasse--Weil zeta function for $\mathbb P^1$.
In particular if $0 < \delta < 1/2$ then this Euler product absolutely converges when $|t_i|\leq q^{\delta}$ for $i = 1, \dotsc, 4$.
Furthermore we have
\[
\lim_{t_i \to 1} \prod_{i = 1}^4 (1-t_i) \cdot \mathsf Z(\mathbf{t}) = (1-q^{-1})^{-4}\prod_{c \in |B|}(1-q^{-|c|})^6(1+6q^{-|c|} + q^{-2|c|}).
\]
\end{prop}

\begin{proof}

First we have
\begin{align*}
&\begin{multlined}\prod_{i = 1}^4 Z_{\mathbb P^1}(q^{-1}t_i)^{-1}
\\ \qquad \times\prod_{c \in |B|}\left( 1 - 6q^{-2|c|} + 8q^{-3|c|} -3q^{-4|c|}+ \sum_{i = 1}^4\sum_{d_i = 1}^\infty t_i^{|c|d_i}q^{-d_i|c|} (1  - 2q^{-|c|} + 2q^{-3|c|} - q^{-4|c|}) \right)\end{multlined}\\
&\begin{multlined}=\prod_{c \in |B|}\prod_{i = 1}^4(1-t_i^{|c|}q^{-|c|}) \\ 
\qquad \times \left( 1 - 6q^{-2|c|} + 8q^{-3|c|}-3q^{-4|c|}+ \sum_{i = 1}^4\frac{t_i^{|c|}q^{-|c|}}{1-t_i^{|c|}q^{-|c|}} (1  - 2q^{-|c|} + 2q^{-3|c|} - q^{-4|c|})\right)\end{multlined}\\
&\begin{multlined}=\prod_{c \in |B|}\Bigg( \left(\prod_{i = 1}^4(1-t_i^{|c|}q^{-|c|}) \right)(1 - 6q^{-2|c|} + 8q^{-3|c|}-3q^{-4|c|})\\
\qquad + \sum_{i = 1}^4\left(\prod_{j \neq i}(1-t_j^{|c|}q^{-|c|})\right)t_i^{|c|}q^{-|c|}(1  - 2q^{-|c|} + 2q^{-3|c|} - q^{-4|c|})\Bigg).\end{multlined}
\end{align*}
The first assertion follows from this expression.
Also one can compute the leading term as
\begin{multline*}
\lim_{t_i \to 1} \prod_{i = 1}^4 (1-t_i)Z_{\mathbb P^1}(q^{-1}t_i)\\ 
\times \prod_{c \in |B|}(1-q^{-|c|})^4
 \left( 1 - 6q^{-2|c|} + 8q^{-3|c|} -3q^{-4|c|}+ \frac{4q^{-|c|}}{1-q^{-|c|}} (1  - 2q^{-|c|} + 2q^{-3|c|} - q^{-4|c|})\right)\\
= (1-q^{-1})^{-4}\prod_{c \in |B|}(1-q^{-|c|})^6(1+6q^{-|c|} + q^{-2|c|}). \qedhere
\end{multline*}
\end{proof}

Next we expand the product: 
\begin{equation*}
\prod_{i = 1}^4 (1-t_i)\mathsf Z(\mathbf{t}) = \sum_{j_1, j_2, j_3, j_4 = 0}^\infty B_{j_1j_2j_3j_4}t_1^{j_1}t_2^{j_2}t_3^{j_3}t_4^{j_4}.
\end{equation*}
We fix $\delta$ satisfying $0 < \delta < 1/2$.  Assume that $|t_i|\leq q^\delta$ for $i = 1, 2, 3, 4$.
Then \cref{prop:formalheightzeta} shows that the above series absolutely converges.
This implies that $|B_{j_1j_2j_3j_4}| = O(q^{-\eta_1(j_1 + j_2 + j_3 + j_4)})$ for any $0 < \eta_1 < \delta$.
This leads to the following proposition:

\begin{prop}
\label{prop:expectedbar}
There exists $\eta_1 > 0$ which does not depend on $\mathbf{k}, q$ such that  
\[
q^{\sum_i k_i}\sum_{(w\leq x) \in (U_{\mathbf{k}} \leq Q^{JB})(k)} \mu_k(w, x)q^{-\gamma(x)} = (1-q^{-1})^{-4}\prod_{c \in |B|}(1-q^{-|c|})^6(1+6q^{-|c|} + q^{-2|c|}) + O(q^{-\eta_1 \min\{k_i\}}),
\]
where the implied $O$-constant does not depend on $\mathbf{k}, q$. 
\end{prop}
\begin{proof}
Recall that $\mathsf{Z}(\mathbf{t})$ can be obtained from the power series $\sum B_{j_{1}j_{2}j_{3}j_{4}} t_{1}^{j_{1}}t_{2}^{j_{2}}t_{3}^{j_{3}}t_{4}^{j_{4}}$ by multiplying by $\prod_{i=1}^{4}(1-t_{i})^{-1}$.  Comparing coefficients, we see that
\[
q^{\sum_i k_i}\sum_{(w\leq x) \in (U_{\mathbf{k}} \leq Q^{JB})(k)} \mu_k(w, x)q^{-\gamma(x)} = \sum_{j_i \leq k_i} B_{j_1j_2j_3j_4.}
\]
\Cref{prop:formalheightzeta} shows that 
\begin{multline*}
\lim_{k_i \to \infty}q^{\sum_i k_i}\sum_{(w\leq x) \in (U_{\mathbf{k}} \leq Q^{JB})(k)} \mu_k(w, x)q^{-\gamma(x)} 
\\ = \sum_{j_1, j_2, j_3, j_4 = 0}^\infty B_{j_1j_2j_3j_4} = (1-q^{-1})^4\prod_{c \in |B|}(1-q^{-|c|})^6(1+6q^{-|c|} + q^{-2|c|})
\end{multline*}
Thus the error term in the proposition statement is given by 
\[
\sum_{j_1 > k_1}|B_{j_1j_2j_3j_4}| + \sum_{j_2 > k_2} |B_{j_1j_2j_3j_4}|+ \sum_{j_3 > k_3}|B_{j_1j_2j_3j_4}| + \sum_{j_4 > k_4} |B_{j_1j_2j_3j_4}| = \sum_{i = 1}^4O(q^{-\eta_1 k_i})
\]
for some $\eta_1 > 0$ and our assertion follows.
\end{proof}

\subsection{Peyre's constant}

Finally, to get an explicit formula for the counting function we must compute Peyre's constant explicitly.  We follow the presentation of \cite{Peyre} and \cite{BT}.
Let $S$ be a split degree $4$ del Pezzo surface over $\mathbb{F}_{q}$ with the anticanonical height.  \cite{Peyre} and \cite{BT} define
\begin{equation*}
c(\mathbb P^1, S, -K_S) = (1-q^{-1})^{-1}\alpha(-K_{S}) \beta(S) \tau_{-K_{S}}(S)
\end{equation*}
where $\alpha(-K_{S})$ is defined as in \cref{defi:alphaconstant}, $\beta(S) = \# H^{1}(\mathbb{F}_{q}(t),\Pic(S_{\overline{\mathbb{F}_{q}(t)}}))$, and $\tau_{-K_{S}}(S)$ is the Tamagawa constant.  The value $\alpha(-K_{S}) = \frac{1}{180}$ is computed by \cite[Theorem 4]{Derenthal07}.  Due to the splitness assumption $\beta(S) = 1$.  Since we are working with a trivial family of split surfaces, the constant $\tau_{-K_{S}}(S)$ has the following simple form:
\begin{equation*}
\tau_{-K_{S}}(S) = q^{2} \left(1 - q^{-1} \right)^{-6} \left( \prod_{c \in |B|} \left( 1- q^{-|c|} \right)^6  \frac{\# S(\mathbb F_{q^{|c|}})}{q^{2|c|}} \right) 
\end{equation*}
as $c \in |B|$ varies over all places where $|c|$ denotes $[k(c):k]$. 
Using $\# S(\mathbb{F}_{q^{|c|}}) = q^{2|c|} + 6 q^{|c|} + 1$, we obtain 
\begin{equation*}
c(\mathbb P^1, S, -K_S) = \frac{q^{2}}{180} \left(1 - q^{-1} \right)^{-7} \left( \prod_{c \in |B|} \left( 1- q^{-|c|} \right)^6  \left( 1 + 6q^{-|c|} + q^{-2|c|} \right)  \right)
\end{equation*}

\subsection{The proofs of the main theorems}

In this section we prove \cref{theo:intromaninconj}. Let $S$ be a split quartic del Pezzo surface over $k = \mathbb F_q$ as before.



Let $\mathfrak P$ be the set of birational morphisms $\rho : S \to \mathbb P^1 \times \mathbb P^1$.
For such a birational morphism $\rho : S \to \mathbb P^1 \times \mathbb P^1$, we define $F, F', E_1, \dotsc, E_4$ as \cref{nota:bignota}.
For each class $\alpha \in \mathrm{Nef}_1(S)$, we define
\[
\mathfrak I_\rho(\alpha) = \frac{1}{32}\min \{ 2F.\alpha - \sum_{i = 1}^4 E_i.\alpha, 2F'.\alpha - \sum_{i = 1}^4 E_i.\alpha\},
\]
and 
\[
\mathfrak R_\rho(\alpha) = \min_i \{ E_i.\alpha\}.
\]
Then we define a rational continous homogeneous piecewise linear function $\ell$ on $\mathrm{Nef}_1(S)$ by
\[
\ell(\alpha) = \max_{\rho \in \mathfrak P} \{\min \{\mathfrak I_\rho(\alpha), \mathfrak R_\rho(\alpha)\}\}.
\]
Then it follows from \cref{prop:nefexpression} that $U = \{ \alpha \in \mathrm{Nef}_1(S) \, | \, \ell(\alpha) > 0\}$ is a dense open cone in $\mathrm{Nef}_1(S)$.

\begin{proof}[Proof of \cref{theo:intromaninconj}:] 
Our goal is to understand the asymptotic formula for
\[
N^{\Nef_1(S)_{\ell, \epsilon}}(B, S, -K_S, d) = \sum_{\alpha \in \Nef_1(S)_{\ell, \epsilon, \mathbb Z},\, -K_S.\alpha \leq d} \# M_{\alpha}(k),
\]
where $\epsilon > 0$ is a small rational number and $\Nef_1(S)_{\ell, \epsilon}$ is the rational polyhedral cone consisting of $\alpha \in \Nef_1(S)$ satisfying $\ell(\alpha) \geq -\epsilon K_S.\alpha$.
We write
\begin{align*}
\frac{N^{\Nef_1(S)_{\ell, \epsilon, \mathbb Z}}(B, S, -K_S, d)}{q^dd^5} &=  (q^dd^5)^{-1}\sum_{\alpha \in \Nef_1(S)_{\ell, \epsilon}, -K_S.\alpha \leq d} \# M_{\alpha}(k)\\
&=  (q^dd^5)^{-1}\sum_{\alpha \in \Nef_1(S)_{\ell, \epsilon, \mathbb Z}, -K_S.\alpha \leq d} q^{-2}(1-q^{-1})^{-2}\# \widetilde{M}_{\alpha}(k)
\end{align*}
We next estimate the sum.  By construction, for any $\alpha \in \Nef_1(S)_{\ell, \epsilon, \mathbb Z}$ we have $$\frac{1}{32}\min\{ 2a(\alpha) - \sum_ik_i( \alpha), 2a'(\alpha) - \sum_ik_i( \alpha), k_i(\alpha)\} \geq \epsilon h(\alpha) = - \epsilon K_S.\alpha,$$
where we choose a birational morphism $\rho : S \to \mathbb P^1 \times \mathbb P^1$ according to \cref{nota:bignota}.
Assume that we have $q^\epsilon > C = \max\{C_1, C_2\}$ and $q > C_3$ where $C_1, C_2, C_3$ are constants defined in \cref{subsec:countingcurves}. (One can take $C = 2^{32}$ and $C_3 > (160)^4$.)
It follows from \cref{theo:Maincountingbybar} that
\begin{equation*}
\begin{multlined}[\displaywidth]
\sum_{\alpha \in \Nef_1(S)_{\ell, \epsilon, \mathbb Z}, -K_S.\alpha \leq d} q^{-2}(1-q^{-1})^{-2}\# \widetilde{M}_{\alpha}(k) = (1-q^{-1})^{-2}\times \\ 
\left( \sum_{\alpha \in \Nef_1(S)_{\ell, \epsilon, \mathbb Z}, -K_S.\alpha \leq d} \left(q^{h(\alpha) + 2+ \sum_i k_i(\alpha)}\left( \sum_{(w\leq x) \in (U_{\mathbf{k}(\alpha)} \leq Q^{JB})(k)} \mu_k(w, x)q^{-\gamma(x)}\right) + O(q^{(1-\delta\epsilon) h(\alpha)}) \right)\right),
\end{multlined}
\end{equation*}
for some $0 < \delta < 1$.
Thus  using \cref{prop:expectedbar} we obtain 
\begin{multline*}
 \sum_{\alpha \in \Nef_1(S)_{\ell, \epsilon, \mathbb Z}, -K_S.\alpha \leq d} q^{-2}(1-q^{-1})^{-2}\# \widetilde{M}_{\alpha}(k)
 \\=  \sum_{\alpha \in \Nef_1(S)_{\ell, \epsilon, \mathbb Z}, -K_S.\alpha \leq d} \left(\tau_{-K_S}(S)q^{h(\alpha)}   + O(q^{(1-\delta'\epsilon )h(\alpha)})\right),
\end{multline*}
for some $0 < \delta' < \min\{\delta, \eta_1\delta\}$ where $\eta_1$ is a constant coming from \cref{prop:expectedbar}.
Thus by using arguemtns using the Ehrhart quasi-polynomials as Theorem~\ref{theo:Upperbounds}, we conclude
\begin{equation*}
\lim_{d\to+ \infty}\frac{N^{\Nef_1(S)_{\ell, \epsilon}}(B, S, -K_S, d)}{q^dd^5} =  (1-q^{-1})^{-1}\alpha(-K_S, \Nef_1(S)_{\ell, \epsilon})\tau_{-K_S}(S). \qedhere
\end{equation*}
\end{proof}

\begin{proof}[Proof of \cref{coro:weakmaninconj}:]
Fix any prime power $q$.  \cref{theo:Upperbounds} shows that for $d$ sufficiently large we have
\begin{equation*}
0 \leq N(\mathbb{P}^{1},S,-K_{S},d) \leq 2 \frac{\alpha(-K_S, \Nef_{1}(X))}{(1-q^{-1})^{7}} q^{d+2}d^{5}
\end{equation*}
In particular, for any finite set of prime powers $q$ we can ensure that \cref{coro:weakmaninconj} holds by choosing $R$ large enough.  Thus it suffices to prove \cref{coro:weakmaninconj} for sufficiently large prime powers $q$.

Fix a constant $C$ as in \cref{theo:intromaninconj}.  Fix a sufficiently large prime power $q$ and let $\epsilon$ be a rational number satisfying $\log(C)/\log(q) < \epsilon < 2\log(C)/\log(q)$.  We define the cone $\Nef_1(S)_{\ell, \epsilon}$ as before and set
\begin{equation*}
\mathcal{C} := \overline{\Nef_{1}(S) \backslash \Nef_1(S)_{\ell, \epsilon}}.
\end{equation*}
\Cref{theo:Upperbounds} shows that for $d$ sufficiently large we have
\begin{equation*}
N^{\Nef_1(S)_{\ell, \epsilon}}(\mathbb{P}^{1},S,-K_{S},d) \leq N(\mathbb{P}^{1},S,-K_{S},d) \leq N^{\Nef_1(S)_{\ell, \epsilon}}(\mathbb{P}^{1},S,-K_{S},d) + 2\frac{\alpha(-K_S,\mathcal{C})q^2}{(1-q^{-1})^{7}} q^{d}d^{5}
\end{equation*}
and thus by \cref{theo:intromaninconj}
\begin{multline*}
\limsup_{d \to \infty}\left| \frac{N(\mathbb P^1, S, -K_S, d)}{q^{d+2} d^{5}}  - (1-q^{-1})^{-1}q^{-2}\alpha(-K_S, \Nef_{1}(S))\tau_{-K_S}(S) \right| \\
\begin{aligned}& \leq \left| \alpha(-K_{S},\Nef_1(S)_{\ell, \epsilon}) - \alpha(-K_S, \Nef_{1}(S) ) \right| \cdot (1-q^{-1})^{-1} q^{-2} \tau_{-K_S}(S) + 2\frac{\alpha(-K_S,\mathcal{C})}{(1-q^{-1})^{7}} \\
& = \alpha(-K_S,\mathcal{C}) \frac{q^{7}}{(q-1)^{7}} \left( 2 + \left( \prod_{c \in |B|} \left( 1- q^{-|c|} \right)^6  \frac{\# S(\mathbb F_{q^{|c|}})}{q^{2|c|}} \right) \right) \\
& \leq \alpha(-K_S,\mathcal{C}) 2^{7} \left( 2 + \left( \prod_{c \in |B|} \left( 1- 2^{-|c|} \right)^6  \frac{\# S(\mathbb F_{2^{|c|}})}{2^{2|c|}} \right) \right)
\end{aligned}
\end{multline*}
Note that when $\epsilon$ is sufficiently small then $\alpha(-K_{S},\mathcal{C})$ is bounded above by a linear function of $\epsilon$ (determined by the shape of $\Nef_{1}(S)$ independently of $q$).  Thus when $q$ is sufficiently large $\alpha(-K_{S},\mathcal{C})$ is bounded by a linear function in $1/\log(q)$.  Altogether, we see that there is an $R$ (not depending on $q$) yielding the desired inequality.
\end{proof}

\appendix

\section{Results for degree \texorpdfstring{$4$}{4} toric surfaces} \label{appe:toric}

Suppose that our ground field is $\mathbb C$ and let $Y$ denote the toric surface obtained by blowing up $\mathbb{P}^{1} \times \mathbb{P}^{1}$ at the four torus-invariant points denoted by $(p_i, p_j)$ for $i, j = 1, 2$. We denote the exceptional divisor above $(p_{i},p_{j})$ by $E_{i, j}$.  In this section we analyze the moduli spaces of rational curves on $Y$.  Note that spaces of rational curves on toric varieties have been well-studied both from the perspective of homological stability (\cite{Guest}) and from the perspective of point counting (\cite{BTtoric98, Bourqui03, Bourqui11}).  However, we need a specific result -- a homology equivalence between the spaces of positive pointed maps and the space of continuous pointed maps -- that does not directly follow from the literature.

\begin{theo} \label{theo:torichomequiv}
For any integral nef curve class $\alpha \in \Nef_1(Y)$, the inclusion
\begin{equation*}
\Top^{+}_{*}(\mathbb{P}^{1},Y)_{\alpha} \to \Top_{*}(\mathbb{P}^{1},Y)_{\alpha}
\end{equation*}
induces isomorphisms of homology groups of degree $\leq \min_{i, j}\{E_{i, j}.\alpha\}$.
\end{theo}

Our plan is as follows.  For any $\alpha \in \Nef_{1}(Y)$, let $M_{\alpha,*}$ denote the moduli space of pointed algebraic morphisms of class $\alpha$.  By mimicking the material in \cref{sec:CJS}, we show that the inclusion $M_{\alpha,*} \to \Top^{+}_{*}(\mathbb{P}^{1},Y)_{\alpha}$ induces an equality of homology groups in some range.  In turn, by \cite{Guest} we also know that the inclusion $M_{\alpha,*} \to \Top_{*}(\mathbb{P}^{1},Y)_{\alpha}$ induces an equality of homology groups in some range.  By combining these two statements we obtain \cref{theo:torichomequiv}.

\subsection{Set-up}

We identify $\mathbb{P}^{1} \times \mathbb{P}^{1} \cong \mathbb{P}(V_{1}) \times \mathbb{P}(V_{2})$ for rank $2$ vector spaces $V_{1},V_{2}$.  For $j=1,2$ we let $\ell_{1,j}, \ell_{2,j}$ denote respectively the lines corresponding to the points $0, \infty$ in $\mathbb{P}(V_{j})$.  Letting $V = V_{1} \oplus V_{2}$, we define the $k$-poscheme
\begin{equation*}
Q' = \{ V, V_1 \oplus \{0\}, \{0\} \oplus V_2\}\cup\{\ell_{i,1} \oplus \ell_{i,2}, \ell_{i,1} \oplus \ell_{3-i,2}, \ell_{i, 1} \oplus \{0\}, \{0\} \oplus \ell_{i, 2}, 0\}_{i \in \{1,2\}}.
\end{equation*}

\subsection{Rational curves on the degree $4$ toric surface over $\mathbb P^1 \times \mathbb P^1$}

Recall that we realize our degree $4$ toric surface $Y$ as a blow-up
\[
\rho : Y \to \mathbb P^1 \times \mathbb P^1.
\]
Let $E_{i, j}$ be the exceptional divisor of $\rho$ at $(p_i, p_j')$ for  $(i, j = 1, 2)$ and $F, F'$ be general fibers of the two conic fibrations coming from $\mathbb P^1 \times \mathbb P^1$. Also recall that for any nef class $\alpha \in \Nef_1(Y)$, we define
\[
a(\alpha) = F.\alpha, \quad a'(\alpha) = F'. \alpha, \quad k_{i, j}(\alpha) = E_{i, j}.\alpha,
\]
for $i, j = 1, 2$.
Note that being nef is equivalent to require all of 
\[
a -k_{i, 1} - k_{i, 2}, a' - k_{1, j} - k_{2, j}, k_{i, j} \qquad (i, j = 1, 2)
\]
are non-negative.

We define $U_{\mathbf{k},*}$ and $Z_{a, \mathbf{k}, *}^{\circ}$ as in \cref{sect:pointedmorzak}.  Using the discussions in \cref{sec:curvesandconicbundle}, one can prove the following propositions:

\begin{prop}
Assume that $2a > \sum k_{i, j}$ and $2a' > \sum k_{i, j}$.  Then the composition
\[
 \phi_{a, \mathbf{k}, *} \circ\psi_{\alpha, *}: M_{\alpha, *} \to Z_{a, \mathbf{k}, *}^{\circ} \to U_{\mathbf{k}, *}.
\]
is smooth, and its restriction over $\phi_{a, \mathbf{k}, *}\circ\psi_{\alpha, *} (M_{\alpha, *} )$ is a Zariski open subset of a bundle of the product $\mathbb P^{2a -\sum k_{i, j}} \times \mathbb P^{2a' - \sum k_{i, j}}$ over $\phi_{a, \mathbf{k}, *}(\psi_{\alpha, *} (M_{\alpha, *} ))$.
\end{prop}

\begin{prop}
\label{prop:dimensionestimatepointed_toric}
Assume that $2a > \sum k_{i, j}$ and $2a' > \sum k_{i, j}$.  Consider the composition
\[
 \phi_{a, \mathbf{k}, *}\circ\psi_{\alpha, *}  : M_{\alpha, *} \to Z_{a, \mathbf{k}, *}^{\circ} \to U_{\mathbf{k}, *}.
\]
Then the dimension of $U_{\mathbf{k}, *} \setminus  \phi_{a, \mathbf{k}, *}\circ\psi_{\alpha, *} (M_{\alpha, *})$ is at most 
\[
\max\left\{\sum_{i, j = 1, 2} 2k_{i, j} - 2a, \sum_{i, j = 1, 2} 2k_{i, j} - 2a' \right\}.
\]
\end{prop}

\subsection{Homological stability and the space of positive continuous maps}

\subsubsection{Semi-topological models and the spaces of positive maps.}

Here we follow the discussion of \cref{sec:CJS}. We work over $\mathbb C$. For $a, a', \mathbf{k}$ as before, we define the space of rational curves $M_{a, a', \mathbf{k}}$ and its torsor $\widetilde{M}_{a, a', \mathbf{k}}$. One can also define their semi-topological models $\mathcal M_{a, a', \mathbf{k}}$ and $\widetilde{\mathcal M}_{a, a', \mathbf{k}}$ as \cref{sec:CJS}. Then \cref{theo:algebraicvssemitopological,theo:semitopologicalvstopological} are valid in this settings of the toric surface $Y$ and their proofs are exactly the same. We do not repeat these in this section.

\subsubsection{A proof of \cref{theo:torichomequiv}}

Let $a, a', \mathbf{k}$ be non-negative integers such that
\[
a -k_{i, 1} - k_{i, 2} > 0, \, a' - k_{1, i} - k_{2, i} > 0, k_{i, j} > 0,
\]
for $i, j = 1, 2$.
We define $I \in \mathbb N$ by
\[
I=   \frac{1}{8} \min\left\{2a -\sum_{i, j}k_{i, j},  2a' -\sum_{i, j}k_{i, j}\right\} - \frac{1}{2}.
\]
Define
\[
E(w < x) := 4m_0(x)  +3 \sum_{i, j} m_{\ell_{i, j}}(x) + 2\sum_j m_{V_j}(x) +2\sum_{i, j} m_{\ell_{i, 1} \oplus \ell_{j, 2}}(x) - 2\sum_{i, j} k_{i, j}
\]
and let $P \subset (U_{\mathbf{k}, *} < \overline{R})$ be the subposet consisting of $w < x$ such that $E(w < x) \leq 2I$.
Repeating the argument for \cref{lemm:Eproper}, we see that $P$ is proper over $U_{\mathbf{k}}$:
\begin{lemm}
\label{lemm:Eproper_toric}
For any positive constant $\mathfrak{T}$ the set of $x \in \mathrm{Hilb}(B)^Q$ such that $4m_0(x)  +3 \sum_{i, j} m_{\ell_{i, j}}(x) + 2\sum_j m_{V_j}(x) + 2\sum_{i, j} m_{\ell_{i, 1} \oplus \ell_{j, 2}}(x)\ \leq \mathfrak{T}$ is closed and downward closed.
\end{lemm}

Let $M_{a, a', \mathbf{k}, *}$ be the space of pointed rational curves of class $(a, a', \mathbf{k})$ on $Y$.
Then we have a natural projection
\[
M_{a, a', \mathbf{k}, *} \to U_{\mathbf{k}, *}
\]
Let $F_1 \subset U_{\mathbf{k}}$ be the complement of the image of the above projection.
It follows from \cref{prop:dimensionestimatepointed_toric} that the real codimension of $F_1$ is greater than $16I$.

Let $(w < x)$ be an element of $P$ and $w < y$ be a essential pair such that $x \prec y$.
Let us consider the space
\[
\Gamma(B, \mathcal O(a)\otimes V_1 \oplus \mathcal O(a')\otimes V_2)_y.
\]
We define
\begin{align*}
a_y & = a- m_0(y) - m_{V_2}(y) - \sum_{i = 1}^2m_{\ell_{i, 2}}(y) \\
 a'_y &= a- m_0(y) - m_{V_1}(y) - \sum_{i = 1}^2m_{\ell_{i, 1}}(y) \\
 k_{y, i, 1} & = m_{\ell_{i, 1}}(y) + m_{\ell_{i, 1} \oplus \ell_{1, 2}}(y) +m_{\ell_{i, 1} \oplus \ell_{2, 2}}(y)\\
  k_{y, i, 2} & = m_{\ell_{i, 2}}(y) + m_{\ell_{1, 1} \oplus \ell_{i, 2}}(y) +m_{\ell_{2, 1} \oplus \ell_{i, 2}}(y)\\
y_j(\ell_{i, 1}) & =  \left(\sum_{c \in B} (m_{\ell_{i, 1}}(f_{y, c}) + m_{\ell_{i, 1} \oplus \ell_{1, 2}}(f_{y, c}) + m_{\ell_{i, 1} \oplus \ell_{2, 2}}(f_{y, c}))[c] \right)_i\\
y_j(\ell_{i, 2}) & =  \left(\sum_{c \in B} (m_{\ell_{i, 2}}(f_{y, c}) + m_{\ell_{1, 1} \oplus \ell_{i, 2}}(f_{y, c}) + m_{\ell_{2, 1} \oplus \ell_{i, 2}}(f_{y, c}))[c] \right)_i.
\end{align*}
Then the above space can be identified with
\[
\Gamma(B, \mathcal O(a_y)\otimes V_1)_{y_1}\oplus \Gamma(B, \mathcal O(a'_y)\otimes V_2)_{y_2}.
\]
We consider 
\[
E_j(w < y) = 2m_0(y) + 2m_{V_{j'}}(y) +2 \sum_{i = 1}^2 m_{\ell_{i, j'}}(y) + \sum_{i = 1}^2 k_{y, i, j} - \sum_{i, j= 1}^2 k_{i, j}
\]
where $\{j, j'\} = \{1, 2\}$. Then since this is monotone increasing with respect to $y$ we have $E_1(w < y) \geq 0$ and $E_2(w < y)\geq 0$. Moreover we have $E = E_1 + E_2$. Thus we have $E_j \leq 2I + 2$. This implies
\[
2a_y - \sum_{i = 1}^2 k_{y, i, 1}  \geq 6I+2, \qquad
2a'_y - \sum_{i = 1}^2 k_{y, i, 2} \geq 6I+2.
\]
The space $\Gamma(B, \mathcal O(a_y)\otimes V_1)_{y_1}$ is unobstructed if $y_1$ is contained in the image of 
\[
Z_{a_y, \mathbf{k}_{y, 1}}^\circ \to U_{\mathbf{k}_{y,1}}.
\]
Then the complement of the image has dimension at most $\sum_{i = 1}^2k_{y_1, i, 1} -6I-2$. Thus the dimension of the set $w$ such that $w < y$ where $y$ corresponds to a point $y_{1}$ in the complement of the image of $Z_{a_y, \mathbf{k}_{y, 1}}^\circ \to U_{\mathbf{k}_{y,1}}$ is at most
\[
m_0(y) + m_{V_1}(y) + m_{V_2}(y) +  \sum_{i = 1}^2m_{\ell_{i, 2}}(y) + \sum_{i = 1}^2k_{y, i, 1} -6I-2 \leq \sum k_i - 4I.
\]
Let $F_2$ be the union of closures of such loci while $T$ runs over all combinatorial types we consider. Note that the real codimension of $F_2$ in $U_{\mathbf{k}}$ is greater than or equal to $8I$.
Similarly we construct $F_3$ for the space $\Gamma(B, \mathcal O(a'_y)\otimes V_2)_{y_2}$.
Let $F = F_1 \cup F_2 \cup F_3$. 

Next we analyze $\kappa(T)$ in our situation.
It follows from our definition that we have
\[
2\gamma(T) = 8m_0(T) + 6 \sum_{i, j} m_{\ell_{i, j}(T)} + 4\sum_j m_{V_j}(T) + 4\sum_{i, j} m_{\ell_{i, 1} \oplus \ell_{j, 2}}(T), 
\]
and
\[
\mathrm{rank}(T) = 3m_0(T) + 2 \sum_{i, j} m_{\ell_{i, j}}(T) + \sum_j m_{V_j}(T) + \sum_{i, j} m_{\ell_{i, 1} \oplus \ell_{j, 2}}(T).
\]
Thus we have
\[
\kappa(T) = 5m_0(T) + 4 \sum_{i, j} m_{\ell_{i, j}}(T) + 3\sum_j m_{V_j}(T) + 3\sum_{i, j} m_{\ell_{i, 1} \oplus \ell_{j, 2}}(T) - 2|\mathrm{Supp}(T)|.
\]
Using this one can prove that when we have $x_1 \prec x_2$, we have $\kappa(x_2) \geq \kappa(x_1) + 1$. This shows that $\kappa(T) \geq \mathrm{rank}(T)$. 
Thus we have $E(T) \leq 2\kappa (T)$. Thus by applying versions of \cref{theo:algebraicvssemitopological,theo:semitopologicalvstopological} for $Y$ we conclude the following theorem:
\begin{theo}
\label{theo:maintoric}
Let $B = \mathbb P^1$ and $Y$ be the toric surface of degree $4$.
We fix an integral nef class $\alpha$ on $S$.
Suppose that we have a birational morphism $\beta : S \to \mathbb P^1 \times \mathbb P^1$ satisfying
\[
2a -\sum_{i, j}k_{i, j} > 0, \, 2a' -\sum_{i, j}k_{i, j} > 0.
\]
We define $I \in \mathbb N$ by
\[
I = \frac{1}{8} \min\left\{2a -\sum_{i, j}k_{i, j},  2a' -\sum_{i, j}k_{i, j}\right\} -\frac{1}{2}.
\]
Then $M_{\alpha, *} \hookrightarrow \mathrm{Top}^+_*(B, S)_{\alpha}$ is homology $\lfloor I \rfloor$-connected.
\end{theo}

Here we record one proposition:

\begin{prop}[{\cite[Proof of Corollary 1.2]{DT24}}] \label{prop:alphainvariance}
Let $X$ be a smooth projective variety, let $\{ E_{i} \}$ be a finite set of disjoint smooth divisors on $X$, and let $\eta: S^{2} \to X$ denote a pointed curve such that $\eta(S^{2}) \cap E_{i} = \emptyset$ for every $i$ and $\eta$ sends at least two points to the basepoint of $X$.

Let $s: C \to X$ be a pointed continuous map.  Suppose we fix a disc around the basepoint $*_{C}$ and let $\pi: C \to S^{2} \vee C$ denote the contraction of the boundary of the disc.  Then the stabilization morphism $s \mapsto (\eta \vee s) \circ \pi$ defines homotopy equivalences
\begin{align*}
\Top^{+}_{*}(C,X)_{\alpha} & \xrightarrow{\eta \vee} \Top^{+}_{*}(C,X)_{\alpha + \eta_{*}\mathbb{P}^{1}} \\
\Top_{*}(C,X)_{\alpha} & \xrightarrow{\eta \vee} \Top_{*}(C,X)_{\alpha + \eta_{*}\mathbb{P}^{1}}
\end{align*}
\end{prop}

\begin{proof}
We can think of $\eta$ as a homotopy from the constant map $S^{1} \to *_{X}$ to itself.  Let $\bar{\eta}$ denote the reversed homotopy, also considered as a map $\bar{\eta}: S^{2} \to X  - (\cup_{i} E_{i})$.  Then the map $\vee \bar{\eta}$ defines a homotopy inverse to $\vee \eta$, since the map
\begin{equation*}
S^{2}  \to S^{2} \vee S^{2} \xrightarrow{\eta \vee \bar{\eta}} X - (\cup_{i} E_{i})
\end{equation*}
is null-homotopic.
\end{proof}

Now we prove our main theorem of this section.

\begin{proof}[Proof of \cref{theo:torichomequiv}]
Let $J = \min\{ k_{i, j}, a -k_{i, 1} - k_{i, 2}, a' - k_{1, j} - k_{2, j}\}$.
By \cite[Theorem 4.1]{Guest}, for any $i \leq J$, we have an isomorphism
\[
H_i^{\mathrm{sing}}(M_{\alpha, *}, \mathbb Z) \cong H_i^{\mathrm{sing}}(\mathrm{Top}_*(B, S)_{\alpha})
\]
which is induced by the inclusion $M_{\alpha, *} \hookrightarrow \mathrm{Top}_*(B, S)_{\alpha}$.
Thus if we let
\[
\ell(\alpha) = \min\{\lfloor I \rfloor -1, J \},
\]
then the inclusion $\mathrm{Top}^+_*(B, S)_{\alpha} \hookrightarrow \mathrm{Top}_*(B, S)_{\alpha}$
is homology isomorphic up to degree $\ell(\alpha)$.
However, by \cref{prop:alphainvariance}, the homotopy types of $\mathrm{Top}^+_*(B, S)_{\alpha}$ and $\mathrm{Top}_*(B, S)_{\alpha}$ do not depend on $a, a'$. If we fix $\{ k_{i,j} \}$ and make $a,a'$ sufficiently large then $\ell(\alpha) = \min\{ k_{i,j} \}$ and our assertion follows.
\end{proof}

\printbibliography

\clearpage


\begin{refsection}

\section{Betti number bounds for the moduli space of maps from a rational curve}
 \label{appendixSS}
 
\smallskip
\begin{center} \textsc{Will Sawin, \ Mark Shusterman\footnote{The Dr. A. Edward Friedmann Career Development Chair in Mathematics}} \end{center}
\medskip

Let $k$ be a separably closed field. Let $X$ be a projective variety over $k$. Fix a very ample line bundle $L$ on $X$. Say a map $g\colon \mathbb P^1 \to X$ has degree $e$ if the pullback of $L$ to $\mathbb P^1$ by $g$ has degree $e$. Let $\Mor(\mathbb P^1,X)$ be the moduli space parameterizing morphisms of degree $e$ from $\mathbb P^1$ to $X$. The goal of this appendix is to prove the following theorem.

\begin{theo}\label{appendix-main} There exists a constant $C$ such that, for $\ell$ a prime invertible in $k$, we have
\[ \sum_{i=0}^\infty \dim H^i_c ( \Mor(\mathbb P^1, X), \mathbb Q_\ell) \leq C^{e+1}.\]
The constant $C$ depends only on $\dim H^0(X, L)$ and the self-intersection $c_1(L)^{\dim X}$. \end{theo}

To prove this, we will use Betti number bounds of Katz for the compactly supported cohomology of intersections of hypersurfaces in affine space. We are faced with the difficulty that $\Mor(\mathbb P^1, X)$ is not usually an intersection of hypersurfaces in affine space. We will replace $\Mor(\mathbb P^1, X)$ by an intersection of hypersurfaces in affine space in two steps: First defining an affine space bundle $\tMor(\mathbb P^1, X)$ on $\Mor(\mathbb P^1, X)$ and then a stratification of $\tMor(\mathbb P^1, X)$ into strata $\tMor^j(\mathbb P^1, X)$, which will be intersections of hypersurfaces in affine space. It is possible to control the way Betti numbers change under these operations, leading to a Betti number bound for $\Mor(\mathbb P^1, X)$.

 First, observe that $L$ defines an embedding $X \to \mathbb P^n$ for the natural number 
 \[
 n = \dim H^0(X,L)-1.
 \]
 
\begin{lemm}\label{bounded-polynomial-degree} Inside $\mathbb P^n$, $X$ is the vanishing locus of some homogeneous polynomials $f_1,\dotsc,f_r$ of degrees $d_1,\dotsc,d_r$ in $n+1$ variables, with $d_1\leq \dotsb \leq d_r$, such that $r$ and $d_r$ are bounded in terms of $\dim H^0(X,L)$ and $c_1(L)^{\dim X}$. 
\end{lemm}
\begin{proof}Inside $\mathbb P^n$,  $X$ is a closed subscheme of degree $c_1(L)^{\dim X}$. Being a closed subscheme, it is the vanishing locus of some homogeneous polynomials. By changing the order of these polynomials, we may assume their degrees are nondecreasing. Using the Chow variety, one sees that all smooth closed $X \subseteq \mathbb P^n$ of degree $c_1(L)^{\dim X}$ are parameterized by a scheme of finite type. One can see that the number $r$ of polynomials defining $X$ and the degrees $d_1,\dotsc, d_r$ of these polynomials are bounded only in terms of $\dim H^0(X, L)$ and the self-intersection $c_1(L)^{\dim X}$ by a Noetherian induction argument over this scheme.
\end{proof}

By Lemma \ref{bounded-polynomial-degree}, it suffices to prove Theorem \ref{appendix-main} where the constant $C$ depends only on $n$ and $d_r$.

We view $\mathbb P^{(n+1)(e+1)-1}$ as the space parameterizing $n+1$-tuples of polynomials $g_0,\dotsc, g_n$ of degree $\leq e$ in the variable $T$, not all $0$, up to scalar multiplication of all polynomials. Here the coordinates of $\mathbb P^{(n+1)(e+1)-1}$ are the coefficients of the polynomials.
\begin{lemm}\label{mor-coordinates} $\Mor(\mathbb P^1, X)$ is isomorphic to the locally closed subscheme of $\mathbb P^{(n+1)(e+1)-1}$ consisting of tuples of polynomials $g_0,\dotsc,g_n$ such that  $f_i(g_0,\dotsc,g_n)=0$ for $i$ from $1$ to $r$, at least one of $g_0,\dotsc, g_n$ has degree $e$, and $g_0,\dotsc, g_n$ share no common factors.  \end{lemm}
\begin{proof} We first construct the map $\Mor(\mathbb P^1,X) \to \mathbb P^{(n+1)(e+1)-1}$.  For  $g \colon \mathbb P^1 \to X$ of degree $e$, the pullback $g^* L$ is necessarily isomorphic to $\mathcal O_{\mathbb P^1}(e)$.  If we fix an isomorphism $g^* L \cong \mathcal O_{\mathbb P^1}(e)$ we obtain a linear map $H^0( X, L) \to H^0( \mathbb P^1, \mathcal O_{\mathbb P^1}(e))$.

Choosing a basis for $H^0( X, L)$, we obtain $n+1$ elements of $H^0( \mathbb P^1, \mathcal O_{\mathbb P^1}(e))$, which cannot all be zero. We can express $H^0( \mathbb P^1, \mathcal O_{\mathbb P^1}(e))$ as the space of polynomials of degree $\leq e$ in one variable $T$ (for example, by first expressing $H^0( \mathbb P^1, \mathcal O_{\mathbb P^1}(e))$  as the space of homogeneous polynomials of degree $e$ in two variables $u,v$ and then sending $u$ to $1$ and $v$ to $T$). This gives a tuple $g_0,\dotsc, g_n$ of polynomials of degree $\leq e$. Varying the isomorphism $g^* L \cong \mathcal O_{\mathbb P^1}(e)$ acts by scalar multiplication on this tuple of polynomials, so the resulting point of $\mathbb P^{(n+1)(e+1)-1}$ is well-defined.

Since the tuple $g_0,\dotsc,g_n$ determines the map $g$, this map is an immersion.



Finally, the image of $\Mor(\mathbb P^1, X) $ inside $ \mathbb P^{(n+1) (e+1)-1}$ consists of tuples of polynomials such that $f_i(g_0,\dotsc,g_n)=0$ for $i$ from $1$ to $r$, at least one of $g_0,\dotsc, g_n$ has degree $e$, and $g_0,\dotsc, g_n$ share no common factors. Indeed the conditions that at least one of $g_0,\dotsc,g_n$ has degree $e$ and $g_0,\dotsc,g_n$ share no common factors are necessary and sufficient for $g_0,\dotsc, g_n$ to define a map $\mathbb P^1 \to \mathbb P^n$ of degree $e$, and this map has image in $X$ if and only if $g_0,\dotsc,g_n$ satisfy the equations of definition $f_1,\dotsc, f_r$ of $X$. \end{proof}

We now define a space $\tMor(\mathbb P^1,X)$, which will be a subscheme of 
\begin{equation} \label{QuotientScheme}
\left ( \left(\mathbb A^{(n+1)(e+1)}- \{0\} \right) \times \mathbb A^{(n+1) e} \right)/\mathbb G_m. 
\end{equation}
We think of $\mathbb A^{(n+1)(e+1)}$ as parameterizing tuples  $g_0,\dotsc, g_n$ of polynomials of degree $\leq e$ in $T$. We think of $\mathbb A^{(n+1)e}$ as parameterizing tuples $h_0,\dotsc, h_n$ of polynomials of degree $\leq e-1$ in $T$. We let $\mathbb G_m$ act by multiplication on the coordinates of $\mathbb A^{(n+1)(e+1)}$ and inverse multiplication on the coordinates of $\mathbb A^{(n+1)e}$, in other words by scalar multiplication on $g_0,\dotsc, g_n$ and inverse scalar multiplication on $h_0,\dotsc,h_n$.

We define $\tMor(\mathbb P^1,X)$ as the locally closed subscheme of \ref{QuotientScheme} defined by the equations $f_i( g_0,\dotsc,g_n)=0$ for $i=1,\dotsc, r$, the equation $\sum_{j=0}^n g_i h_i=1$, and the condition that at least one of the $g_i$ has degree exactly $e$. Observe that the first two conditions are closed and the last condition is open.  The fact that the $f_i$ are homogeneous means that the first condition is invariant under the $\mathbb G_m$ action, and the fact that $\sum_{j=0}^n g_j h_j=1$ has each $g_j$ appearing multiplied by an $h_j$ means that the second condition is invariant under the $\mathbb G_m$ action.

Note that the conditions defining $\tMor(\mathbb P^1,X)$ already imply that $g_0,\dotsc,g_n$ are not all zero. The only reason we use $\mathbb A^{(n+1)(e+1)}- \{0\} $ instead of $\mathbb A^{(n+1)(e+1)}$ in the definition is to express $\tMor(\mathbb P^1,X)$ as a subscheme of a scheme and not a substack, that happens to also be a scheme, of a stack.

\begin{lemm}\label{map-exists} 
There exists a map $\tMor(\mathbb P^1, X)\to \Mor(\mathbb P^1, X)$ given by 
\[
((g_0,\dotsc,g_n),(h_0,\dotsc,h_n)) \mapsto (g_0,\dotsc,g_n).
\] 
\end{lemm}
\begin{proof}It suffices to check that if $g_0,\dotsc,g_n, h_0,\dotsc,h_n$ satisfy the conditions of $\tMor(\mathbb P^1, X)$ then $g_0,\dotsc,g_n$ satisfy the conditions of $ \Mor(\mathbb P^1, X)$. The only condition for which this is not immediate is the condition that $g_0,\dotsc,g_n$ have no common factors, but the  condition  $\sum_{j=0}^n g_j h_j=1$ ensures that $g_0,\dotsc,g_n$ have no common factors. \end{proof}

We say a map $Z \to Y$ is an affine space bundle of relative dimension $d$ if it is locally in the Zariski topology on $Y$ isomorphic to $\mathbb A^d\times Z \to Z$, with transition maps affine transformations.

We will check that the map of Lemma \ref{map-exists} is an affine space bundle of relative dimension $(n-1)e$, using a series of lemmas. We first explain the utility of affine space bundles to us, by the following lemma.

\begin{lemm}\label{affine-bundle-cohomology} Let $\pi \colon Z \to Y$ be an affine space bundle of relative dimension $d$. Then we have $H^i_c ( Y, \mathbb Q_\ell) = H^{i+2d}_c(Z,\mathbb Q_\ell(d))$ for all $i$. \end{lemm}
\begin{proof} We have a natural trace map $R \pi_! \mathbb Q_\ell \to \mathbb Q_\ell [ -2d](d)$. We can check that this map is an isomorphism locally, which reduces to the case that $Z = \mathbb A^d \times Y$, where proper base change reduces us to the case that $Y$ is a point, which is just the computation of compactly supported cohomology of affine space. This isomorphism together with the Leray spectral sequence with compact supports gives $H^i_c ( Y, \mathbb Q_\ell) = H^{i+2d}_c(Z,\mathbb Q_\ell(d))$. \end{proof}

Let $\hMor(\mathbb P^1,X)$ be the subset of \ref{QuotientScheme} consisting of points satisfying the conditions that $f_i(g_0,\dotsc,g_n)=0$ for $i$ from $1$ to $r$, at least one of $g_0,\dotsc, g_n$ has degree $e$, and $g_0,\dotsc, g_n$ share no common factors.

\begin{lemm}\label{easier-map-bundle} The map $\hMor(\mathbb P^1,X)\to \Mor(\mathbb P^1,X)$ sending $((g_0,\dotsc,g_n),(h_0,\dotsc,h_n))$ to $(g_0,\dotsc,g_n)$ is an affine space bundle of relative dimension $(n+1)e$. \end{lemm}

\begin{proof} Over the intersection of $\Mor(\mathbb P^1,X)$ with any standard affine open in $\mathbb P^{(n+1) (e+1)-1}$, the space $\hMor(\mathbb P^1,X)$ becomes isomorphic to $\mathbb A^{(n+1)e}  \times \Mor(\mathbb P^1,X)$ via the coordinates $h_0,\dotsc,h_n$, and the change of coordinates is given by scalar multiplication which is an affine transformation. \end{proof}

\begin{lemm}\label{inverse-image-affine} Let $Y$ be a scheme, $Z$ an affine space bundle on $Y$ of relative dimension $a$, and $Z'$ an affine space bundle on $Y$ of relative dimension $b$. Let $q \colon Z\to Z'$ a map over $Y$ and $s$ a section of $Z'$. If the map $q$ is fiberwise surjective then  $q^{-1}(s)$ is an affine space bundle on $Y$ of relative dimension $a-b$.\end{lemm}

\begin{proof} Working locally, we may assume $Z= \mathbb A^a \times Y$ and $Z'= \mathbb A^b\times Y$, so the map $q$ is given by an $ a\times b$ matrix together with an additional $1\times b$ vector giving the translation. Fiberwise surjectivity at a point is witnessed by the nonvanishing of one of the $b\times b$ minors of the matrix. Over an open neighborhood of that point where the same $b\times b$ minor is nonvanishing, we can use the $a-b$ coordinates not appearing in the minor as coordinates for $q^{-1}(s)$, giving an isomorphism of $q^{-1}(s)$ to affine space. The remaining coordinates are given by affine functions of these coordinates, showing that the transition maps are affine.\end{proof}

\begin{lemm}\label{map-is-bundle} The map $\tMor(\mathbb P^1, X)\to \Mor(\mathbb P^1, X)$ of Lemma \ref{map-exists} is an affine space bundle of relative dimension $(n-1)e$. \end{lemm}

\begin{proof} View $\mathbb A^{2e}$ as the space parameterizing polynomials of degree $\leq 2e-1$. There is a map 
\[
q \colon \hMor(\mathbb P^1 ,X) \to  \mathbb A^{2e} \times \Mor(\mathbb P^1,X)
\]
given by
\[ 
q((g_0,\dotsc,g_n),(h_0,\dotsc,h_n)) = \left(\sum_{j=0}^n g_j h_j, (g_0,\dotsc,g_n)\right).  
\]

Since $\tMor(\mathbb P^1,X)$ is by definition the inverse image of the section $1$ under this map, by Lemmas \ref{easier-map-bundle} and \ref{inverse-image-affine} it suffices to check that $q$ is fiberwise surjective.

For this purpose, fix a tuple of polynomials $g_0,\dotsc,g_n$, one of degree $e$, with no common factors.  

Choose $j_0$ such that $g_{j_0}$ has degree exactly $e$. Then all polynomials of degree $\leq 2e - 1$ that are multiples of $g_{j_0}$ have the form $g_{j_0} h$ for $h$ of degree $\leq e-1$. Since $g_0,\dotsc, g_n$ have no common factors, every polynomial may be written as $\sum_{j=0}^n g_j h_j$ for some $h_j$, not necessarily of degree $\leq e-1$. The image of $\sum_{j=0}^n g_j h_j$ in $k[T]/g_{j_0}$ is determined only by $h_j$ mod $g_{j_0}$. We may replace $h_j$ by a polynomial of degree $<e$ with the same residue mod $g_{j_0}$ by long division with remainder. It follows that every polynomial is congruent mod $g_{j_0}$ to a polynomial of the form  $\sum_{j=0}^n g_j h_j$ for some $h_j$ of degree $\leq e-1$. If the original polynomial has degree $2e-1$, the difference between that polynomial and $\sum_{j=0}^n g_j h_j$ is a multiple of $g_{j_0}$ of degree $\leq 2e-1$, hence of the form $g_{j_0} h$ for some $h$ of degree $\leq e-1$, and subtracting $h$ from $h_{j_0}$ gives the needed $h_0,\dotsc, h_n$ to check surjectivity. \end{proof}

Using Lemmas \ref{affine-bundle-cohomology} and \ref{map-is-bundle}, we will reduce to bounding the cohomology  of $\tMor(\mathbb P^1,X)$.

To do this, we will further stratify $\tMor(\mathbb P^1,X)$.

Let $\tMor^j(\mathbb P^1,X)$ be the subset of $\tMor(\mathbb P^1,X)$ consisting of tuples $(g_0,\dotsc,g_n),(h_0,\dotsc,h_n)$ where $g_0,\dotsc, g_{j-1}$ have degree $<e$ and $g_j$ has degree exactly $e$. The union of $\tMor^j(\mathbb P^1,X)$ for $j$ from $m$ to $n+1$ is a closed subset of $\tMor(\mathbb P^1,X)$ for each $m$, so this indeed gives a stratification.

\begin{lemm}\label{stratification-bound}
We have
\[ \sum_{i=0}^\infty \dim H^i_c ( \tMor(\mathbb P^1,X), \mathbb Q_\ell)  \leq  \sum_{j=0}^n \sum_{i=0}^\infty \dim H^i_c ( \tMor^j(\mathbb P^1,X), \mathbb Q_\ell). \] \end{lemm}
\begin{proof} This follows immediately from the spectral sequence associated to the stratification of $\tMor(\mathbb P^1,X)$ into $\tMor^j(\mathbb P^1,X)$. In other words, this is the spectral sequence associated to the filtration of  $H^i_c ( \tMor(\mathbb P^1,X), \mathbb Q_\ell) $ into $H^i_c ( \bigcup_{j=m}^n \tMor^j (\mathbb P^1,X), \mathbb Q_\ell) $ whose associated graded is calculated using the excision long exact sequence. \end{proof}

\begin{lemm}\label{is-complete-intersection} $\tMor^j(\mathbb P^1,X)$ is isomorphic to the vanishing locus in $\mathbb A^{ (n+1)(2e+1)-1-j}$ of 
\[
2e+ \sum_{i=1}^r (d_i e+1)
\]
equations of degree $\leq \max(d_r,2)$. \end{lemm}

\begin{proof} Each point of $\tMor^j(\mathbb P^1,X)$ corresponds to an orbit of tuples $g_0,\dotsc, g_n, h_0,\dotsc,h_n$ under the $\mathbb G_m$-action, but since by definition of $\tMor^j(\mathbb P^1,X)$ the coefficient of $T^e$ in $g_j$ is nonzero, this orbit has a unique member where the coefficient of $T^e$ in $g_j$ is $1$. We can thus express a point of $\tMor^j(\mathbb P^1,X)$ using coordinates the coefficients of the polynomials $g_0,\dotsc, g_n, h_0,\dotsc, h_n$, except that we do not need to include the coefficients of $T^e$ in $g_0,\dotsc, g_{j-1}$ (as these are always $0$) or the coefficient of $T^e$ in $g_j$ (as this is always $1$). This gives an embedding of $\tMor^j(\mathbb P^1,X)$ into $\mathbb A^{ (n+1) (e+1)-1-j + (n+1)e}$. Inside $\mathbb A^{ (n+1) (e+1)-1-j + (n+1)e}$, $\tMor^j(\mathbb P^1,X)$ is defined by the conditions that $f_i(g_0,\dotsc,g_n)=0$ for $i$ from $1$ to $r$ and that $\sum_{j=0}^n g_j h_j=1$. These are both closed conditions.

More precisely, the condition that $f_i(g_0,\dotsc, g_n)=0$ may be expressed as $d_i e +1 $ equations in the coefficients of $g_0,\dotsc,g_n$ of degree $\leq d_i$, since $f_i(g_0,\dotsc,g_n)$ is a polynomial of degree $\leq d_i e$ and each of its coefficients is a polynomial of degree $\leq d_i$ in the coefficients of $g_0,\dotsc,g_n$ which must be set to $0$. Similarly, the condition that $\sum_{j=0}^n g_j h_j=1$ may be expressed as $2e$ equations in the coefficients of $g_0,\dotsc, g_n, h_0,\dotsc,h_n$ of degree $\leq 2$, since $\sum_{j=0}^n g_j h_j$ is a polynomial of degree $\leq 2e-1$ and each of its coefficients is a polynomial of degree $\leq 2$ in the coefficients of $g_0,\dotsc, g_n, h_0,\dotsc,h_n$ which must be set to either $0$ or $1$.

Combining these, $\tMor^j(\mathbb P^1,X)$ is the intersection in $\mathbb A^{ (n+1) (e+1)-1-j + (n+1)e} = \mathbb A^{ (n+1)(2e+1)-1-j}$ of  $2e+ \sum_{i=1}^r (d_i e+1)$ equations of degree $\leq \max(d_r,2)$.\end{proof}

\begin{lemm}\label{fundamental-betti-bound}We have \[ \sum_{i=0}^\infty \dim H^i_c ( \tMor^j(\mathbb P^1,X), \mathbb Q_\ell)  \leq 3 (\max(d_r,2)+2)^{ (n+1)(2e+1)-1-j  + 2e+ \sum_{i=1}^r (d_i e+1)}.\] 
\end{lemm}

\begin{proof} This follows from \cite[Theorem 12]{katzbetti}, taking the polynomial $f$ of Katz to vanish, the polynomials $F_1,\dotsc,F_r$ to be the equations defining $\tMor^j(\mathbb P^1,X)$ described by Lemma \ref{is-complete-intersection}, and $s=0$ so that there are no polynomials $G_1,\dotsc,G_s$.  \end{proof}

\begin{prop}\label{TildeBound} We have
\[ \sum_{i=0}^\infty \dim H^i_c( \tMor(\mathbb P^1,X), \mathbb Q_\ell) \leq 4 (\max(d_r,2)+2)^{n +r}  \left( (\max(d_r,2)+2)^{ 2n+4 + \sum_{i=1}^r d_i} \right)^e  .\]
\end{prop}

\begin{proof}By Lemmas \ref{stratification-bound}, and \ref{fundamental-betti-bound}, we have
\begin{equation*}
\begin{split}
\sum_{i=0}^\infty \dim H^i_c( \tMor(\mathbb P^1,X), \mathbb Q_\ell) &\leq   \sum_{j=0}^n \sum_{i=0}^\infty \dim H^i_c ( \tMor^j(\mathbb P^1,X), \mathbb Q_\ell) \\
&\leq \sum_{j=0}^n  3 (\max(d_r,2)+2)^{ (n+1)(2e+1)-1-j  + 2e+ \sum_{i=1}^r (d_i e+1)} \\
&\leq  \sum_{j=0}^\infty 3 (\max(d_r,2)+2)^{ (n+1)(2e+1)-1-j  + 2e+ \sum_{i=1}^r (d_i e+1)} \\
&= \frac{3 }{ 1- (\max(d_r,2)+2)^{-1}} (\max(d_r,2)+2)^{ (n+1)(2e+1)-1  + 2e+ \sum_{i=1}^r (d_i e+1)} \\
&\leq 4 (\max(d_r,2)+2)^{ (n+1)(2e+1)-1  + 2e+ \sum_{i=1}^r (d_i e+1)} \\
&=  4 (\max(d_r,2)+2)^{n +r}  \left( (\max(d_r,2)+2)^{ 2n+4 + \sum_{i=1}^r d_i} \right)^e.
\end{split}
\end{equation*}
 \end{proof}

\begin{proof}[Proof of Theorem \ref{appendix-main}] 
This follows from Lemmas \ref{TildeBound}, \ref{affine-bundle-cohomology}, \ref{map-is-bundle}, and \ref{bounded-polynomial-degree} if we take  
\[
C = \max ( (\max(d_r,2)+2)^{ 2n+4 + \sum_{i=1}^r d_i} ,  4 (\max(d_r,2)+2)^{n +r} ). 
\]
In fact, since $d_1,\dotsc,d_r\geq 1$, the ratio 
\[
\frac{(\max(d_r,2)+2)^{ 2n+4 + \sum_{i=1}^r d_i}}{  4 (\max(d_r,2)+2)^{n +r}}
\]
 is at least 
 \[
 \frac{(\max(d_r,2)+2)^{n+4}}{ 4} \geq \frac{4^5}{ 4} > 1
 \]
 so $(\max(d_r,2)+2)^{ 2n+4 + \sum_{i=1}^r d_i} $ is larger and we can simply take \[C=(\max(d_r,2)+2)^{ 2n+4 + \sum_{i=1}^r d_i} .\qedhere\] \end{proof}

One can try to improve the value of $C$ using ideas from \cite{WanZhang}.

\bigskip

\

\noindent
{\bf Acknowledgments:}
While working on this research, Will Sawin was supported NSF grant DMS-2101491, and was a Sloan Research Fellow. Mark Shusterman's research is co-funded by the European Union (ERC, Function Fields, 101161909). Views and opinions expressed are however those of the authors only and do not necessarily reflect those of the European Union or the European Research Council. Neither the European Union nor the granting authority can be held responsible for them.

\printbibliography[title={References (for Appendix~\ref{appendixSS})}]

\end{refsection}

\end{document}